\documentclass[12pt,a4paper,english]{article}
\usepackage[T1]{fontenc}
\usepackage[applemac]{inputenc}
\usepackage{babel}
\usepackage{amsthm}
\usepackage{amsmath}
\usepackage{amssymb}
\usepackage{amsfonts}
\usepackage{mathpazo}
\usepackage{newtxtext}
\usepackage{hyperref}
\usepackage[all]{xy}
\usepackage{tikz}
\usetikzlibrary{matrix,arrows,decorations.pathmorphing,positioning}
\usepackage{geocft}

\numberwithin{equation}{section}

\begin{document}

\centerline{\LARGE{Geometric Class Field Theory}} 
\vspace{4pt} 
\centerline{\LARGE{with Bounded Ramification}} 
\vspace{25pt}
\centerline{\large{Henrik Russell}} 
\vspace{30pt}


\begin{abstract}
Let $U$ be a smooth quasi-projective variety 
over a field $k$ that is finite, the algebraic closure of a finite field 
or algebraically closed of characteristic 0. 
Let $X$ be a suitable projective compactification of $U$ 
such that the boundary is a divisor, 
and $\mdl$ an effective divisor on $X$ with support in $X \setminus U$. 
We consider a relative Chow group $\CHm{X}{\mdl}$ of modulus $\mdl$ 
(defined in a geometric way), 
the Albanese variety $\Albm{X}{\mdl}$ of $X$ of modulus $\mdl$ 
and 
the Abel-Jacobi map with modulus 
$\abjacm{X}{\mdl}: \CHmo{X}{\mdl} \lra \Albm{X}{\mdl}(k)$. 
We show that there is a 1-1 correspondence 
between relative Cartier divisors on $X$ 
and compatible systems of relative Cartier divisors on curves in $X$. 
This so called \emph{Skeleton Theorem} allows us to prove a 
\emph{Roitman Theorem with Modulus}: 
the \emph{Abel-Jacobi map with modulus} 
is an isomorphism on torsion parts 
over an algebraic closure of a finite field 
or an algebraically closed field of characteristic 0. 
We obtain a Reciprocity Law and an Existence Theorem 
for abelian coverings of $X$ over a finite field 
with ramification bounded by $\mdl$. 
All of this is done for a log as well as for a non-log version. 
\end{abstract}

\setcounter{tocdepth}{2}
\tableofcontents{} 

\setcounter{section}{-1}

\section{Introduction} 
\label{Intro}

Let $X$ be a projective variety over a field $\fld$, 
$\mdl$ an effective Cartier divisor on $X$ (with multiplicities), 
assume $X \setminus \mdl$ to be smooth. 

\bigskip 
\textbf{Geometric Class Field Theory:} 
Assume that the base field $\fld$ is finite. 
The motivating goal of this paper was to describe abelian coverings of $X$ 
whose ramification is bounded by $\mdl$ 
($\see$Definition \ref{bounded ramification})
in terms of a relative Chow group with modulus $\CHm{X}{\mdl}$ 
(defined in a geometric way, $\see$Definition \ref{DefCHoMod}). 
One of the main results 
is an isomorphism of finite groups 
\[ \CHmo{X}{\mdl} \iso \fundGgeom{X}{\mdl} 
\] 
($\see$Theorem \ref{reciprocityDiagram}, 
referred to as ``Reciprocity Law'' for historical reasons). 
Here $\fundGabm{X}{\mdl}$ is the {abelian fundamental group  
of $X$ of modulus $\mdl$} 
($\see$Definition \ref{fundGroupMod}). 
The upper $0$ denotes the degree $0$ part of $\CHm{X}{\mdl}$, 
respectively the geometric part of $\fundGabm{X}{\mdl}$ 
classifying those abelian coverings of $X$ that arise from a ``geometric situation'', 
i.e.\ not from extending the base field. 
If necessary, $X$ has to be blown up outside of $X\setminus\mdl$ 
to a suitable projective variety ($\see$Definition \ref{suitable_Def}, 
explained below). 
A similar result, but derived via different methods,  
is contained in the work of Kerz and Saito \cite{KeSa}. 

Main ingredients for the proof of the Reciprocity Law 
are Wiesend's (tame) class field theory \cite{KeSc09} 
and an affine version of the Roitman Theorem \cite{Roit} 
(see below). 

\bigskip 
\textbf{Roitman Theorem with Modulus:} 
The generalized Albanese variety with modulus $\Albm{X}{\mdl}$ from \cite{Ru13} 
is defined if $X$ is a suitable projective variety ($\see$Definition \ref{suitable_Def}) 
over a perfect base field, 
and $\Albm{X}{\mdl}$ admits a canonical map 
$\abjacm{X}{\mdl}: \CHmo{X}{\mdl} \lra \Albm{X}{\mdl}(k)$, 
the Abel-Jacobi map with modulus. 
If $X$ is smooth (a weaker assumption would be sufficient, 
see  $\see$Definition \ref{strongSuitable}) 
and the base field $\fld$ is an algebraic closure of a finite field 
or algebraically closed of characteristic $0$, 
the Abel-Jacobi map with modulus 
is an isomorphism on torsion parts: 
\[ \CHm{X}{\mdl}^{\tor} \iso \Albm{X}{\mdl}(\fld)^{\tor} 
\] 
($\see$Theorem \ref{RoitmanThm}). 
Let $\ACHm{X}{\mdl}$ denote the affine part of $\CHm{X}{\mdl}$ 
($\see$Definition \ref{Def_ACHm}), 
and  $\Lm{X}{\mdl}$ the affine part  (= largest affine subgroup) of $\Albm{X}{\mdl}$. 
For any suitable projective $X$, 
if the base field is an algebraic closure of a finite field 
or algebraically closed of characteristic $0$, 
or if $X \setminus \mdl$ is smooth and the base field is finite, 
the map on the affine part of $\CHm{X}{\mdl}$ 
induced by the Abel-Jacobi map with modulus 
is an isomorphism
\[ \abjacmaff{X}{\mdl}: \ACHm{X}{\mdl} \iso \Lm{X}{\mdl}(\fld) 
\] 
($\see$Theorem \ref{affineRoitmanThm}). 
For the proof of this ``affine Roitman Theorem'' 
we use duality and 

\bigskip 
\textbf{Skeleton Divisors:} 
Let $\Lm{X}{\mdl}$ be the affine part of the smooth connected 
algebraic group $\Albm{X}{\mdl}$ 
and $\Fmor{X}{\mdl}$ its Cartier dual. 
Here $\Fmor{X}{\mdl}$ is represented by a formal subgroup of 
the group-sheaf $\Divf_X$ of relative Cartier divisors on $X$ 
($\see$Definition \ref{DefFm(X,D)}). 
For a curve $C$ in $X$, 
the pull-back of relative Cartier divisors induces a map 
$\Fmor{X}{\mdl} \lra \Emor{C}{\mdl_{\Cgst}}$. 
Extending the notion of $\Fmor{C}{\mdl}$ in such a way 
that it is defined for singular curves $C$ in $X$ 
($\see$\Point \ref{dual_ACHm}), 
we can consider compatible systems of relative Cartier divisors on curves 
$\pn{\sD_C}_C \in \varprojlim_C \Emo{C}{\mdl_{\Cgst}}$ 
($\see$Definition \ref{Def_skeletonDiv}). 
Those will be called \emph{skeleton relative divisors} 
in analogy to the terminology \emph{2-skeleton sheaves} from \cite{EK}.
We show that the natural homomorphism 
from $\Fmor{X}{\mdl}$ to 
$\varprojlim_C \Emor{C}{\mdl_{\Cgst}}$, 
where $C$ ranges over the various curves in $X$, 
is an isomorphism 
\[ \skelm{X}{\mdl}: \Fmor{X}{\mdl} \iso \varprojlim_C \Emor{C}{\mdl_{\Cgst}} 
\] 
($\see$Theorem \ref{SkelThm}, referred to as ``Skeleton Theorem''), 
which establishes a 1-1 correspondence between relative Cartier divisors on $X$ 
and compatible systems of relative Cartier divisors 
on the various curves of $X$. 
This is the most substantial and most difficult part of the story, 
for a proof modulo technicalities see page \pageref{Idea}. 

\newpage 

Philosophically, this approach is related to 
a conjecture of Deligne 
on the existence of lisse $\ol{\Qrat}_{\ell}$-sheaves on $X$, 
$\seecite$\cite[Question 1.2]{EK}: 
A compatible system of lisse $\ol{\Qrat}_{\ell}$-sheaves on subcurves of 
$X \setminus \mdl$ 
whose ramification is bounded in terms of $\mdl$ 
comes from a lisse $\ol{\Qrat}_{\ell}$-sheaf on $X$. 
The link between relative Cartier divisors 
and lisse $\ol{\Qrat}_{\ell}$-sheaves 
is as follows: 
A lisse $\ol{\Qrat}_{\ell}$-sheaf on $X$ with ramification bounded by $\mdl$ 
is a continuous representation of the fundamental group with modulus 
$\fundG{X,\mdl}$ ($\see$Definition \ref{fundGroupMod}) 
on finite dimensional $\ol{\Qrat}_{\ell}$-vector spaces. 
(The group of lisse $\ol{\Qrat}_{\ell}$-sheaves of rank one 
is hence a dual of $\fundG{X,\mdl}$.)
On the other hand, 
if $R$ is a finite dimensional $k$-algebra, 
a relative Cartier divisor $\sD \in \Fmor{X}{\mdl}\pn{R}$ 
gives rise to a rational map $\rmp{\sD}$ 
of modulus $\leq \mdl$ ($\see$Definition \ref{DefMod}) 
from $X$ to an algebraic group $\Gp\pn{\sD}$, 
corresponding to the canonical section of the line bundle 
$\sO_{X \tens R}\pn{\sD}$ on $X \tens_k R$ 
($\see$\Point \ref{pairingFml}). 
Then $\rmp{\sD}$ factors through a homomorphism of torsors 
$\Albm{X}{\mdl} \lra \Gp\pn{\sD}$, 
thus $\sD$ can be thought of as a representation of 
the Albanese variety with modulus $\Albm{X}{\mdl}$. 
Due to the natural connection between $\Albm{X}{\mdl}$ and $\fundGgeo{X,\mdl}$ 
($\see$Corollary \ref{fundGm-Albm-exSeq}), 
the Skeleton Theorem becomes an analog 
of the rank one case of the conjecture of Deligne. 

Conversely, since Chow groups of $0$-cycles are compatible with curves 
by construction, 
any proof of the Reciprocity Law also covers the rank one case 
of Deligne's conjecture via duality
(cf.\ \cite[Theorem II]{KeSa}). 

\bigskip 
\textbf{Duality of Divisors and 0-Cycles:} 
The basic concept of this work is duality of group objects 
and the interplay of group objects associated to a variety 
and those associated to curves in that variety. 
For this aim, one needs to understand the structure of Chow groups of $0$-cycles with modulus, which is investigated in Section \ref{ChowMod}. 
A pairing between the affine part of these Chow groups and relative Cartier divisors is constructed in Section \ref{AbelJacobiMap}. 
The pairing is based on the local symbol from \cite{S59} on curves, 
and extended to higher dimensions via restriction to curves. 

This pairing becomes an essential tool for the proof of the Skeleton Theorem: 
in Subsection \ref{rigidity} it is used to formulate rigidity properties of relative Cartier divisors, 
i.e.\ an element $\sD \in \Fmor{C}{\mdl}$ 
is uniquely determined by the values on a fixed system of $0$-cycles,  
a \emph{(pro-)basis of $\ACHm{C}{\mdl}$} 
($\see$Definition \ref{basis-ACHm} resp.\ \ref{pro-basis-ACHm}), 
under this pairing 
($\see$Rigidity Lemma \ref{Basis-determines_0} resp.\ \ref{Basis-determines_p}). 

\bigskip
\textbf{Suitable Varieties:} 
The interplay between varieties and 
families of lower dimensional subvarieties 
often makes it necessary to consider singular (sub)varieties as well. 
In order to allow for the construction of generalized Albanese varieties, 
we use a process called \emph{suitabilization} 
($\see$\ref{suitabilization}), 
which is based on blowing up and works for any projective variety over any perfect field 
($\see$Proposition \ref{suitable_exists}), 
avoiding resolution of singularities. 
This is the content of Section \ref{Sec:SuitableVar}.

\bigskip
\textbf{Log versus Non-Log Version:} 
There are different ways in which ramification is captured by the modulus $\mdl$, 
according to different filtrations of the Witt group due to Brylinski \cite{Br} and Matsuda \cite{Mda}, 
which have their origin in Kato's class field theory. 
Here the name ``log'' is related to the fact that Brylinski's filtration can be realized as 
the preimage of certain differentials with log poles 
(see e.g.\ \cite[4.6]{KR10}). 
Each version has its advantages and disadvantages: 
the non-log version is geometrically more intuitive and easier to implement 
(especially if $X$ is singular), 
while the log version is compatible with the notions of tame and wild ramification  
(cf.\ Leitfaden Appendix \ref{Leitfaden}, Section \ref{ChowMod}). 

This paper was initially written as a pure log version, like \cite{KR10} and \cite{Ru13} (in contrast to e.g.\ \cite{KeSa}, which is non-log). 
It turned out that for all objects a non-log version could be defined in a similar fashion, 
and all assertions of this paper hold for both versions. 
Therefore I treated both versions in parallel. 
This should increase the possibility to compare the results 
to those obtained by other methods. 

\bigskip 
I included a ``Leitfaden'' as an appendix 
that reflects the way of procedure and contains some more intermediate results. 
It might also serve as an orientation 
in case the reader gets lost. 


\vspace{20mm} 

\textbf{Acknowledgements.} 
I am much indebted to Kazuya Kato for our joint work regarding the modulus 
\cite{KR10}, which is fundamental for this work. 
Marc Levine directed my attention to Bloch's ``easy'' moving lemma 
and the reciprocity trick, which became an essential tool in this work; 
thank you very much.  
I thank Shuji Saito for pointing out a mistake in a former version of this paper. 
I had helpful discussions with 
Michel Brion, H\'el\`ene Esnault, Toshiro Hiranouchi, Kazuya Kato, Moritz Kerz, 
Marc Levine, Kay R\"ulling, Shuji Saito and Rin Sugiyama, 
and I am very thankful to all of them for questions, comments and hints. 
Moreover, I thank a referee for a virtually endless list of recommendations of 
how to make the exposition more readable.

\newpage 

\subsection{Notations and Conventions} 
\label{Terminology}

We use the notations and conventions from \cite{Ru13}: 

Let $k$ be a ring. 
$\Ab$ denotes the category of abelian groups, 
$\Algk$ the category of $k$-algebras, 
$\Artk$ the category of artinian $k$-algebras. 
$\Fctr\pn{\mathsf{A},\mathsf{B}}$ denotes the category of functors from 
category $\mathsf{A}$ to category $\mathsf{B}$. 
We refer to the objects of $\Fctr\pn{\Algk,\Ab}$ as \emph{$k$-group functors} 
and to the objects of $\Fctr\pn{\Artk,\Ab}$ as \emph{formal $k$-group functors}. 

A \emph{finite ($k$-)ring} means an artinian $k$-algebra. 

The process of \emph{completion}, 
i.e. replacing an object by (a limit of) finite approximations, 
is for functors expressed by 
\;$\wh{\lul}: \Fctr\pn{\Algk,\Ab} \lra \Fctr\pn{\Artk,\Ab}$, 
induced by the natural inclusion \;$\Artk \inj \Algk$. 
This completion has a left adjoint 
\;$\wt{\lul}: \Fctr\pn{\Artk,\Ab} \lra \Fctr\pn{\Algk,\Ab}$, 
given by 
\[ \wt{\fmlG}(R) = \varinjlim_{\textrm{finite } S \subset R} \fmlG(S) 
\] 
(cf. \cite[Rmk.~1.1]{Ru13}). 
Via this map $\wt{\lul}$ we will consider the category of formal $k$-group functors 
as a subcategory of the category of $k$-group functors, 
omitting the $\wt{\phantom{a}}$. 

By an \emph{algebraic group} we always mean a \emph{commutative} algebraic group. 
Likewise, an \emph{affine group} or a \emph{formal group} always denotes 
a commutative one. 
We consider algebraic $k$-groups as $\fppf$-sheaves on $\Algk$ 
with values in $\Ab$. 

The Weil restriction of an $R$-group scheme $G$ from $R$ to $k$, 
where $R$ is a $k$-algebra, 
is denoted by $\WeilRes{R}{k} G$, 
and sometimes (as a functor) 
by $\WRS{G}{R}{k}$. 
We denote $\Lin_R := \WeilRes{R}{k} \Gm$, $\Trz_R := \WeilRes{R_{\red}}{k} \Gm$ 
and $\Upf_R := \ker\pn{\Lin_R \ra \Trz_R}$. 

Let 
$\Expaht\pn{w,t} = \exp \bigpn{ -\sum_{r\geq 0} \frac{t^{p^{r}}}{p^{r}} \,\Phi_r\pn{w}}$ 
be the Artin-Hasse exponential in the sense of \cite[V, \S4, No.~4]{DG}, 
where $k$ is a field of characteristic $p > 0$, $R$ is a $k$-algebra, 
$w = \pn{w_0, w_1, \ldots} \in \Witt\pn{R}$ is a Witt vector 
with $w_i \in \Nil(R)$ for all $i$ and $w_i = 0$ for almost all $i$, 
$t$ is a variable and 
$\Phi_r\pn{w} = \sum_{i=0}^{r} p^i w_i^{p^{r-i}}$. 
Then we define \;$\Expah\pn{w} := \Expaht\pn{w,1}$. 
\footnote{
This notation differs slightly from the one in \cite[Exm.~1.11]{Ru13}, 
where the notation \;$\Expah\pn{w,t} = \Expaht\pn{w,t}$ was used. 
}

A \emph{variety} is a separated reduced scheme of finite type over a field, 
not necessarily irreducible. 
For two subschemes $C$ and $D$ of a scheme $X$, 
the relation $C \iscp D$ means that $C$ intersects $D$ properly. 

The support of a divisor $D$ is denoted by $\vrt{D}$. 
The pull-back of a Cartier divisor $D$ to a scheme $Y$ is denoted by $D \cut Y$. 

If $X$ is a scheme, then $\sK_X$ denotes the sheaf of total quotient rings 
of the structure sheaf $\sO_{X}$. 
If $X$ is irreducible, then $K_{X}$ denotes its function field. 

\medskip 

A tuple $\pn{x_i}_{1 \leq i \leq n}$ is considered as a row-vector 
(e.g.\ in matrix equations).

\newpage 

\section{Suitable Varieties} 
\label{Sec:SuitableVar}


\bDef 
\label{suitable_Def}
A projective variety $X$ is called \emph{suitable} 
if it satisfies the following properties: 

\begin{tabular}{rl} 
(S1) & Every rational map from $X$ to an abelian variety $A$ \\ 
 & extends to a morphism from $X$ to $A$. \\ 
(S2) & $X$ is normal. 
\end{tabular} 
\eDef 

\bRmk 
\label{Pic-abelian}
A suitable projective variety $X$ over a perfect field satisfies: 

\begin{tabular}{rl} 
(S3) & The Picard variety $\Picor{X}$ of $X$ is an abelian variety. 
\end{tabular} 
\eRmk 

\bPf 
According to \cite[Thm.~5.4]{K}, 
the identity component $\Pic_{X}^0$ of the Picard scheme is projective. 
As the base field $k$ is perfect, 
the underlying reduced scheme $\Picor{X}$ is an abelian variety, 
which gives condition (S3). 
\ePf 

\bPrp 
\label{smooth-impl-suitable}
Every smooth variety is suitable. 
\ePrp 

\bPf 
A rational map to an abelian variety is defined at every smooth point 
by \cite[II, \S1, Thm.~2]{La}. This gives (S1). 
Smooth points are regular by \cite[Lem.\ 6.26]{GW}, 
and regular points are normal by \cite[Cor.\ 6.39]{GW}. 
This yields (S2). 
\ePf 

\bPrp 
\label{suitable_exists}
For every projective variety $X$ over a perfect field 
there exists a birational projective morphism $\mfn: Y \ra X$ 
such that $Y$ is suitable 
and which is an isomorphism over the regular locus of $X$. 
\ePrp 

\bPf 
We give an explicit construction of \,$\mfn: Y \ra X$\, 
in \Point \ref{suitabilization}. 
\ePf 

\bPnt 
\label{suitabilization}
Let $X$ be a projective variety over a perfect field. 
We are going to construct a birational projective morphism $\mfn: \Xt \ra X$ 
with the desired properties of Proposition~\ref{suitable_exists}, 
called a \emph{suitabilization of $X$}. 

Consider the universal rational map $\alp: X \dra \Albfct{X}$\
from $X$ to the Albanese of $X$ in the sense of functions 
($\seecite$\cite[II, \S3]{La}). 
As $\Albfct{X}$ is an abelian variety, it is projective 
($\seecite$\cite[II, Thm.~8.12]{P}). 
Let $\emb: \Albfct{X} \lra \Prj^r$ be an embedding of $\Albfct{X}$ 
into a projective space. 
Then the composition $\emb \circ \alp: X \dra \Prj^r$ 
lifts to a morphism $\bet: V \lra \Prj^r$, 
where $\sig: V \lra X$ is the blowing up of a suitable ideal sheaf $\sJ$ on $X$ 
($\seecite$\cite[II, Exm.~7.17.3]{H}). 
The closed subscheme corresponding to $\sJ$ has support equal to 
$X \setminus U$, 
where $U$ is the open set on which $\alp$ 
is regular. 
In particular $U \supset X_{\reg}$, 
i.e. the support of $\sJ$ is disjoint from the regular locus of $X$. 
\[ \xymatrix{ V \ar[d]_{\sig} \ar[r]^{\bet}  &  \Prj^r \\ 
                     X \ar@{-->}[r]^-{\alp}  &  \Albfct{X} \ar[u]_{\emb} 
   }
\] 
The image of $\bet$ is the closure of the image of $\iota \circ \alp$ in $\Prj^r$, 
which is contained in the image of $\Albfct{X}$ in $\Prj^r$: 
\[ \bet\pn{V} \,\subset\, \ol{\bet\pn{\sig^{-1}U}} \,=\, \ol{\emb \alp\pn{U}} 
   \,\subset\, \emb\bigpn{\Albfct{X}} \laurin 
\] 
Thus $\bet$ induces a morphism $\alp_V: V \lra \Albfct{X}$ 
extending $\alp \circ \sig|_{\sig^{-1}U}$. 
Since $\alp: X \dra \Albfct{X}$ is defined by the universal mapping property 
for rational maps from $X$ to abelian varieties, 
every rational map $\phe: X \dra A$ factors as 
$\phe = h \circ \alp: X \dra \Albfct{X} \lra A$ 
and hence extends to a morphism 
$\phe_V = h \circ \alp_V: V \lra \Albfct{X} \lra A$. 
\[ \xymatrix{ V \ar[dd]_{\sig} \ar[dr]^{\alp_V} \ar@/^1.3pc/[drr]^{\phe_V}  &  & \\ 
    & \Albfct{X} \ar[r]^-{h \;}  &  A \\ 
   X \ar@{-->}[ur]_{\alp} \ar@{-->}@/_1.3pc/[urr]_{\phe}  &  & \\ 
   }
\] 
Let $\nu: Y \lra V$ be the normalization of $V$, 
so $Y$ satisfies (S2). 
A morphism from $V$ to $A$ yields a morphism from $Y$ to $A$ 
via pull-back $\nu^*$. 
Then $\mfn:= \sig \circ \nu: Y \lra V \lra X$ is a birational morphism 
that satisfies condition (S1). \\ 
Thus \,$\Xt := Y$\, is a suitabilization of $X$. 
\ePnt 

\bRmk 
\label{birInvariant}
Let $X$ be a projective variety over a perfect field $k$, 
and let $\Xt$ be its suitabilization. 
By construction we have 
\[ \Alba{\Xt} = \Albfct{X} 
\laurin 
\] 
Thus $\Alba{\Xt}$ as well as its dual $\Picor{\Xt}$ 
are birational invariants among suitable varieties. 
\eRmk

\newpage 

\section{Relative Chow Group with Modulus} 
\label{ChowMod}

\subsection{Definition and Basic Properties}
\label{CHm_Basics}

Let $X$ be a projective variety over a field $k$, 
not necessarily smooth or irreducible. 
Let $\mdl$ be an effective Cartier divisor, a \emph{modulus} for $X$, 
defined as follows. 

In the sequel, we will treat log and non-log theory in a parallel manner. 
For the non-log theory, $\mdl$ lives on $X$. 
For the log theory, we need a locally factorial projective variety $\Va$ over $k$ 
that admits a morphism $\va: X \ra \Va$, 
and $\mdl$ is an effective Cartier divisor on $\Va$ 
that intersects $\va(X)$ properly. 
(The reason for the need of $\Va$ will become clear in the definition 
of the relative Chow groups with modulus, Definition \ref{DefCHoMod}, 
where we use the notion (\ref{L}) below.)
For example, if $\emb: X \inj \Prj^m$ is an embedding into a projective space, 
we could choose $\Va = \Prj^m$ and $\va = \emb$; 
this picture we should have in mind. 
If $X$ is locally factorial (e.g.\ if $X$ is smooth), 
we can choose $\Va = X$ and $\va = \id$. 
When fixing a modulus $\mdl$ for $X$, we will always assume that $\mdl$ 
comes with a fixed morphism $\va: X \ra \Va$ 
(even if this is not mentioned explicitly). 

\bNot 
\label{Not_curve}
A \emph{variety} means a separated reduced scheme of finite type 
over a field, not necessarily irreducible. 
A \emph{subvariety} of a scheme is a closed subscheme 
which is a variety.  
By a \emph{curve} we mean a 
purely 1-dimensional separated reduced scheme of finite type, 
not necessarily normal or irreducible. 
By a \emph{curve in $X$} we mean a 
closed subscheme of $X$ that is a curve. 

If $Y$ is a projective variety, 
then \,$\sig: \Yt \ra Y$\, denotes a suitabilization. 
For $f \in \sK_Y$ (= total quotient ring), 
let $\wt{f} := \sig^{*}f$ denote the image of $f$ in $\sK_{\Yt}$. 
If $Y$ is a subvariety of $X$, 
we write \,$\iot{Y}: Y \inj X$\, 
for the embedding of $Y$ into $X$. 
If $C$ is a curve, then $\Ct \ra C$ denotes its normalization 
(and this is a suitabilization). 

If \,$\morphism: Y \ra X$\, is a morphism of varieties 
whose image $\morphism\pn{Y}$ meets a Cartier divisor $D$ on $X$ properly, 
then \,$D \cut Y$\, denotes the pull-back of $D$ to $Y$. 
The reduced part of a divisor $D$ is denoted by $D_{\red}$. 
The Weil divisor associated to a Cartier divisor $D$ 
is denoted by $\bt{D}$. 
We define for a Cartier divisor $\mdl$ on $\Va$ 
the Weil divisor $\mdl_{Y}$ on $Y$ 
in the \emph{log} theory by
\bMyEqn \label{L} \mytag{L} 
\mdl_{Y} := \bigbt{\mdl \cut Y}_{\red} + \bigbt{\pn{\mdl-\mdl_{\red}} \cut Y} 
\laurink 
\eMyEqn 
and in the \emph{non-log} theory by 
\bMyEqn \label{NL} \mytag{NL} 
\mdl_{Y} := \bigbt{\mdl \cut Y} 
\laurin 
\eMyEqn 
For a curve $C$ in $X$, we can consider $Y = \Ct$ 
(= normalization of C), and write 
$\mdl_{\Ct} = \sum n_{\pnt} [{\pnt}]$, 
where ${\pnt}$ ranges over the closed points of $\Ct$, 
for integers $n_{\pnt} \geq 0$ and $n_{\pnt} = 0$ for almost all ${\pnt}$. 
Then ``$f \equiv 1 \mod \mdl_{\Ct}$'' is defined as 
$\val_{\pnt}(1 - f) \geq n_{\pnt}$ for all $\pnt \in C$, 
where $\val_{\pnt}$ is the valuation attached to the point $\pnt$. 
\eNot 

\bDef 
\label{DefCHoMod}
Let $\ZOm{X}{\mdl}$ 
be the group of 0-cycles on $X \setminus \mdl_X$, set 
\[ \Rlnpm{X}{\mdl} 
= \lrst{ \pn{C,f} \left| 
\begin{array}{l}
C\textrm{ an irreducible curve in $X$ meeting } \mdl_X \textrm{ properly}  \\ 
f \in \sK_C^* \hspace{2mm} \textrm{ such that } \hspace{2mm} 
\wt{f} \equiv 1 \mod \mdl_{\Ct} 
\end{array}
\right.
} 
\] 
and let $\Rlnm{X}{\mdl}$ be the subgroup of 
$\ZOm{X}{\mdl}$ 
generated by the elements $\dv\pn{f}_{C}$ with 
$\pn{C,f}\in\Rlnpm{X}{\mdl}$. \\ 
The \emph{relative Chow group $\CHm{X}{\mdl,\va}$ 
of $X$ of modulus $\mdl$ along $\va$} 
(we will omit the $\va$ in the following) 
is defined as 
\[ \CHm{X}{\mdl} = \frac{\ZOm{X}{\mdl}}{\Rlnm{X}{\mdl}} \laurin 
\] 
We define 
\[ \Rlnpm{X}{0 \mdl} = \left\{ \pn{C,f}\left|\begin{array}{l}
C\textrm{ an irreducible curve in $X$ meeting } \mdl_X \textrm{ properly} \\ 
f \in \sK_C^* \hspace{2mm} \textrm{such that} \hspace{2mm} 
f \in \sO_{\Ct,\pnt}^* \hspace{5mm} \forall\, \pnt \in \mdl \cut \Ct 
\end{array}\right.\right\} 
\] 
and let $\Rlnm{X}{0 \mdl}$ be the subgroup of $\ZOm{X}{\mdl}$ 
generated by the elements $\dv\pn{f}_{C}$ with 
$\pn{C,f}\in\Rlnpm{X}{0 \mdl}$. 
Then we define 
\[ \CHm{X}{0 \mdl} 
:= \frac{\ZOm{X}{\mdl}}{\Rlnm{X}{0 \mdl}} \laurin 
\] 
We let $\CHmo{X}{\mdl}$ (resp.\ $\CHmo{X}{0 \mdl}$) 
be the subgroup of $\CHm{X}{\mdl}$ (resp.\ of $\CHm{X}{0 \mdl}$) 
of cycles $\zeta$ with $\deg\zeta|_{W}=0$ for all connected components
$W$ of $X\setminus\mdl$. 
\eDef  

\bRmk 
\label{usualCHo}
Note that for $\mdl = 0$ 
the relative Chow group with modulus becomes the usual Chow group: 
\bDpl \CHm{X}{0} = \CHO{X} \laurin 
\eDpl 
\eRmk 

\bRmk 
\label{Bloch}
If $X$ is \emph{smooth}, 
by Bloch's so called ``easy moving lemma'' \cite[Prop.~2.3.1]{Bl05} it holds 
\[ \CHm{X}{0 \mdl} = \CHO{X} 
\laurin 
\] 
\eRmk 


\bPrp 
\label{FunctorCHoMod}
Let \,$\morphism: Y \ra X$\, be a morphism of projective varieties over $k$, 
and consider \,$\va \circ \morphism: Y \ra \Va$\, 
as the defining morphism for a modulus on $Y$ 
(in the log theory). 
Let $\mdll$ be an effective Cartier (= Weil) divisor on $\Va$ 
with $\mdl \leq \mdll$ in the log theory, 
respectively let $\mdll$ be an effective Cartier divisor on $Y$ 
such that there is an embedding $[\morphism^*\mdl] \subset [\mdll]$ 
of underlying subschemes on $Y$ in the non-log theory. 
Then $\morphism$ induces a homomorphism of groups 
\[ \morphism_*: \CHm{Y}{\mdll} \lra \CHm{X}{\mdl} 
\laurin 
\] 
(Here \,$\mdl \leq \mdll$\, is defined as usual by the following condition: 
If \,$\mdl = \sum n_v [v]$\, and \,$\mdll = \sum m_v [v]$\, 
for integers $n_v, m_v \geq 0$, 
where $v$ ranges over all points of codimension 1 in $\Va$, 
then \,$\mdl \leq \mdll$\, if and only if \,$n_v \leq m_v$\, for all $v$.) 
\ePrp 

\bRmk[Functoriality of $\CHm{X}{\mdl}$]
\label{FunctorialityCHoMod}
The assertion of Proposition \ref{FunctorCHoMod} 
is equivalent to the following two statements: 

(a) Any morphism $\morphism: Y \ra X$ 
in the category of projective schemes over $k$ 
induces a homomorphism 
\[ \morphism_*: \CHm{Y}{\mdl} \lra \CHm{X}{\mdl}
\]
respectively 
\[ \morphism_*: \CHm{Y}{\morphism^*\mdl} \lra \CHm{X}{\mdl} 
\laurin 
\] 

(b) Any relation $\mdl \leq \mdll$ of divisors on $\Va$, 
respectively $[\mdl] \subset [\mdll]$ of Cartier-divisors on $X$, 
induces a homomorphism 
\[ \CHm{X}{\mdll} \lra \CHm{X}{\mdl} 
\laurin 
\] 

Then $\CHm{\lul}{\mdl}$ becomes a \emph{covariant functor} 
from the category of projective schemes over $\Va$ resp.\ $X$ 
to the category of abelian groups, 
whereas $\CHm{X}{\lul}$ becomes a \emph{contravariant functor} from the category 
of Cartier divisors on $\Va$ resp.\ $X$ to the category of abelian groups. 
Here the morphisms between Cartier divisors on $\Va$ resp.\ $X$ 
are given by embeddings between the underlying subschemes of $\Va$ resp.\ $X$. 
\eRmk 

\bPf 
The second statement of Remark \ref{FunctorialityCHoMod} 
is clear since 
$\ZOm{X}{\mdll} \subset \ZOm{X}{\mdl}$ and 
$\Rlnm{X}{\mdll} \subset \Rlnm{X}{\mdl}$ for $\mdl \leq \mdll$ 
resp.\ $[\mdl] \subset [\mdll]$. 

Let $\morphism: Y \ra X$ be a morphism of varieties 
and $\mdll \geq \mdl$ an effective divisor on $\Va$ 
(in the log theory), 
respectively $\mdll$ an effective Cartier-divisor on $Y$ with 
$[\mdll] \supset [\morphism^*\mdl]$ (in the non-log theory). 
In view of the above we are reduced to the case $\mdll = \mdl$, 
respectively $\mdll = \morphism^*\mdl$. 
For $C \subset Y$ we have $\pn{\morphism^*\mdl}_{\Ct} = \mdl_{\Ct}$ 
(in the non-log theory, while the corresponding assertion 
in the log theory is trivial by definition).  
For $f \in \sK_C^*$ we have 
$\morphism_* \dv(f)_C = 0$ if $\dim C > \dim \morphism(C)$, 
and $\morphism_* \dv(f)_C = \dv\bigpn{\Norm\pn{f}}_{\morphism(C)}$ 
if $\dim C = \dim \morphism(C)$, 
where $\Norm\pn{f} \in \sK_{\morphism(C)}^*$ is the norm of $f$ 
(see \cite[Chap.~1, Prop.~1.4]{F}). 
Since $\wt{f} \equiv 1 \mod \mdl_{\Ct}$ implies 
$\wt{\Norm\pn{f}} \equiv 1 \mod \mdl_{\wt{\morphism(C)}}$ 
(cf.\ \cite[III, No.~2, proof of Prop.~4]{S59}), 
the push-forward of cycles 
$\phe_*: \ZOm{Y}{\mdll} \lra \ZOm{X}{\mdl}$, 
as it maps $\Rlnm{Y}{\mdll} \lra \Rlnm{X}{\mdl}$, 
induces a homomorphism $\CHm{Y}{\mdll} \lra \CHm{X}{\mdl}$. 
\ePf

\newpage 
\subsection{Description of the Chow Group for Singular Curves} 
\label{CHm_singular}

\bDef 
\label{Def_CHmSing}
Let $C$ be a curve in $X$ and \;$\nu: \Ct \ra C$\; its normalization.  
Let $\Sing\pn{C,\mdl} \subset \Sing(C)$ 
be the subset of the singular locus of $C$ 
consisting of those irreducible components of $\Sing(C)$ 
that intersect $\mdl_X$ properly 
(i.e.\ do not meet $\mdl_X$, since $\Sing(C)$ is of dimension 0). 
Define 
\[ \KZm{C}{\mdl_{\Cgst}} = 
	\ker\Bigpn{\Z_0\bigpn{\Ct \tms_C \Sing\pn{C,\mdl}} 
	\lra \Z_0\bigpn{\Sing\pn{C,\mdl}}} 
\] 
where the map is given by the push-forward of cycles $\nu_*$. 

Then $\QHm{C}{\mdl_{\Cgst}}$ is defined to be the cokernel 
of the natural homomorphism 
from $\KZm{C}{\mdl_{\Cgst}}$ to $\CHm{\Ct}{\mdl}$: 
\[ \QHm{C}{\mdl_{\Cgst}} = \coker\Bigpn{\KZm{C}{\mdl_{\Cgst}} \lra \CHm{\Ct}{\mdl}} 
\laurin
\] 
Similarly we define 
\[ \QHm{C}{0 \mdl_{\Cgst}} 
= \coker\Bigpn{\KZm{C}{\mdl_{\Cgst}} \lra \CHm{\Ct}{0 \mdl_{\Ct}}} 
\laurin 
\] 
\eDef 



\bRmk[Functoriality of $\QHm{C}{\mdl}$] 
\label{QHm-functoriality}
If $C$ and $Z$ are curves in $X$ 
and $\zeta: Z \lra C$ is a morphism over $X$, 
then $\zeta$ is an embedding. 
We get an induced morphism $\wt{\zeta}: \Zt \lra \Ct$ 
on their normalizations. 
Then the push-forward of cycles $\wt{\zeta}_*$ 
induces a homomorphism $\CHm{\Zt}{\mdl} \lra \CHm{\Ct}{\mdl}$ 
by Remark \ref{FunctorialityCHoMod} (a), 
and $\wt{\zeta}_*$ maps $\KZm{Z}{\mdl}$ to $\KZm{C}{\mdl}$, 
since $\zeta$ is injective. 
Thus we obtain a functoriality homomorphism 
\[ \zeta_*: \QHm{Z}{\mdl} \lra \QHm{C}{\mdl} 
\laurin 
\] 
Suppose $\mdll \geq \mdl$ is an effective Cartier divisor on $\Va$, 
resp.\ $C$. 
Then \,$\Sing\pn{C,\mdll} \subset \Sing\pn{C,\mdl}$, 
inducing a homomorphism \,$\KZm{C}{\mdll} \lra \KZm{C}{\mdl}$. 
By Remark \ref{FunctorialityCHoMod} (b) 
we have a homomorphism 
\,$\CHm{\Ct}{\mdll} \lra \CHm{\Ct}{\mdl}$. 
Putting together we obtain a homomorphism 
\[ \QHm{C}{\mdll} \lra \QHm{C}{\mdl} 
\laurin 
\] 
\eRmk 

\bRmk 
\label{Char_CHmSing}
If $C$ and $Z$ are two curves in $X$ that meet each other and $\mdl$ properly, 
then $\QHm{C \cup Z}{\mdl_{\CcupZgst}}$ is characterized by the push-out 
\[ \xymatrix{ 
	& \CH_0\bigpn{\Ct \tms_X \Zt, \mdl} 
	\ar[dl] \ar[dr]  & \\ 
	\QHm{C}{\mdl_{\Cgst}} \ar[dr] & & \QHm{Z}{\mdl_{\Zgst}} \ar[dl] \\ 
	& \QHm{C \cup Z}{\mdl_{\CcupZgst}} \laurin 
}
\] 
\eRmk 

\bPrp 
\label{CHmSing-motivation}
Suppose that one of the following conditions holds: 

\begin{tabular}{rl} 
{\rm (QH1)} & The curve $C$ is normal. \\ 
{\rm (QH2)} & The base field $k$ is algebraically closed. 
\end{tabular} 

\noindent 
Then the groups from Definition \ref{DefCHoMod} 
and Definition \ref{Def_CHmSing} coincide: 
\[ \CHm{C}{\mdl} = \QHm{C}{\mdl} 
\laurin 
\]  
\ePrp 

\bPf 
If $C = \Ct$, then $\KZm{C}{\mdl} = 0$, 
thus $\QHm{C}{\mdl} = \CHm{C}{\mdl}$. 

\noindent 
If the base field $k$ is algebraically closed, 
then all residue fields at closed points are equal to $k$. 
Thus the push-forward of cycles 
$\sig_*: \ZOm{\Ct}{\mdl} \lra \ZOm{C}{\mdl}$ 
is surjective. 
Set $S := \Sing\pn{C,\mdl}$ and $\St := S \tms_C \Ct$. 
We obtain the following commutative diagram, 
in which all diagonal and vertical sequences are exact: 
\[ \xymatrix{ 
0 \ar[dr] & \KZm{C}{\mdl} \ar[d] \ar@{=}[dr] & & \Ker\CHO{\sig} \ar[d] & 0 \\ 
0 \ar[dr] & \ZOm{\St}{\mdl} \ar@{->>}[d] \ar[dr] & 
\KZm{C}{\mdl} \ar[d] \ar[ur]^{\kap} & \CHm{\Ct}{\mdl} \ar@{->>}[d] \ar[ur] & 0 \\ 
 & \ZOm{S}{\mdl} \ar[dr]|!{[d];[r]}\hole & 
\ZOm{\Ct}{\mdl} \ar@{->>}[d] \ar[dr] \ar[ur] 
 & \CHm{C}{\mdl} \ar[ur] & \\ 
 & \Rlnm{\Ct}{\mdl} \ar@{->>}[d] \ar[ur] & 
\ZOm{C}{\mdl} \ar[dr] \ar[ur]|!{[u];[r]}\hole & 
\ZOm{\Ct}{\mdl_{\Ct} \cup \St} \ar@{=}[d] \ar[dr] & \\ 
0 \ar[ur] & \Rlnm{C}{\mdl} \ar[d] \ar[ur] & & 
\ZOm{C}{\mdl_C \cup S} \ar[dr] & 0 \\ 
0 \ar[ur] & 0 & & & 0 
} 
\] 
The arrow \;$\kap: \KZm{C}{\mdl} \lra \Ker\CHO{\sig}$\; 
is surjective by the Snake Lemma. 
This shows 
\;$\CHm{C}{\mdl} 
= \coker\bigpn{\KZm{C}{\mdl} \lra \CHm{\Ct}{\mdl}}$. 
\ePf 

\subsection{Affine Part of the Chow Group with Modulus}

\bDef 
\label{Def_ACHm}
Define $\ACHm{X}{\mdl}$ to be the kernel of the natural homomorphism 
\;$\CHm{X}{\mdl} \lra \CHm{X}{0\mdl}$: 
\begin{eqnarray*} 
\ACHm{X}{\mdl} 
 & = & \ker\Bigpn{\CHm{X}{\mdl} \lra \CHm{X}{0\mdl}} \\ 
 & = & \Rlnm{X}{0 \mdl} \big/ \Rlnm{X}{\mdl} \laurin 
\end{eqnarray*} 
In the same way we define $\AQHm{C}{\mdl}$ for a curve $C$ in $X$: 
\begin{eqnarray*} 
\AQHm{C}{\mdl} 
 & = & \ker\Bigpn{\QHm{C}{\mdl} \lra \QHm{C}{0\mdl}} 
\laurin 
\end{eqnarray*} 
\eDef 

\bRmk  
\label{kerCHm}
The map $\CHm{X}{\mdl} \lra \CHm{X}{0 \mdl}$ 
is obviously surjective. 
We have thus exact sequences 
\[ 0 \lra \ACHm{X}{\mdl} \lra \CHm{X}{\mdl} \lra \CHm{X}{0 \mdl} \lra 0 
\] 
and 
\[ 0 \lra \ACHm{X}{\mdl} \lra \CHmo{X}{\mdl} \lra \CHmo{X}{0 \mdl} \lra 0 
\] 
where the superscript ``$0$'' refers to 0-cycles of degree 0 
on every connected component of $X$. 

The functoriality of $\CHm{X}{\mdl}$ and $\CHm{X}{0 \mdl}$ 
implies that $\ACHm{X}{\mdl}$ 
admits the same functorial properties as in Remark \ref{FunctorialityCHoMod}.  

Likewise, as $\QHm{C}{\mdl} \lra \QHm{C}{0 \mdl}$ is surjective 
(see in the proof of Lemma \ref{CHm-surj}, Step 2), 
the functoriality of $\QHm{C}{\mdl}$ and $\QHm{C}{0 \mdl}$ 
implies that $\AQHm{C}{\mdl}$ 
admits the same functorial properties as in Remark \ref{QHm-functoriality}.  
\eRmk 

\subsection{Chow Group as Direct Limit over Curves}
\label{CHm_as_limit}

Suppose the base field $k$ is finite or algebraically closed. 
Let $\mdl$ and $\mdll$ be effective divisors on $\Va$ resp.\ $X$. 
If $k$ is not algebraically closed, 
we suppose furthermore that $X$ is smooth outside of $\vrt{\mdl_X}$. 

\bPrp 
\label{CHm_indLim}
Given the assumptions as above, 
then we have the following presentations of 
\,$\CHm{X}{\mdl}$ and $\ACHm{X}{\mdl}$ as inductive limits: 
\begin{eqnarray*} 
\CHm{X}{\mdl} & = & \varinjlim_{C} \QHm{C}{\mdl_{\Cgst}} \\ 
\ACHm{X}{\mdl} & = & \varinjlim_{C} \AQHm{C}{\mdl_{\Cgst}} 
\end{eqnarray*} 
where $C$ ranges over all curves in $X$ 
that meet $S := \vrt{\mdl_X + \mdll_X}$ properly. 
These curves form a directed system 
with transition maps given by inclusion: 
for two curves $C_1$ and $C_2$ the curve $Z = C_1 \cup C_2$ in $X$ 
satisfies $C_1 \subset Z$, $C_2 \subset Z$. 
\ePrp 

\bPf 
For all 0-cycles $z,t$ on $X$ there is a 
hypersurface $H_1$ containing $z$ and avoiding $t$, 
which is smooth on the smooth locus of $X$. 
This is evident for $k$ algebraically closed. 
For $k$ finite see \cite[Cor.~1.6]{Ga}. 
If $t$ contains a point of every irreducible component of $S$, 
then this implies that $H_1$ meets $S$ properly. 
By induction we find for every $z \in \ZOm{X}{\mdl}$ a 
curve $C = H_{\dim X -1}$ in $H_{\dim X -2}$ 
that contains $z$ and meets $S$ properly. 
If $k$ is not algebraically closed, 
then by assumption on $X$ for this case 
we may further suppose that C is smooth outside of $\vrt{\mdl_X}$. 
Hence the map 
\,$\varinjlim \QHm{C}{\mdl_{\Cgst}} \lra \CHm{X}{\mdl}$\, 
is surjective. 

Next we show injectivity. 
Suppose $\kap \in \ker\bigpn{\QHm{C}{\mdl_{\Cgst}} \xra{\iota} \CHm{X}{\mdl}}$. 
We have $\QHm{C}{\mdl} = \CHm{\Ct}{\mdl} \big/ \KZm{C}{\mdl}$. 
Let $\wt{z} \in \ZOm{\Ct}{\mdl}$ be a representative of $\kap$. 
There are $\pn{Z_i, f_i} \in \Rlnpm{X}{\mdl}$ 
such that $\sum \dv\pn{f_i}_{Z_i} = \pn{\iota \nu}_* \wt{z} \in \ZOm{X}{\mdl}$. 
Let $Z := \bigcup Z_i$, 
let $\del$ be the Cartier divisor on $\coprod Z_i$ formed by the $f_i$, 
and $\wt{\del} := \nu^*\del \in \Div\pn{\Zt} \cong \ZO{\Zt}$ 
the pull-back of $\del$ to the normalization $\Zt$ of $Z$. 
Moving the representative $\wt{z}$ of $\kap$ by rational equivalence 
$\mod \mdl_{\Ct}$, 
we may assume that the curves $Z_i$ intersect $C$ properly: 
indeed, if $Z_j$ is a component of $C$ for some $j$, 
then $\bigbt{\dv\pn{\wt{f}_j}_{\wt{Z}_j}} = 0 \in \CHm{\Ct}{\mdl}$, 
hence we can replace $\wt{z}$ by $\wt{z} - \dv\pn{\wt{f}_j}_{\wt{Z}_j}$, 
thus erasing $\pn{Z_j, f_j}$ from the list of involved pairs $\pn{Z_i, f_i}$. 
Then $C \cap Z \subset \Sing\pn{C \cup Z}$, 
and thus \;$\wt{z} - \wt{\del} \in \KZm{C \cup Z}{\mdl_{\CcupZgst}}$. 
If \,$\eta: C \lra C \cup Z$\, denotes the embedding, 
then \,$\eta_*\kap = \bt{\wt{z}} = \bt{\wt{\del}} = 0 
\in \QHm{C \cup Z}{\mdl}$, 
showing that the image of $\kap$ in \,$\varinjlim \QHm{C}{\mdl_{\Cgst}}$\, is $0$. 

It remains to show the statement about $\ACHm{X}{\mdl}$. 
Consider the following commutative diagram with exact rows: 
\[ \xymatrix{ 
0 \ar[r] & \varinjlim \AQHm{C}{\mdl_{\Cgst}} \ar[r] \ar[d] & 
\varinjlim \QHm{C}{\mdl_{\Cgst}} \ar[r] \ar[d]^{\wr} & 
\varinjlim \QHm{C}{0 \mdl_{\Cgst}} \ar[r] \ar[d]^{\wr} 
& 0 \\ 
0 \ar[r] & \ACHm{X}{\mdl} \ar[r] & \CHm{X}{\mdl} \ar[r] & \CHm{X}{0 \mdl} \ar[r] & 0 
} 
\] 
where the vertical maps in the middle and on the right are isomorphisms 
by the statement about $\CHm{X}{\mdl}$ and in particular $\CHm{X}{0 \mdl}$. 
Then the vertical map on the left is an isomorphism by the Five Lemma. 
\ePf 

\vspace{\vs} 

In view of Proposition \ref{CHmSing-motivation} 
we can reformulate the result of Proposition \ref{CHm_indLim} as follows: 

\bPrp 
\label{CHm_indLim_algclosed}
Given the assumptions as above, 
if the base field $k$ is algebraically closed, then 
\begin{eqnarray*} 
	\CHm{X}{\mdl} & = & \varinjlim_{C} \CHm{C}{\mdl_{\Cgst}} \\ 
	\ACHm{X}{\mdl} & = & \varinjlim_{C} \CHm{C}{\mdl_{\Cgst}} 
\end{eqnarray*} 
where $C$ ranges over all curves in $X$ 
that meet $\vrt{\mdl_X}$ properly. 
\ePrp

\bLem 
\label{CHm-surj}
Given the assumptions as above, 
then the functoriality maps \\  
\;$ \CHm{X}{\mdll} \lra \CHm{X}{\mdl} 
$\; 
and 
\;$ \ACHm{X}{\mdll} \lra \ACHm{X}{\mdl} 
$\;  \\ 
are surjective for $\mdll \geq \mdl$. 
\eLem 

\bPf 
\textbf{Step 1:} 
First we can reduce to the case that $X$ is a curve $C$, 
according to Proposition \ref{CHm_indLim} 
and the exactness of the inductive limit. 
(Here it will suffice to show the statement for 
$\QHm{C}{\lull}$ instead of $\CHm{C}{\lull}$.) 

\textbf{Step 2:} 
Next we reduce to the case of a normal curve $\Ct$ 
via the following commutative diagram with exact rows: 
\[ \xymatrix{ 
& \KZm{C}{\mdll} \ar[r] \ar[d] & 
\CHm{\Ct}{\mdll_{\Ct}} \ar[r]^{\rho_{\mdll}} \ar[d]^-{\wt{\del}} & 
\QHm{C}{\mdll} \ar[r] \ar[d]^-{\del} & 0 \phantom{\laurin} \\ 
& \KZm{C}{\mdl_{\Cgst}} \ar[r] & \CHm{\Ct}{\mdl_{\Ct}} \ar[r]^{\rho_{\mdl}} & 
\QHm{C}{\mdl_{\Cgst}} \ar[r] & 0 \laurin 
} 
\] 
If $\wt{\del}$ is surjective, 
then $\del \circ \rho_{\mdll} = \rho_{\mdl} \circ \wt{\del}$ is surjective, 
hence $\del$ is surjective. 

\textbf{Step 3:} 
Due to Remark \ref{kerCHm} and Remark \ref{Bloch}, 
for a normal curve (= smooth curve, as the base field is perfect) 
we have a decomposition of $\CHm{C}{\lull}$ 
as extension of the group $\CHO{C}$ by $\ACHm{C}{\lull}$. 
Then by the Five Lemma, 
the assertion follows from the corresponding statement 
for $\ACHm{C}{\lull}$ instead of $\CHm{C}{\lull}$. 

\textbf{Step 4:} 
For each irreducible component $C$ 
of a smooth curve (= normal curve, as over a perfect field), 
the description of $\ACHm{C}{\mdl_{\Cgst}}$ 
from Definition \ref{Def_ACHm} 
and the Approximation Lemma ($\see$\cite[I, \S3]{S66}) yield 
\[ 
\ACHm{C}{\mdl_{\Cgst}} 
\;=\; \frac{\bigcap_{\pnt \in \vrt{\mdl}} \sO_{C,\pnt}^*} 
       {\bigcap_{\pnt \in \vrt{\mdl}} \pn{1 + \fm_{C,\pnt}^{n_{\pnt}}}} 
       \Bigg/ \extfld^* 
\;=\; \Biggpn{\prod_{\pnt \in \vrt{\mdl}} 
       \frac{\sO_{C,\pnt}^*}{1 + \fm_{C,\pnt}^{n_{\pnt}}}} 
       \Bigg/ \extfld^* 
\] 
where $\extfld := \H^0\pn{C,\sO_C}$. 
Then the statement of the Lemma becomes obvious for 
$\mdll = \sum_{\pnt \in \vrt{\mdll}} m_{\pnt} \bt{\pnt} 
\geq \sum_{\pnt \in \vrt{\mdl}} n_{\pnt} \bt{\pnt} = \mdl$. 
\ePf


\subsection{Affine Chow Group as Sum over Curves} 
\label{AuxiliaryResults}

Let $X$ be a 
projective variety 
over an algebraically closed 
field $k$. 
Let $\mdl$ be an effective divisor on $\Va$ resp.\ $X$, 
set $\Sm := \Supp\pn{\mdl \cut X}$. 

The goal of this subsection is to find a presentation of $\ACHm{X}{\mdl}$
as the sum of images of $\ACHm{C}{\mdl}$ for certain curves $C$, 
where the sum ranges over a set of curves that is supposed to be small 
and the curves are desired to be ``simple'', i.e.\ of bounded degree. 

\bNot 
\label{DefC1}
Let $\ZLm{X}{\mdl}$ be the subgroup of 
$\underset{C \subset X} \dsum K_{C}^*$
generated by 
$\Rlnpm{X}{\mdl}$, 
let \,$\partial: \ZL{X} \lra \ZO{X}$\, denote the boundary map, i.e. 
\bDpl  
\partial\lrpn{C,f} = \dv(f)_C 
\hspace{10mm} \textrm{for $C \subset X$ a curve and $f \in \sK_C^*$} 
\laurin 
\eDpl 
\eNot 

\bLem[Reciprocity Trick] 
\label{tameLift}
Let $X$ be a projective surface, 
$C$ a complete intersection curve in $X$, 
not necessarily reduced or pure-dimensional.  
Let $f \in \sK_C^*$ such that $f \in \sO_{C,\pntt}^*$ 
for all $\pntt \in \pn{C \isec \mdl_X} \cup \st{\textrm{generic points of $C$}}$. 
Then 

{\rm (a)} There exist \,$\gam \in \sK_X^*$ with \,$\dv(\gam) =: \Gam \supset C$ \\
      and \,$\tha \in \sK_X^*$ with \,$\tha|_C \equiv f \mod \mdl_{\Ct}$\,
      such that 
      \[ \bigbt{\partial\lrpn{\Gam,\tha|_{\Gam}}} 
          = \bigbt{\partial\lrpn{C,f}} \in \ACHm{X}{\mdl} 
          \laurin 
      \] 

{\rm (b)} \;For \,$\Tha := \dv(\tha)$\, we have 
\[ \partial\lrpn{\Tha,\gam|_{\Tha}} 
\,=\, \partial\lrpn{\Gam,\tha|_{\Gam}} 
\laurin 
\] 
\eLem 

\bPf 
(a) 
By assumption $C$ is a complete intersection curve, 
i.e.\ there is a line bundle $\sL$ and a global section $s \in \H^0(\sL)$ 
such that $C = \Zero(s)$ is the divisor of zeroes. 
First suppose that $\sL$ is sufficiently ample, 
i.e. there is an embedding of $X$ 
into a projective space $\Prj^m$ 
such that $\sL$ comes from $\sO_{\Prj^m}(n)$ for $n$ sufficiently large. 
Then we find global sections $t, t' \in \H^0(\sL)$ 
such that $\Zero(t)$ and $\Zero(t')$ do not have a common component 
and 
for $\tha := t / t' \in \sK_X$ it holds 
$\tha|_C \equiv f \mod \mdl_{\Ct}$.
Set $E := \Zero(t - t') = \Zero(\tha - 1)$, then $\tha|_E = 1$. 
Set $\gam := s / (t - t') \in \sK_X^*$ and $\Gam = \dv(\gam) = C - E$. 
Then $\bigbt{\partial\lrpn{\Gam,\tha|_{\Gam}}} 
          = \bigbt{\partial\lrpn{C,f}} \in \ACHm{X}{\mdl}$ 
as required. 

If $\sL$ is not sufficiently ample, 
we choose an ample line bundle $\sA$ on $X$. 
Then for some $n \in \Nat$ we find global sections 
$h, h' \in \H^0\pn{\sA^{\tens n}}$ 
such that $\Zero(h)$ and $\Zero(h')$ do not have a common component 
and $\eta := h / h' \in \sK_X$ satisfies $\eta|_C \equiv f \mod \mdl_{\Ct}$ 
as above. 
The curve $A := \Zero(h - h')$ $= \Zero(\eta - 1)$ is an ample divisor on $X$ 
with $\eta|_A = 1$. 
Then $C' := C + A$ is sufficiently ample, 
and for $f' := \eta|_{C'}$ it holds 
\,$\bt{\dv\pn{f'}_{C'}} = \bt{\dv\pn{f}_C} \in \ACHm{X}{\mdl}$. 
Then we replace $\pn{C,f}$ by $\pn{C',f'}$ and proceed as above. 

(b) The sequence 
\[ \K_2\pn{\sK_X} \overset{\ft} \lra \ZL{X} \overset{\partial} \lra \ZO{X} 
\] 
is a complex ($\seecite$\cite[(1.1)]{B}), 
where $\ft$ is the tame symbol, given by 
\[ \ft\lrst{\gam,\tha} = \tha|_{\dv(\gam)} - \gam|_{\dv(\tha)} 
\laurin 
\] 
Thus 
\;$ 
\partial\lrpn{\Gam,\tha|_{\Gam}} - \partial\lrpn{\Tha,\gam|_{\Tha}} 
\,=\, \partial\,\bigpn{\ft\lrst{\gam,\tha}} 
=\, 0 
\,\in \ZO{X} 
\laurin 
$
\ePf 

\bLem 
\label{div_C=div_Z}
Let \,$C$\, and \,$Z$\, be two curves in $X$ 
with 
\,$C \tms_X \bt{\mdl_X} = Z \tms_X \bt{\mdl_X}$, 
where $\bt{\mdl_X}$ denotes the (non-reduced) subscheme of $X$ 
given by $\mdl_X$. 
Then $\ACHm{C}{\mdl_{\Cgst}}$ and $\ACHm{Z}{\mdl_{\Zgst}}$ 
generate the same subgroup in $\ACHm{X}{\mdl}$: 
\[ {\iot{C}}_* \ACHm{C}{\mdl_{\Cgst}} = {\iot{Z}}_* \ACHm{Z}{\mdl_{\Zgst}} 
   \subset \ACHm{X}{\mdl} 
   \laurin 
\] 
In other words: 
Let $f \in \sK_C^*$ with $f \in \sO_{C,\pntt}^*$\,
for all $\pntt \in C \cap \bt{\mdl_X}$. 
Then there is $g \in \sK_Z^*$ with $g \in \sO_{Z,\pntt}^*$\, 
for all $\pntt \in Z \isec \bt{\mdl_X}$, such that 
\[ \bigbt{\dv(f)_C} = \bigbt{\dv(g)_Z} \in \ACHm{X}{\mdl} \laurin 
\] 
\eLem 

\bPf 
\textbf{Step 1:} Assume that $X$ is a surface (not necessarily smooth). 
By Lemma \ref{tameLift} we find 
$\gam \in \sK_X^*$ with \,$\dv(\gam) =: \Gam \supset C$ and 
$\tha \in \sK_X^*$ with \,$\tha|_C \equiv f \mod \mdl_{\Ct}$ such that 
$\bigbt{\partial\lrpn{\Gam,\tha|_{\Gam}}} 
= \bigbt{\partial\lrpn{C,f}} \in \ACHm{X}{\mdl}$. 
Similarly we find           
$\sig \in \sK_X^*$ with \,$\dv(\sig) =: \Sig \supset Z$ 
such that for $\tha|_Z =: g$ we have 
$\bigbt{\partial\lrpn{\Sig,\tha|_{\Sig}}} 
= \bigbt{\partial\lrpn{Z,g}} \in \ACHm{X}{\mdl}$. 
It remains to show 
$\bigbt{\partial\lrpn{\Gam,\tha|_{\Gam}}} 
= \bigbt{\partial\lrpn{\Sig,\tha|_{\Sig}}} \in \ACHm{X}{\mdl}$. 

By Lemma \ref{tameLift} (b) it holds 
$\partial\lrpn{\Gam,\tha|_{\Gam}} = \partial\lrpn{\Tha,\gam|_{\Tha}}$ and 
$\partial\lrpn{\Sig,\tha|_{\Sig}} = \partial\lrpn{\Tha,\sig|_{\Tha}}$, 
where $\Tha := \dv(\tha)$. 
The assumption on $C$ and $Z$ implies that 
we can choose $\gam$ and $\sig$ such that 
$\Gam \isec \bt{\mdl_X} = \Sig \isec \bt{\mdl_X}$, 
i.e.\ $\gam \equiv \sig \mod \mdl$. 
Then $\partial\lrpn{\Tha,\gam|_{\Tha}} \equiv \partial\lrpn{\Tha,\sig|_{\Tha}} 
\mod \Rlnm{X}{\mdl}$, 
which yields 
$\bigbt{\partial\lrpn{\Gam,\tha|_{\Gam}}} 
= \bigbt{\partial\lrpn{\Tha,\gam|_{\Tha}}} 
= \bigbt{\partial\lrpn{\Tha,\sig|_{\Tha}}} 
= \bigbt{\partial\lrpn{\Sig,\tha|_{\Sig}}}$. 

\textbf{Step 2:} Now we reduce the higher dimensional case 
to the case of a surface. Let $d := \dim X$. 
By assumption $C$ and $Z$ are complete intersection curves: 
$C = \Zero\pn{\gam_1, \ldots, \gam_{d-1}}$ and 
$Z = \Zero\pn{\sig_1, \ldots, \sig_{d-1}}$, 
where the $\gam_i$ and $\sig_j$ are global sections of a line bundle. 
Reordering $\gam_1, \ldots, \gam_{d-1}$ 
we can form the curves 
$C_i := \Zero\pn{\sig_1, \ldots, \sig_{i-1}, \gam_{i}, \ldots, \gam_{d-1}}$ 
for $i = 1, \ldots, d$. 
Then $C_1 = C$ and $C_{d} = Z$. 
By construction, $C_i$ and $C_{i+1}$ both lie in the surface 
$\Zero\pn{\sig_1, \ldots, \sig_{i-1}, \gam_{i+1}, \ldots, \gam_{d-1}}$. 
Thus by Step 1, with $f_1 := f$ we find inductively $f_i \in \sK_{C_i}^*$ 
such that it holds 
$\bigbt{\dv(f_{i-1})_{C_{i-1}}} = \bigbt{\dv(f_{i})_{C_{i}}} \in \ACHm{X}{\mdl}$ 
for $i = 2, \ldots, d-1$. 
Hence with $g := f_{d-1}$ this shows 
$\bigbt{\dv(f)_{C}} = \bigbt{\dv(g)_{Z}}$. 
\ePf 


\vspace{\vs} 
We fix a projective embedding of $X$, 
so we can talk about the degree of subschemes. 
For the computation of 
$\ACHm{X}{\mdl}  = \sum_C {\iot{C}}_*\ACHm{\Crv}{\mdl_{\Crvgst}}$, 
the following moving lemma allows us 
to restrict the index set of the sum to the set of curves of bounded degree: 

\bLem[Moving] 
\label{moving}
Let $C$ be a complete intersection curve in $X$, 
not necessarily reduced ore pure-dimensional, 
that meets $\mdl_X$ properly. 
Let $f \in \sK_C^*$ such that $f \in \sO_{C,\pntt}^*$\, 
for all $\pntt \in \pn{C \isec \mdl_X} \cup \st{\textrm{generic points of $C$}}$. 
There exist complete intersection curves $Z_i$ in $X$ 
of degree $\leq \degbnd := \deg X \cdot \pn{\deg\mdl_X}^{\dim X - 1}$ 
and elements $\pn{Z_i,g_i} \in \ZLm{X}{0 \mdl}$ 
such that 
\[ \bigbt{\dv(f)_C} = \biggbt{\sum_i \dv(g_i)_{Z_i}} \in \ACHm{X}{\mdl} \laurin 
\] 
\eLem 

\bPf 
Let \,$\emb: X \inj \Prj^m$\, be the fixed embedding into a projective space. 
We can extend $\mdl_X$ to a divisor that comes from $\Prj^m$; 
more precisely, we find an effective divisor $\mdll$ on $\Prj^m$ 
of degree \,$\deg \mdll \leq \deg \mdl_X$, 
meeting $X$ properly  
and such that $\mdl_X \leq \mdll \cut X$. 

According to Lemma \ref{degZf} 
we may assume that $f \in \sK_C^*$ is a representative of 
$\bigbt{\dv(f)_{\Ct}} \in \ACHm{\Ct}{\mdl_{\Ct}}$ such that 
\begin{align*} 
\deg \Zero\pn{f}_{\Ct} 
& \,\leq\, \deg\pn{\mdll \cut \Ct} \\ 
\vee \hspace{-1mm} \parallel \hspace{8mm} 
& \hspace{15mm} \parallel \\ 
\deg \Zero\pn{f}_{C} 
& \hspace{7mm} 
\deg\pn{\mdll \cut C} 
\,\leq\, \deg \mdl_X \cdot \deg C 
\end{align*} 
(for \,$\deg\pn{\mdll \cut \Ct} = \deg\pn{\mdll \cut C}$ see \cite[9.1]{BLR}), 
thus 
\[ \deg \Zero\pn{f}_{C} \,\leq\, \deg \mdl_X \cdot \deg C 
\laurin 
\] 

Let $\tha \in \sK_X^*$ be a lift of $f \in \sK_C^*$. 
For the divisor of zeroes of $f$ on $C$ we have the following relations 
\[ \deg \Zero\pn{f}_{C} \,=\, \deg \Zero\pn{\tha|_C} 
\,=\, \deg \bigpn{\Zero\pn{\tha} \cut C} 
\,=\, \frac{\deg \Zero\pn{\tha} \cdot \deg C}{\deg X} 
\laurin 
\] 
Hence 
\[\deg \Zero(\tha) \,\leq\, \deg \mdl_X \cdot \deg X 
\laurin 
\] 
In the same way we obtain for the divisor of poles 
\[ \deg \Pole(\tha) \,\leq\, \deg \mdl_X \cdot \deg X 
\laurin 
\] 
As 
\[ \dv(\tha) \,=\, \Zero(\tha) - \Pole(\tha) 
\laurink 
\] 
the result follows now by the reciprocity trick from Lemma \ref{tameLift}: 
\[ \bigbt{\dv(f)_C} = \bigbt{\partial\bigpn{\dv(\tha),\gam|_{\dv(\tha)}}} 
= \bigbt{\partial\bigpn{\Zero(\tha),\gam|_{\Zero(\tha)}} 
    - \partial\bigpn{\Pole(\tha),\gam|_{\Pole(\tha)}}} 
\] 
if $X$ is a surface. 
If $X$ is higher dimensional, 
we derive the assertion by induction over a regular sequence for $C$ 
and the surface case, as in  Step 2 of the proof of Lemma \ref{div_C=div_Z}. 
\ePf 

\bLem 
\label{degZf}
Let $C$ be a smooth projective curve over a field, 
$\mdl_{C}$ an effective divisor on $C$. 
Let \,$\morphism: C \ra \Prj^m$ be a morphism to a projective space, and 
$\mdll$ be an effective divisor on \,$\Prj^m$ meeting \,$\morphism\pn{C}$ properly 
such that \,$\mdl_{C} \leq \mdll \cut C$. 

For any \,$\alp \in \ACHm{C}{\mdl_{C}}$\, 
there is \,$f \in \sK_{C}^*$\, 
with \,$\bt{\dv(f)_{C}} = \alp \in \ACHm{C}{\mdl_{C}}$\, 
such that \,$\deg \Zero(f) \leq \deg\pn{\mdll \cut C}$. 
\eLem 

\bPf 
The effective divisor $\mdll$ on $\Prj^m$ is defined by a homogeneous 
element $\eps \in k[X_0, \ldots, X_m]$, and 
$\mdll \cut C = \Zero(\eps|_C) = \sum_{\pnt \in \vrt{\mdll \cut C}} n_{\pnt} \bt{\pnt}$. 
Now consider the following surjections ($\see$Lemma \ref{CHm-surj}): 
\[ \ACHm{C}{\mdl_{C}} \lsurleft \ACHm{C}{\mdll \cut C} 
	\lsurleft \frac{\bigcap_{\pnt \in \vrt{\mdll \cut C}} \sO_{C,\pnt}^*} 
       {\bigcap_{\pnt \in \vrt{\mdll \cut C}} \pn{1 + \fm_{C,\pnt}^{n_{\pnt}}}} 
= \biggpn{\frac{\sO_C}{(\eps|_C)}}^* 
\laurin 
\] 
Then for any $\alp \in \ACHm{C}{\mdl_{C}}$ 
there is $f \in \bigcap_{\pnt \in \vrt{\mdll \cut C}} \sO_{C,\pnt}^*$ 
with $\bt{\dv(f)_{C}} = \alp \in \ACHm{C}{\mdl_{C}}$ 
and such that 
\,$\deg \Zero(f) \leq \deg \Zero(\eps|_C) = \deg\pn{\mdll \cut C}$. 
\ePf 

\bPrp 
\label{ACH bounded intersection}
The affine part $\ACHm{X}{\mdl}$ 
of the relative Chow group with modulus of $X$ 
is generated by the groups $\ACHm{C}{\mdl_{\Cgst}}$ 
for complete intersection (=: ci) curves $C$ in $X$ 
of degree \,$\deg C \leq \degbnd := \deg X \cdot \pn{\deg\mdl}^{\dim X - 1}$: 
\[ \ACHm{X}{\mdl}  \;=  
\sum_{\substack{C \subset X \,\textrm{ci curve} \\ C \iscp \mdl \\ 
\deg C \leq \degbnd}} {\iot{C}}_* \ACHm{C}{\mdl_{\Cgst}} 
\laurin 
\] 
\ePrp 

\bPf 
Let $C$ be a curve in $X$. 
For every point $\pntt \in C$ 
there is an open affine neighbourhood $U \subset X$ of $\pntt$ 
and a complete intersection curve $C_U \subset U$ containing $C \cap U$. 
Write $C_U = \bigcap_i H_{i,U}$ for some hypersurfaces $H_{i,U}$ on $U$. 
Let $H_i$ be the closure of $H_{i,U}$ in $X$. 
Then $C_{(U)} := \bigcap_i H_i$ is a complete intersection curve on $X$. 
Here $C_U$ and $C_{(U)}$ need not be pure-dimensional. 
As $C$ is quasi-compact, there is a finite cover $\fU$ of $C$ 
by open affine subsets $U \subset X$. 
Then $C \subset \bigcup_{U \in \fU} C_{(U)}$, 
hence 
\[ \ACHm{C}{\mdl} \;\subset\; \sum_{U \in \fU} {\iot{C}}_* \ACHm{C_{(U)}}{\mdl} 
\laurin 
\] 
The rest of the statement follows directly from Lemma \ref{moving}. 
\ePf 

\vspace{\vs} 

\bNot 
\label{truncated curve}
For $e \in \Nat$ set 
\[ \Trc{X}{\mdl}{\degbnd} = 
\bigst{C \tms_X \bt{\mdl_X} \;\big|\; 
	C \subset X \textrm{ curve}, \,C \iscp \mdl_X, \,\deg C \leq \degbnd} 
\] 
i.e.\ the set of a proper intersections of 
curves of degree $\leq \degbnd$ 
with the subscheme $\bt{\mdl_X}$ of $\mdl_X$ in $X$. 
This is a set of finite subschemes in $X$, 
which we will regard as \emph{truncated curves with respect to $\mdl$}. 
\eNot 

\bDef 
For $e \in \Nat$ and $T \in \Trc{X}{\mdl}{\degbnd}$ 
let 
\[\fC_{T} = 
	\bigst{ C \;\big|\; C \subset X \textrm{ curve such that }
					C \tms_X \bt{\mdl_X} = T 
	}
	\laurin 
\]
Set 
\[ {\iot{T}}_* \ACHm{T}{\mdl} := {\iot{C}}_* \ACHm{C}{\mdl} 
	\hspace{5mm} \textrm{for some } C \in \fC_T 
\]
as a subgroup of $\ACHm{X}{\mdl}$, 
where $\iot{C}: C \inj X$ is the embedding. 
\eDef

The definition of ${\iot{T}}_* \ACHm{T}{\mdl}$ 
does not depend on the choice of $C \in \fC_T$, 
according to Lemma \ref{div_C=div_Z}. 
Then Proposition \ref{ACH bounded intersection} yields 

\bThm 
\label{ACH generated by truncations}
The group $\ACHm{X}{\mdl}$ is generated as follows: 
\[ \ACHm{X}{\mdl}  \;=  
\sum_{T \in \Trc{X}{\mdl}{\degbnd}} 
{\iot{T}}_* \ACHm{T}{\mdl_{\CTgst}}  
\laurin 
\] 
where \;$\degbnd := \deg X \cdot \pn{\deg\mdl_X}^{\dim X - 1}$. 
\eThm

\newpage 

\section{Modulus of a Rational Map to a Torsor} 
\label{Modulus}

Let $X$ be a normal proper variety over a perfect field $k$. 

\bDef 
\label{DefMod}
Let $\phe:X\dra P$ be a rational map from $X$ 
to a torsor $P$ under a smooth connected algebraic group $G$. 
Let $L$ be the affine part of $G$ 
in the canonical decomposition according to the theorem of Chevalley \cite{C}, 
and $\Upt$ the unipotent part of $L$ 
in the decomposition according to \cite[XVII, 7.2.1]{SGA3}. 
The \emph{modulus of $\phe$} from \cite[\S3]{KR10} 
is the effective divisor 
\[ \modu\pn{\phe} = \sum_{\codim(\pnt)=1} \modu_{\pnt}(\phe) \;\mdl_{\pnt} 
\] 
where $\pnt$ ranges over all points of codimension 1 in $X$, 
and $\mdl_{\pnt}$ is the prime divisor associated to $\pnt$, 
and $\modu_{\pnt}(\phe)$ is defined as follows. 

First assume $k$ is algebraically closed. 
Then we identify the torsor $P$ 
with the group $G$ acting on it. 
For each $\pnt\in X$ of codimension 1, 
the canonical map 
$L\pn{\sK_{X,\pnt}}/L\pn{\sO_{X,\pnt}} \lra 
  G\pn{\sK_{X,\pnt}}/G\pn{\sO_{X,\pnt}}$ 
is bijective, see \cite[No.~3.2]{KR10}. 
Take an element $l_{\pnt} \in L\pn{\sK_{X,\pnt}}$ 
whose image in $G\pn{\sK_{X,\pnt}}/G\pn{\sO_{X,\pnt}}$ 
coincides with the class of $\phe \in G\pn{\sK_{X,\pnt}}$. 
If $\chr(k) = 0$, let $(u_{\pnt,i})_{1\leq i\leq s}$ be the image of $l_{\pnt}$ 
in $\Ga(\sK_{X,\pnt})^s$ under $L\to\Upt\cong (\Ga)^s$. 
If $\chr(k) = p > 0$, 
let $(u_{\pnt,i})_{1\leq i\leq s}$  be the image of $l_{\pnt}$ in 
$\Witt_r(\sK_{X,\pnt})^s$ under  $L \to \Upt \subset \pn{\Witt_r}^s$. 
Let $\fil_n$, $\filf_n$ and $\filfb_n$ be the filtrations of $\Ga$ resp.\ $\Witt_r$ 
from Definition \ref{filtration_Witt}. 

In the \emph{log} theory, the multiplicities of the modulus are defined by 
\[ \modu_{\pnt}(\phe) = \left\{ 
   \begin{array}{ll} 
      0 & \;\textrm{if } \; \phe \in G\pn{\sO_{X,\pnt}} \\ 
      1+\max\bigst{\nty_{\pnt}(u_{\pnt,i}) \;|\; 1\leq i\leq s} 
         & \;\textrm{if } \; \phe \notin G\pn{\sO_{X,\pnt}} 
   \end{array} 
   \right. 
\] 
where for $u \in \Ga\pn{\sK_{X,\pnt}}$ resp. $\Witt_r\pn{\sK_{X,\pnt}}$ 
we set 
\[ \nty_{\pnt}(u) = \left\{ 
   \begin{array}{ll} 
      \min\bigst{n\in\Nat \;|\; u\in \fil_{n} \Ga\pn{\sK_{X,\pnt}}} 
          & \;\textrm{if} \; \chr(k) = 0 \phantom{F_{F_{F_{F_{F}}}}}\\ 
      \min\bigst{n\in\Nat \;|\; u\in \filf_{n} \Witt_{r}\pn{\sK_{X,\pnt}}} 
         & \;\textrm{if} \; \chr(k) = p > 0 \laurin 
   \end{array} 
   \right. 
\] 
In the \emph{non-log} theory, the multiplicities of the modulus are defined by 
\[ \modu_{\pnt}(\phe) = \left\{ 
   \begin{array}{ll} 
      0 & \;\textrm{if } \; \phe \in G\pn{\sO_{X,\pnt}} \\ 
      \max\bigst{\ntyb_{\pnt}(u_{\pnt,i}) \;|\; 1\leq i\leq s} 
       & \;\textrm{if } \; \phe \notin G\pn{\sO_{X,\pnt}} 
   \end{array} 
   \right. 
\] 
where for $u \in \Ga\pn{\sK_{X,\pnt}}$ resp. $\Witt_r\pn{\sK_{X,\pnt}}$ 
we set 
\[ \ntyb_{\pnt}(u) = \left\{ 
   \begin{array}{ll} 
      \min\bigst{n \in \Nat\setminus\st{0} \;|\; u\in \fil_{n-1} \Ga\pn{\sK_{X,\pnt}}} 
          & \;\textrm{if} \; \chr(k) = 0 \phantom{F_{F_{F_{F_{F}}}}}\\ 
      \min\bigst{n \in \Nat\setminus\st{0} \;|\; u\in \filfb_{n} \Witt_{r}\pn{\sK_{X,\pnt}}} 
          & \;\textrm{if} \; \chr(k) = p > 0 \laurin 
   \end{array} 
   \right. 
\] 
Note that in the case of characteristic 0, 
log and non-log modulus coincide. 

The definition of the multiplicities $\modu_{\pnt}(\phe)$ 
makes use of an isomorphism $\Upt \cong (\Ga)^s$ 
resp.\ of an embedding $\Upt \subset \pn{\Witt_r}^s$, 
which are not canonical. 
However, the definition is independent of the choice of those isomorphisms, 
see \cite[Thm.~3.3]{KR10}. 

For an arbitrary perfect base field $k$ we obtain \,$\modu\pn{\phe}$\, 
by means of a Galois descent from \,$\modu\pn{\phe \tens_k \clfld}$, 
where $\clfld$ is an algebraic closure of $k$, 
see \cite[No.~3.4]{KR10}. 
\eDef

\newpage 

\section{Albanese Variety with Modulus} 
\label{subsec:AlbMod}

Let $X$ be a suitable projective variety over a perfect field $\fld$, 
let $\mdl$ be a modulus for $X$, as in Section \ref{ChowMod}. 
If the base field $\fld$ is finite or algebraically closed 
(we only deal with these cases in this paper), 
every $\fld$-torsor admits a $\fld$-rational point, 
so in this case one can identify a $\fld$-torsor $\Trr{1}$ 
with the $\fld$-group $\Trr{0}$ acting on it. 
Therefore we will omit the superscript in this case. 

\bDef 
\label{CatMr} 
(Cf.\ \cite[Def.~2.8]{Ru13}.)
A \emph{category} $\Mr$ \emph{of rational maps from} $X$
\emph{to torsors} is a category satisfying the following conditions: 
The objects of $\Mr$ are rational maps $\phe: X \dra P$, 
where $P$ is a torsor for a smooth connected commutative algebraic group. 
The morphisms of $\Mr$ between two objects $\phe: X \dra P$
and $\psi: X \dra Q$ are given by the set of 
those homomorphisms of torsors 
$h:P \lra Q$ 
such that $h \circ \phe = \psi$. 
\eDef 

\bDef 
\label{DefAlb(X,D)_2}
(Cf.\ \cite[Def.~3.12]{Ru13}.)
The category $\Mrm{X}{\mdl}$  
is the category of those rational maps 
$\phe$ from $X$ to torsors 
such that \;$\modu\lrpn{\phe} \leq \mdl$. 
The universal object of $\Mrm{X}{\mdl}$ (if it exists) 
is denoted by $\albbm{1}{X}{\mdl}: X \dra \Albbm{1}{X}{\mdl}$ 
and called the \emph{Albanese of $X$ of modulus $\mdl$}. 
\eDef 

\bDef 
\label{DefDivf}
(Cf.\ \cite[No.~2.1]{Ru13}.)
Let $\Divf_X: \Algk \lra \Ab$ be the group-functor of relative Cartier divisors 
(from the category of $k$-algebras to the category of abelian groups). 
For a finite $k$-ring $R$ it holds 
\[ \Divf_X\pn{R} = 
   \Gam\bigpn{X\tens R, \lrpn{\sK_X \tens_k R}^*/\lrpn{\sO_X \tens_k R}^*}  
   \laurin 
\] 
\eDef 

\bDef 
\label{DefFm(X,D)}
(Cf.\ \cite[Def.~3.14]{Ru13}.)
Let $\Fm{X}{\mdl}$ denote the formal subgroup of $\Divf_X$ 
characterized as follows:  
Set 
\[ \lrpn{\Fm{X}{\mdl}}_{\et} = 
   \left\{ B\in\Divf_{X}(k) \;\big|\; \Supp(B) \subset \Supp(\mdl_X) \right\} 
   \laurin 
\] 
\noindent 
Set $\mdlinf = \lrpn{\mdl-\mdl_{\red}} \cut X$. 
For $\chr(k)=0$ set 
\[ \lrpn{\Fm{X}{\mdl}}_{\inf} = \exp \lrpn{ \Gac \tens_k 
   \Gam\bigpn{\sO_{X}\pn{\mdlinf} \big/ \sO_{X}}} 
\] 
where \;$\exp$\; is the exponential map, 
and $\Gac$ is the completion of $\Ga$ at $0$. 

\noindent 
For $\chr(k)=p>0$ in the \emph{log} theory set 
\[ \pn{\Fm{X}{\mdl}}_{\inf} = 
   \Expah \Biggpn{ \sum_{r > 0} \Wcfl{r} \tens_{\Witt(k)} 
   \Gam\Bigpn{\filf_{\mdlinf} \Witt_r(\sK_{X}) \Big/ \Witt_r(\sO_{X})}} 
   \laurin 
\] 
\noindent 
For $\chr(k)=p>0$ in the \emph{non-log} theory set 
\[ \pn{\Fm{X}{\mdl}}_{\inf} = 
   \Expah \Biggpn{ \sum_{r > 0} \Wcfl{r} \tens_{\Witt(k)} 
   \Gam\Bigpn{\filfb_{\mdl} \Witt_r(\sK_{X}) \Big/ \Witt_r(\sO_{X})}} 
\] 
where 
$\Expah$ denotes the Artin-Hasse exponential 
(as in Subsection \ref{Terminology}), 
$\Wittc$ is the completion of the Witt group $\Witt$ at $0$ 
from \cite[Example 1.11]{Ru13}, 
and $\filf_{\mdl} \Witt_r(\sK_X)$ resp.\ $\filfb_{\mdl} \Witt_r(\sK_X)$ 
are the filtrations from \cite[Definition 3.2 resp.\ 3.8]{Ru13}. 

Note that 
\;$\filfb_{\mdl} \Witt_r(\sK_{X}) = \filf_{\mdlinf} \Witt_r(\sK_{X})$\; 
if $\fld(\pnt)$ is perfect for every generic point $\pnt$ of $\vrt{\mdl}$, 
so in particular if $\chr(k)=0$, 
$\seecite$\cite[4.7 (2)]{KR10}. 
\eDef 

\bDef 
Let $\Fmor{X}{\mdl} = \Fm{X}{\mdl} \tms_{\Picf_X} \Picorf{X}$ 
be the part of $\Fm{X}{\mdl}$ that is 
mapped to the Picard variety $\Picor{X}$ 
under the class map \,$\Divf_X \lra \Picf_X$. 

\eDef 

\bPrp 
\label{Fm(X,D)algebraic}
The formal groups $\Fm{X}{\mdl}$ and $\Fmor{X}{\mdl}$ are dual-algebraic. 
\ePrp 

\bPf
In the log theory:  
See \cite[Proposition 3.15]{Ru13}. \\ 
In the non-log theory: 
As \,$\filfb_{\mdl} \Witt_r(\sK_{X}) \subset \filf_{\mdl} \Witt_r(\sK_{X})$\, 
by \cite[Def.~3.8]{Ru13}, 
$\pn{\Fm{X}{\mdl}}_{\textrm{non-log}}$ is a formal subgroup of 
$\pn{\Fm{X}{\mdl+\mdl_{\red}}}_{\textrm{log}}$. 
Since the latter formal group is dual-algebraic by the proof in the log theory, 
the same is true for $\pn{\Fm{X}{\mdl}}_{\textrm{non-log}}$ 
by \cite[Lemma 1.17]{Ru13}. 
Then also $\pn{\Fmor{X}{\mdl}}_{\textrm{non-log}}$ is dual-algebraic, 
by the same argument as in the log theory.  
\ePf 

\bLem 
\label{mod(phe)-im(trafo_phe)}
Let $\phe:X\dra G$ be a rational map from $X$ to a smooth connected 
algebraic group $G$. Then the following conditions are equivalent: 

\begin{tabular}{rl} 
\hspace{\kukwad} {\rm (i)} & $\modu\pn{\phe} \leq \mdl_X$\laurink \\ 
\hspace{\kukwad} {\rm (ii)} & $\im\pn{\trafo_{\phe}} \subset \Fm{X}{\mdl}$\laurin 
\end{tabular} 
\eLem 

\bPf
In the log theory: See \cite[Lemma 3.16]{Ru13}. \\ 
In the non-log theory: 
The proof is similar to the one in the log theory, 
but we have to replace $(ii)_{\inf,\pnt}$ by 

\begin{tabular}{ll} 
\hspace{\kukwad} ${\rm (ii)}_{\inf,\pnt}$ & 
If $n_{\pnt} = 0$, 
then $u_{\pnt,i} \in \Ga\pn{\sO_{X,\pnt}}$ resp.\ $\Witt_r\pn{\sO_{X,\pnt}}$
for $1\leq i\leq a$. \\ 
 & If $n_{\pnt} > 0$, 
then $\ntyb_{\pnt}\pn{u_{\pnt,i}} \leq n_{\pnt}$. 
\end{tabular} 

\hfill 
\ePf

\bThm 
\label{AlbMod-construction_2}
The Albanese $\albbm{1}{X}{\mdl}: X \dra \Albbm{1}{X}{\mdl}$ 
of \,$X$ of modulus $\mdl$ exists. 
The algebraic group $\Albbm{0}{X}{\mdl}$ acting on $\Albbm{1}{X}{\mdl}$ 
is dual (as a 1-motive) to the 1-motive 
$\bigbt{\Fmor{X}{\mdl} \lra \Picor{X}}$. 
(For 1-motives and duality see \cite[No.\ 1.2]{Ru13}.) 
\eThm 

\bPf 
See \cite[Thm.s~3.18 and 3.19]{Ru13}. 
The proof makes use of Lemma \ref{mod(phe)-im(trafo_phe)}, 
the non-log version of which is added above. 
\ePf 

\bPrp 
\label{phe^Sig}
A rational map $\phe: X \dra \Trr{1}$ from $X$ to a torsor $\Trr{1}$
under an algebraic group $\Trr{0}$ 
induces a map 
\begin{eqnarray*}
  \phe^{\Sig}: \Z_0(U)^0 & \lra & \Trr{0}(k) \\ 
  \sum \nbr_i \, \pntt_i & \lmt & \sum \nbr_i \, \phe(\pntt_i), 
  \hspace{8mm} 
  l_i \in \Zint, \; \pntt_i \in U \;\textrm{closed} 
\end{eqnarray*} 
where $U \subset X$ is the set on which $\phe$ is defined. 
\ePrp 

\bPf 
Let $\clfld$ be an algebraic closure of $\fld$, 
denote $\ol{U} := U \tens_{\fld} \clfld$. 
We have a map $\ZOo{U} \lra \ZOo{\ol{U}} \lra \Trr{0}\pn{\clfld}$ 
as defined above. 
The image $\phe^{\Sig}\pn{z}$ of any $z \in \ZOo{U}$, 
as a point of $\Trr{0}\pn{\clfld}$, 
is geometrically irreducible, and defined over $\fld$, 
hence an element of $\Trr{0}\pn{\fld}$. 
\ePf 

\bDef 
\label{MrCH(X,D)_2} 
Let $\MrCHm{X}{\mdl}$ be the category of those 
rational maps to torsors $\phe:X \dra \Trr{1}$ 
whose associated map $\phe^{\Sig}: \Z_0(U)^0 \lra \Trr{0}(k)$ 
factors through a homomorphism of groups $\CHmo{X}{\mdl}\lra \Trr{0}(k)$. 
%

We refer to the objects of $\MrCHm{X}{\mdl}$ as 
rational maps from $X$ to torsors 
\emph{factoring through $\CHm{X}{\mdl}$}. 
(Cf.\ \cite[Def.~3.28]{Ru13}.) 
The universal object of $\MrCHm{X}{\mdl}$ (if it exists) 
is denoted by $\albbch{1}{X}{\mdl}: X \dra \Albbch{1}{X}{\mdl}$. 
The algebraic group $\Albbch{0}{X}{\mdl}$ acting on $\Albbch{1}{X}{\mdl}$
is called the \emph{universal quotient of $\CHmo{X}{\mdl}$}. 
\eDef 

\bRmk 
\label{extend to X-D}
A rational map $\phe: X \dra G$ factoring through $\CHm{X}{\mdl}$ 
always extends to a morphism $X\setminus\mdl \lra G$, 
if the base field is algebraically closed. 
This follows from \cite[Thm.~3.29]{Ru13}, 
since in \cite[Lem.~3.30]{Ru13} 
one can replace condition (ii) 
by the corresponding condition for 
``all curves $C$ intersecting $X \setminus U$ properly'', 
where $U$ is an open dense subset of $X$ that $\phe$ is defined on. 
Then a rational map $\phe: X \dra G$ factoring through \,$\CHm{X}{\mdl}$\, 
as defined in Definition \ref{MrCH(X,D)_2} is of modulus $\leq \mdl$, 
which shows that $\phe$ is defined on $X \setminus \mdl$. 
\eRmk 

\bDef 
\label{DefFch}
Let $\Fch{X}{\mdl}$ be the formal subgroup of $\Divf_X$ defined as follows: 
\[  \Fch{X}{\mdl} 
    \,=\, \bigcap_{C} \lrpn{ \llul\cut\Ct }^{-1} \Fm{\Ct}{\mdl_{\Ct}} 
\] 
where $C$ ranges over all curves in $X$ that intersect $\mdl$ properly, 
$\Decf_{X,C}$ 
is the subfunctor of $\Divf_{X}$ 
consisting of those relative Cartier divisors on $X$ 
whose support (see \cite[Def.~2.2]{Ru13}) 
meets $C$ properly 
and $\lul \cut \Ct: \Decf_{X,C} \lra \Divf_{\Ct}$ 
is the pull-back of relative Cartier divisors from $X$ to $\Ct$. 
\eDef 

\bPrp 
\label{AlbCH-construction}
Assume the base field $\fld$ is algebraically closed. 
The universal quotient $\Albch{X}{\mdl}$ of $\CHmo{X}{\mdl}$, 
if it exists (as an algebraic group), 
is dual (in the sense of 1-motives) to the 1-motive 
$\bigbt{\Fchor{X}{\mdl} \lra \Picor{X}}$. 
\ePrp 

\bPf 
It is sufficient to show that the category $\MrCHm{X}{\mdl}$ 
is equal to the category $\Mr_{\Fch{X}{\mdl}}$ 
of those rational maps that induce a transformation to $\Fch{X}{\mdl}$ 
(see \cite[Thm.~2.16 and Rmk.~2.18]{Ru13}). 
For this aim we show that for a morphism $\phe: X\setminus\mdl \lra G$ 
the following conditions are equivalent: 

\begin{tabular}{rll} 
(i) & $\phe\lrpn{\dv(f)_C} = 0$ \hspace{10mm} 
    & $\forall \,(C,f) \in \Rlnpm{X}{\mdl} \laurink$ \\ 
(ii) & $\pn{\phe|_{\Ct},f}_{\pnt} = 0$ 
     & $\forall \,(C,f) \in \Rlnpm{X}{\mdl}, 
          \hspace{1mm} \forall \,\pnt \in \vrt{\mdl_{\Ct}} \laurink$ \\ 
(iii) & $\modu\bigpn{\phe|_{\Ct}} \leq \mdl_{\Ct}$ 
      & $\forall \,C \subset X, \,C \iscp \mdl \laurink$ \\ 
(iv) & $\im\bigpn{\trafo_{\phe|_{\Ct}}} \subset \Fm{\Ct}{\mdl_{\Ct}}$ 
      & $\forall \,C \subset X, \,C \iscp \mdl \laurink$ \\ 
(v) & $\im\pn{\trafo_{\phe}} \cut \Ct \subset \Fm{\Ct}{\mdl_{\Ct}}$ 
     & $\forall \,C \subset X, \,C \iscp \mdl \laurink$ \\ 
(vi) & $\im\pn{\trafo_{\phe}} \subset \Fch{X}{\mdl} \laurin$ & 
\end{tabular} 

\noindent 
(i)$\Llra$(ii) see \cite[III, \S1]{S59}\laurink 
(ii)$\Llra$(iii) see \cite[No.~6.1--3]{KR10}\laurink 
(iii)$\Llra$(iv) is \cite[Lem.~3.16]{Ru13}\laurink 
(iv)$\Llra$(v) is evident\laurink 
(v)$\Llra$(vi) by Definition \ref{DefFch}\laurin 

\hfill 
\ePf 

\bThm 
\label{ThmCHoMod_2}
Assume the base field $\fld$ is algebraically closed. 
The category $\Mrm{X}{\mdl}$ of rational maps of modulus $\leq\mdl$ 
is equal to the category $\MrCHm{X}{\mdl}$ of rational maps 
factoring through $\CHm{X}{\mdl}$. 
In particular we have \;$\Albm{X}{\mdl} = \Albch{X}{\mdl}$\; 
and \;$\Fmor{X}{\mdl} = \Fchor{X}{\mdl}$. 
\eThm 

\bPf 
See \cite[Thm.~3.29]{Ru13}. 
The proof makes use of Lemma \ref{Mdl_restr_to_curves}, 
the non-log version of which is proven below. 
\ePf 

\bLem 
\label{Mdl_restr_to_curves}
Assume $\fld$ is algebraically closed. 
Let $\phe:X\dra G$ be a rational map from $X$ 
to a smooth connected algebraic group $G$, 
regular away from a closed proper subset $S \varsubsetneq X$. 
Then the following conditions are equivalent: 

\begin{tabular}{rl} 
\hspace{\kukwad} {\rm (i)} & $\modu\pn{\phe} \leq \mdl_X$\laurink \\ 
\hspace{\kukwad} {\rm (ii)} & $\modu\bigpn{\phe|_{\Ct}} \leq \mdl_{\Ct} 
\hspace{15mm} \forall \,C \subset X, \,C \iscp S$\laurin 
\end{tabular}
\eLem 

\bPf
In the log theory: See \cite[Lemma 3.30]{Ru13}. \\ 
In the non-log theory: 

(i)$\Lra$(ii)
Let $C$ be a curve in $X$ intersecting $S$ properly. 
By Definition \ref{DefMod} of the modulus, 
it is easy to see that 
\,$\modu\pn{\phe} \leq \mdl_X$\, 
implies 
\;$\modu\bigpn{\phe|_{\Ct}} \leq \mdl_{\Ct}$. 

\medskip 
(ii)$\Lra$(i) 
We show $\neg {\rm (i)} \Lra \neg {\rm (ii)}$. 
Suppose $\mdll := \modu\pn{\phe_X} \not \leq \mdl_X$. 
Then there is a point $\pnt \in \Supp\pn{\mdll}$ of codimension 1 in $X$ 
such that $\multy_{\pnt}(\mdll) > \multy_{\pnt}\pn{\mdl_X}$, 
where $\multy_{\pnt}$ is the multiplicity at $\pnt$. 
If $\chr(k) = 0$, it is easy to see that a general curve $C$ in $X$ 
intersecting $\Es_{\pnt}$ in a regular closed point $x$ satisfies 
\[ \modu_{x}\bigpn{\phe|_{\Ct}} 
= \multy_{x}\bigpn{\pn{\mdll-\mdll_{\red}}\cut \Ct} + 1 
> \multy_{x}\pn{\mdl \cut \Ct} = \multy_{x}\pn{\mdl_{\Ct}} 
\laurin 
\] 
Therefore we suppose $\chr(k)=p>0$. 
Using the notation of Definition \ref{DefMod}, let 
\,$\pn{u_{\pnt,i}}_{1\leq i\leq a} 
\in \Witt_r\pn{\sK_{X,\pnt}}^a$\, 
be a representative of the unipotent part of the class of 
$\phe\in G\pn{\sK_{X,\pnt}}$ in $G\pn{\sK_{X,\pnt}}/G\pn{\sO_{X,\pnt}} 
= L\pn{\sK_{X,\pnt}}/L\pn{\sO_{X,\pnt}}$. 
Then $\modu_{\pnt}\pn{\phe} = \ntyb_{\pnt}\pn{u_{\pnt,i}}$ 
for some $1 \leq i \leq a$. 
Write $\ntyb_{\pnt}\pn{u_{\pnt,i}} =: n + 1$. 
Let $t \in \fm_{X,\pnt}$ be a uniformizer at $\pnt$. 
If \,$\qfill_{(n+1)\pnt} \pn{u_{\pnt,i}} \in \Qfillb_{(n+1)\pnt}$ 
($\see$\cite[Def.~3.8]{Ru13}), 
let \,$\sum_{\nu} \Frob^{\nu} \tens \,\oma_{\nu} \tens t^{-(n+1)} \in 
k[\Frob] \tens_k \Kot_{\Es,\pnt} \tens_{\sO_{X,\pnt}} \fm_{X,\pnt}^{-(n+1)}$ 
be a representative of 
$\qfill_{(n+1)\pnt} \pn{u_{\pnt,i}}$; 
if not, let \;$\sum_{\nu} \Frob^{\nu} \tens \,\oma_{\nu} \tens t^{-n} \in 
k[\Frob] \tens_k \Kot_{X,\pnt}\pn{\log\pnt} \tens_{\sO_{X,\pnt}} \fm_{X,\pnt}^{-n}$
be a representative of $\qfill_{n\pnt} \pn{u_{\pnt,i}} \in \Qfill_{n\pnt}$ 
(see \cite[Def.~3.7]{Ru13}). 
Choose a regular closed point $x \in \Es_{\pnt}$ 
such that $t$ is a local equation for $\Es_{\pnt}$ at $x$
and $\oma_{\nu}$ is regular and $\neq 0$ at $x$ for some $\nu$. 
We may assume that $\dim X = 2$ 
via cutting down by hyperplanes through $x$ transversal to $\Es$.  
Let $s \in \fm_{X,x}$ be a local parameter at $x$ 
that gives a uniformizer of $\sO_{\Es,x}$. 
Define a curve $C_{e}$ locally around $x$ 
by the equation \;$t = s^{e}$ \; for $e \geq 1$. 
As $\mdll-\mdll_{\red}$ is locally defined by the equation $t^{n} = 0$, 
then 
\[ \multy_{x}\bigpn{\mdl_{\wt{C_{e}}}} 
\,\leq\, \multy_{x}\bigpn{\pn{\mdll-\mdll_{\red}}\cut C_{e}} 
\,=\, \dim_k \frac{\sO_{X,x}}{\pn{t^{n}, t - s^e}} \,=\, ne 
\laurin 
\]
If \,$\qfill_{(n+1)\pnt} \pn{u_{\pnt,i}} \in \Qfillb_{(n+1)\pnt}$, 
we write \;$\oma_{\nu} = g\, \der s$\; (for some $\nu$) 
with $g \in \sO_{X,\pnt}$ 
and we may assume $g(x) \neq 0$. 
The restriction of \,$t^{-(n+1)} \oma_{\nu}$\, to \,$C_{e}$\, is 
\begin{eqnarray*} 
t^{-(n+1)} \oma_{\nu}|_{C_{e}} 
& = & s^{-(n+1)e} g\, \der s \\ 
& = & s^{1-(n+1)e} g\, \der\log s \laurin 
\end{eqnarray*} 
Otherwise, if \,$\qfill_{n\pnt} \pn{u_{\pnt,i}} \in 
 \Qfill_{n\pnt} \setminus \Qfillb_{n\pnt}$, 
we have \,$\oma_{\nu} = g\, \der s + h\, \der\log t$\, 
with $g,h \in \sO_{X,\pnt}$ and $h(x) \neq 0$. 
The restriction of \,$t^{-n} \oma_{\nu}$\, to \,$C_{e}$\, is 
\begin{eqnarray*} 
t^{-n} \oma_{\nu}|_{C_{e}} 
& = & s^{-ne} g\, \der s +  s^{-ne} h\, \der\log s^e \\ 
& = & s^{1-ne} g\, \der\log s + e \,s^{-ne} h\, \der\log s \laurin 
\end{eqnarray*} 
Then the class of \,$t^{-(n+1)} \oma_{\nu}|_{C_{e}}$\, 
resp.\ \,$t^{-n} \oma_{\nu}|_{C_{e}}$\, 
is \,$\neq 0$\, in 
\[ 
\left\{ 
\begin{array}{l}
\Kot_{C_{e},x}(\log x) \tens_{\sO_{C_{e},x}} 
\fm_{C_{e},x}^{1-(n+1)e} \big/ \fm_{C_{e},x}^{2-(n+1)e} 
\hspace{4mm} \textrm{ if }\; \qfill_{(n+1)\pnt} \pn{u_{\pnt,i}} \in \Qfillb_{(n+1)\pnt 
 _{\phantom{F_{F_{F}}}}} \\ 
\Kot_{C_{e},x}(\log x) \tens_{\sO_{C_{e},x}} 
\fm_{C_{e},x}^{-ne} \big/ \fm_{C_{e},x}^{1-ne}
\hspace{6mm} \textrm{ if }\; \qfill_{n\pnt} \pn{u_{\pnt,i}} \in 
 \Qfill_{n\pnt} \setminus \Qfillb_{n\pnt} \textrm{ and } p \nmid e 
\laurin 
\end{array}
\right. 
\] 
\cite[Lemma 3.10]{Ru13} 
assures that the modulus of $\phe|_{C_{e}}$ 
is computed from the restriction 
of $\qfil_{(n+1)\pnt} \pn{u_{\pnt,i}}$ 
resp.\ $\qfil_{n\pnt} \pn{u_{\pnt,i}}$ 
to $C_{e}$, 
for $e$ sufficiently large such that \,$ne-1 > \lfloor ne/p\rfloor$ 
(this is satisfied for $e > 2$). 
Thus we have 
\begin{align*} 
\modu_{x} \pn{\phe|_{C_{e}}} 
= \ntyb_{x} \pn{u_{\pnt,i}|_{C_{e}}} 
& = \left\{ 
\begin{array}{cl}
(n+1)e & \textrm{ if }\; \qfill_{(n+1)\pnt} \pn{u_{\pnt,i}} \in \Qfillb_{(n+1)\pnt
_{\phantom{F_{F_{F}}}}} \\ 
ne + 1 & \textrm{ if }\; \qfill_{n\pnt} \pn{u_{\pnt,i}} \in 
 \Qfill_{n\pnt} \setminus \Qfillb_{n\pnt} \textrm{ and } p \nmid e 
\end{array}
\right. \\ 
& >\; ne 
\;\geq\; \multy_{x}\bigpn{\mdl_{\wt{C_{e}}}} 
\end{align*} 
which shows that 
\;$\modu_{x} \pn{\phe|_{C_{e}}} \not \leq \mdl_{\wt{C_{e}}}$. 
\ePf 

\bCor 
\label{factorCHm}
For any perfect base field $\fld$, 
a rational map $\phe: X \dra \Trr{1}$ 
to a $\fld$-torsor $\Trr{1}$ for an algebraic $\fld$-group $\Trr{0}$ 
of modulus $\leq \mdl$ 
factors through $\CHm{X}{\mdl}$ in the sense of Definition \ref{MrCH(X,D)_2}. 
\eCor 

\bPf 
The implications (i) $\Lra$ (iii) $\Lra$ (ii) 
in the proof of \cite[Thm.~3.29]{Ru13} 
hold for any perfect base field $\fld$, 
for this direction we do not need the assumption 
that $\fld$ is algebraically closed. 
\ePf 

\bCor 
\label{Alb=QuotCH}
For any perfect base field $\fld$, 
the group of $\fld$-rational points of 
the Albanese group with modulus $\Albbm{0}{X}{\mdl}(\fld)$ 
is a quotient of $\CHmo{X}{\mdl}$, 
and this formation is compatible with the universal map 
$\albbm{1}{X}{\mdl}$ 
of $\Mrm{X}{\mdl}$. 
\eCor 

\bPf 
As the universal map 
\[ \albbm{1}{X}{\mdl}: X \lra \Albbm{1}{X}{\mdl} 
\] 
of $\Mrm{X}{\mdl}$ 
is of modulus $\leq \mdl$, 
it factors through $\CHm{X}{\mdl}$ 
by Corollary \ref{factorCHm}. 
Now $\Albbm{1}{X}{\mdl}$, as the universal object of $\Mrm{X}{\mdl}$, 
is generated by $X$, 
hence the induced map 
\[ \ZOm{X}{\mdl} \lra \Albbm{0}{X}{\mdl}\pn{\fld} 
\] 
is surjective 
(cf.\ \cite[No.\ 1]{SeAlba}). 
Thus there is an epimorphism 
\[ \CHmo{X}{\mdl} \lsur \Albbm{0}{X}{\mdl}\pn{\fld} 
\laurin 
\] 
\hfill 
\ePf 

\bCor 
\label{Fm=Fm^CH}
The equality \;$\Fmor{X}{\mdl} = \Fchor{X}{\mdl}$\; holds 
for any perfect base field. 
\eCor 

\bPf 
Both formal groups are compatible with Galois descent. 
\ePf

\vspace{\vs} 
\noindent 
Now we consider an arbitrary projective variety $X$, 
not necessarily suitable. 

\bThm 
\label{Alb^CH--X-sing}
Assume $\fld$ is algebraically closed. 
The universal quotient $\Albch{X}{\mdl}$ of \,$\CHmo{X}{\mdl}$ 
exists for every projective variety $X$ over $\fld$. 
Moreover, $\Albch{X}{\mdl}$ is covariant functorial in $X$ 
and contravariant functorial in $\mdl$. 
\eThm 

\bPf 
Let $\sig: \Xt \ra X$ be the suitabilization of $X$. 
Since $\sig$ is birational, 
any rational map $\phe: X \dra G$ from $X$ to an algebraic group $G$ 
can be viewed 
as a rational map $\phe \circ \sig: \Xt \dra G$ from $\Xt$ to $G$. 
As $k$ is algebraically closed, 
the push-forward of cycles $\ZOm{\Xt}{\mdl} \lra \ZOm{X}{\mdl}$ is surjective, 
and the induced map $\CHm{\Xt}{\mdl} \lra \CHm{X}{\mdl}$ is well defined 
by Proposition \ref{FunctorCHoMod}, making 
$\CHm{X}{\mdl}$ a quotient of $\CHm{\Xt}{\mdl}$. 
In this way, 
the category $\MrCHm{X}{\mdl}$ becomes a subcategory of $\MrCHm{\Xt}{\mdl}$. 
Therefore, if $\phe: X \dra G$ is an object of $\MrCHm{X}{\mdl}$ 
and $G$ is generated by $X$, 
then $G$ is a quotient of the universal object $\Albch{\Xt}{\mdl}$ 
of $\MrCHm{\Xt}{\mdl}$. 
Thus, according to \cite[No.~2, Cor.\ of Thm.~2]{SeAlba}, 
the category $\MrCHm{X}{\mdl}$ admits a universal object $\Albch{X}{\mdl}$. 
The functoriality of $\Albch{X}{\mdl}$ follows easily from the 
defining universal property of $\Albch{X}{\mdl}$. 
\ePf

\newpage 

\section{Duality of Divisors and 0-Cycles} 
\label{AbelJacobiMap}

Let $X$ be a projective variety over an algebraically closed field $k$. 
Let $\mdl$ be a modulus for $X$, 
as in Section \ref{ChowMod}. 
In this section we construct a pairing 
\[ \pair{\llul,\lull}_{X,\mdl}: \Fmor{X}{\mdl} \tms \ACHm{X}{\mdl} \lra \Gm 
\] 
($\see$Proposition \ref{Pairing_X}), 
that will become crucial for the Skeleton Theorem ($\see$Theorem \ref{SkelThm}) 
and for the Roitman Theorem with Modulus ($\see$Theorem \ref{RoitmanThm}).

\subsection{Duality on Curves} 
\label{sub:Duality}

Let $\Crv$ be a smooth proper 
curve over a field $k$, 
and let $\mdl = \sum_{\pnt\in S} n_{\pnt} \bt{\pnt}$ be an effective divisor on $\Crv$, 
where $S$ is a finite set of closed points on $\Crv$ 
and $n_{\pnt}$ are integers $\geq 1$ for $\pnt\in S$. 

\bPnt 
\label{Curve-business}
If $k$ is a perfect field, 
the Albanese group with modulus $\Albbm{0}{\Crv}{\mdl}$ 
of $\Crv$ of modulus $\mdl$ 
coincides with the Jacobian group with modulus $\Jaccm{0}{\Crv}{\mdl}$ 
of Rosenlicht-Serre 
(see \cite[Thm.~3.25]{Ru13}). 
While the log and non-log versions of $\Albbm{0}{X}{\mdl}$ 
do not coincide in general, 
they do in the curve case over a perfect field. 
This is due to the fact that the two filtrations of the Witt group from Definition \ref{filtration_Witt} 
coincide if the residue fields of codimension 1 points are perfect ($\seecite$\cite[4.7 (2)]{KR10}). 
The generalized Jacobian $\Jaccm{0}{\Crv}{\mdl}$ is an extension 
\bDpl  0 \lra \Llm{0}{\Crv}{\mdl} \lra \Jaccm{0}{\Crv}{\mdl} \lra \Jaccc{0}{\Crv} \lra 0 
\eDpl 
of the classical Jacobian $\Jaccc{0}{\Crv} \cong \Pic^0_{\Crv}$ of $\Crv$, 
an abelian variety,  
by an affine algebraic group $\Llm{0}{\Crv}{\mdl}$ 
(see \cite[V, \S3]{S59}). 
Taking $k$-valued points, this exact sequence becomes 
the short exact sequence from Remark \ref{kerCHm}: 
\begin{align*} 
\Jaccc{0}{\Crv}\pn{k} & \,=\, \CHOo{\Crv} \laurink \\ 
\Jaccm{0}{\Crv}{\mdl}\pn{k} & \,=\, \CHmo{\Crv}{\mdl} \laurink \\ 
\Llm{0}{\Crv}{\mdl}\pn{k} & \,=\, \ACHm{\Crv}{\mdl} 
 \,=\,  \frac{\prod_{\pnt\in S} \sO_{\Crv,\pnt}^*}
            {\extfld^* \tms \prod_{\pnt\in S}\bigpn{1+\fm_\pnt^{n_\pnt}}} 
\,\laurink 
\end{align*} 
where $\extfld := \H^0\pn{C,\sO_C}$, 
where 
$\fm_{\pnt}$ denotes the maximal ideal at $\pnt\in\Crv$ 
(see \cite[I, No.~1]{S59}). 
Our duality construction of $\Albbm{0}{\Crv}{\mdl}$ 
($\see$Theorem \ref{AlbMod-construction_2}) 
yields that $\Jaccm{0}{\Crv}{\mdl}$ is dual to 
$\bigbt{\Fmo{\Crv}{\mdl} \lra \Pic^0_{\Crv}}$, 
in particular $\Llm{0}{\Crv}{\mdl}$ is Cartier dual to $\Fmo{\Crv}{\mdl}$. 
\ePnt 


\bDef 
\label{Def_perfPair}
Let $\fmlG$ be a formal group and $A$ be an abstract abelian group. 
A \emph{pairing} 
\bDpl  
\fmlG \tms A \lra \Gm 
\eDpl 
is a collection of bi-homomorphisms 
\;$  
\fmlG(R) \tms A \lra \Gm(R) 
\laurin 
$\; 
Such a pairing is called \emph{perfect}, 
if it induces an isomorphism 
\;$ 
A \cong \Homabk\pn{\fmlG, \Gm}
\laurin 
$\; 

If $\fmlG$ is dual to an affine group-variety and 
$A = L\pn{k}$ denotes the group of $k$-valued points of an affine group-variety $L$ 
over an algebraically closed field $k$, 
the condition above is equivalent to saying that 
the pairing induces an isomorphism of sheaves 
\;$ \fmlG \cong \Homfabk\pn{L, \Gm} 
\laurink 
$ 
according to Cartier duality (see e.g.\ 
\cite[Thm.~1.5]{Ru13}) 
and the fact that 
varieties over an algebraically closed base field $k$ 
are determined by its $k$-valued points. 
\eDef 

\bDef 
\label{Pairing_curve,S}
Let $\Divf_{\Crv}^{S}$ denote the subfunctor of $\Divf_{\Crv}$ 
of relative Cartier divisors with support in $S$. 
Then we define a pairing 
\begin{align*}
\pair{\llul,\lull}_{C,S}: \Divsfc{\Crv}{S} \;\tms\; \Rlnm{\Crv}{0 \mdl} \; & \lra \; \Gm \\ 
\bigpair{\sD,\dv(f)}_{\Crv,S} & = \,\prod_{\pnt \in S} \,\pn{\sD,f}_{\pnt} 
\end{align*} 
induced by the local symbols 
\,$\pn{\llul,\lull}_{\pnt} : G\pn{\sK_{\Crv}} \tms \Oc_{\Crv,\pnt}^* \lra G 
$\, 
at $\pnt \in S$, where $G$ is an algebraic group, 
as in \cite[III, \S1]{S59}. 
For a relative divisor $\sD \in \Divf_{\Crv}(R)$, 
where $R$ is a finite $k$-ring, 
choose a local section 
of $\sD$ in a neighbourhood $U_{\pnt}$ of $\pnt$ in $\Crv$, 
which is an element of $\Gamma\bigpn{U\tens R,\sO_{\Crv\tens R}^*} 
= \Gm\bigpn{\sO_{\Curv}(U)\tens_k R}$, 
where $U = U_{\pnt} \setminus S$. 
This defines a rational map \,$\phi^{\sD,\pnt}: \Crv \dra \WeilRes{R}{k} \Gm$\, 
from $C$ to the Weil restriction of $\GmS{R}$ from $R$ to $k$. 
Then we define 
\[ \pn{\sD,f}_{\pnt} := \pn{\phi^{\sD,\pnt},f}_{\pnt} 
\laurin 
\]
Here the local symbol $\pn{\phi^{\sD,\pnt},f}_{\pnt}$ is independent 
of the choice of the local section 
of $\sD$ in $U_{\pnt}$: 
another local section would lead to a rational map 
$\psi^{\sD,\pnt} = \phi^{\sD,\pnt} \cdot \gam$, 
where $\gam: \Crv \dra \WeilRes{R}{k} \Gm$ is regular at $\pnt$. 
Then, due to \cite[Prop.\ 3.15]{Ru08}, 
$\pn{\phi^{\sD,\pnt},f}_{\pnt}$ and $\pn{\psi^{\sD,\pnt},f}_{\pnt}$ 
differ only by $\pn{\gam,f}_{\pnt} =  \gam(\pnt)^{\val_{\pnt}(f)} = 1$, 
as $f \in \Oc_{\Crv,\pnt}^*$ is regular at $\pnt$. 
\eDef 

\bPnt 
\label{pairingFml_curve}
Let $\Divf_{\Crv}^{S,0}$ be the subfunctor of 
$\Divf_{\Crv}^{0} = \Divf_{\Crv} \tms_{\Picf_X} \Picf_X^0$ 
of relative Cartier divisors of degree 0 with support in $S$. 
Let $\sD \in \Divf_{\Crv}^{S,0}(R)$, $R$ a finite $k$-ring. 
Let $\Gp\pn{\sD} \in \Extabk\pn{\Alba{\Crv}, \Lin_R} \cong \Pic_{\Alba{\Crv}}^0(R) 
\cong \Pic_{\Crv}^0(R)$ 
be the algebraic group corresponding to $\sO_{\Crv \tens R}(\sD)$, 
where $\Lin_R := \WeilRes{R}{k} \Gm$  
is the Weil restriction of $\GmS{R}$ from $R$ to $k$. 
(Here we used \cite[Prop.\ 1.13]{Ru13}.) 
The canonical 1-section of $\sO_{\Crv \tens R}(\sD)$ induces a rational map 
\,$\rmp{\sD}: \Crv \dra \Gp\pn{\sD}$, 
which is regular away from $S$. 
If $h \in \sO_{\Crv,\pnt}^*$ for some closed point $\pnt \in \Crv$, 
then by \cite[Lemma~3.16]{Ru08} 
the local symbol $\pn{\rmp{\sD},h}_{\pnt}$ 
lies in the fibre of $\Gp\pn{\sD}$ over \,$0 \in \Alba{\Crv}$, which is $\Lin_R$, 
and $\pn{\rmp{\sD},h}_{\pnt} = \pn{\sD,h}_{\pnt}$. 
Then for $f \in \sO_{C,S}^*$ we have 
\begin{eqnarray*} 
\bigpair{\sD, \dv\pn{f}}_{\Crv,S} 
 & = & \prod_{\pnt \in S} \bigpn{\rmp{\sD}, f}_{\pnt} 
 \; = \; \prod_{\pnt \in \vrt{\dv\pn{f}}} 
            \bigpn{\rmp{\sD}, f}_{\pnt}^{-1} \\ 
 & = & \prod_{\pnt \in \vrt{\dv\pn{f}}} 
            \rmp{\sD}(\pnt)^{-\val_{\pnt}\pn{f}} 
 \; = \; \rmp{\sD} \bigpn{\dv\pn{f}}^{-1} 
\laurin 
\end{eqnarray*}
\ePnt 

\bPrp 
\label{duality_curve}
In the curve case over a perfect field, the Cartier duality \\ 
\;$\Fmo{\Crv}{\mdl} = \bigpn{\Llm{0}{\Crv}{\mdl}}^{\vee}$\; 
is expressed by the perfect pairing 
\begin{align*} 
\pair{\llul,\lull}_{C,\mdl_{\Cgst}}: \Fmo{\Crv}{\mdl} \;\tms\; \ACHm{\Crv}{\mdl} \; & \lra \; \Gm \\ 
\Bigpair{\sD,\bigbt{\dv(f)}}_{\Crv,\mdl} & = \,\prod_{\pnt \in S} \,\pn{\sD,f}_{\pnt} 
\end{align*} 
induced by the pairing from Definition \ref{Pairing_curve,S}. 
\ePrp 

\bPf 
By means of a Galois descent we can reduce to the case 
of an algebraically closed base field. 
$\ACHm{\Crv}{\mdl}$ is compatible with Galois descent 
since it is given by the affine part of the Jacobian with modulus. 

Let $\sD \in \Divf_{\Crv}^{S,0}(R)$, with $R$ a finite $k$-ring. 
Suppose $\pn{\sD,f}_{\pnt} = 0$ for all $\pnt \in S$, $f \in \Oc_{\Crv,\pnt}^*$. 
Then by \Point \ref{pairingFml_curve} we have 
\;$\rmp{\sD}\bigpn{\dv(f)} = 0$\; for all \,$\dv(f) \in \Rlnm{\Crv}{0 \mdl}$, 
i.e.\ \,$\pn{\rmp{\sD}}^{\Sig}: \Z_0(\Crv \setminus S)^0 \lra \Gm(R)$\, 
factors through 
$\Z_0(\Crv \setminus \mdl)^0 \big/ \Rlnm{\Crv}{0 \mdl} = \CHOo{\Crv}$ 
due to Bloch's moving ($\see$Remark \ref{Bloch}). 
This implies $\rmp{\sD}$ factors through $\Jac_{\Crv}$, 
which is an abelian variety, 
therefore extends to a morphism defined on the whole of $\Crv$. 
This means that $\rmp{\sD}$ is a global section of $\sO_{\Crv}(\sD)$, 
which is only possible if \,$\sD = 0$. 

Thus the local symbol induces a monomorphism of formal groups 
\bDpl 
\Divsfc{\Crv}{S,0} 
\,\linj\, \dsum_{\pnt \in S} \Homfabk\Bigpn{\Oc_{\Crv,\pnt}^*,\Gm} 
  = \Homfabk\Bigpn{\prod_{\pnt \in S} \Oc_{\Crv,\pnt}^*,\Gm} 
  \laurin 
\eDpl 
By construction of $\Fmo{\Crv}{\mdl}$ ($\see$Definition \ref{DefFm(X,D)}) 
one sees, according to \cite[\S6, Proposition 6.4 (3)]{KR10} 
or by the equivalence of conditions (ii) and (iv) 
in the proof of Proposition \ref{AlbCH-construction}, 
that $\Fmo{\Crv}{\mdl}$ annihilates 
$k^* \tms \prod_{\pnt\in S}\bigpn{1+\fm_\pnt^{n_\pnt}}$ 
with respect to the local symbol: 
\[ \Fmo{\Crv}{\mdl} = \Ann \bigpn{\Rlnm{\Crv}{\mdl}} \laurin 
\]  
Thus the local symbol induces a monomorphism $\lam\sig:$ 
\[ 
   \Fmo{\Crv}{\mdl} \,\lra \, 
   \Homfabk\Biggpn{\frac{\prod_{\pnt\in S} \sO_{C,\pnt}^*}
            {k^* \tms \prod_{\pnt\in S}\bigpn{1+\fm_\pnt^{n_\pnt}}} , \Gm} 
   = \Homfabk\bigpn{\Lm{\Crv}{\mdl},\Gm}   
\laurin 
\] 
It remains to show surjectivity of $\lam\sig$. 
Using $\Fmo{C}{\mdl} \cong \Lm{C}{\mdl}^{\vee}$, 
this map yields a monomorphism $\Lm{C}{\mdl}^{\vee} \linj \Lm{C}{\mdl}^{\vee}$. 
By duality, this corresponds to an epimorphism $\Lm{C}{\mdl} \lsur \Lm{C}{\mdl}$. 
As the unipotent part of $\Lm{C}{\mdl}$ has a filtration with quotients of type $\Ga$, 
such an epimorphism is necessarily an isomorphism on the unipotent part of 
$\Lm{C}{\mdl}$. 
This shows that the map $\lam\sig$ is an an isomorphism 
on infinitesimal parts. 
The surjectivity of $\lam\sig$ on \'etale parts is easily checked by hand, 
using the explicit computation from \Point \ref{PairingExplicit}. 
\ePf 

\smallskip 

\bPnt[Explicit computation of the duality pairing from Proposition \ref{duality_curve}] 
\label{PairingExplicit}
It is sufficient to describe the log version, 
since the non-log version coincides with the log version 
in the curve case (cf.\ \Point \ref{Curve-business}). 
The Cartier duality pairing 
\;$\pair{\llul,\lull}_{C,\mdl}: \Fmo{C}{\mdl} \tms \Llm{0}{C}{\mdl} \lra \Gm$\; 
is expressed by local symbols. 
Assume first 
that the base field $k$ is algebraically closed. 
For $\sD \in \Fmo{C}{\mdl}(R) \subset \Divf_{\Crv}(R)$, $R$ a finite $k$-ring, 
choose a local section of $\sD$ in a neighbourhood of $\pnt$ in $\Crv$ 
for every $\pnt \in S$. 
Then locally around $\pnt$ the pairing $\pair{\llul,\lull}_{C,\mdl}$ 
can be computed by explicit formulas for the local symbols 
\;$\pn{\llul,\lull}_{C,\pnt}: G\pn{\sK_{C,\pnt}} \tms \sK_{C,\pnt}^* \lra \Gm\pn{k}$\; 
for $G = \Gm, \Ga, \Witt_r$ 
via the following commutative diagrams: 

For the \'etale part of $\Fmo{C}{\mdl}$, 
as it is determined by its $k$-valued points, 
it suffices to consider 
\[ \xymatrix{ 
\pn{\Fmo{C}{\mdl}}_{\et}\pn{\fld} 
\tms \ACHm{C}{\mdl} \ar[rr]^-{\pair{\llul,\lull}_{C,\mdl}} 
&& \Gm\pn{\fld} \\ 
\Gm\pn{\sK_C} \tms \bigcap_{\pnt} \sO_{C,\pnt}^* 
\ar[rr]_-{\pn{\llul,\lull}_{C,\pnt}} 
\ar[u]^{\dv \tms \res} 
&& \Gm\pn{\fld} \ar@{=}[u] 
} 
\] 
\hspace{9mm}
where \;$\pn{f,\alp}_{C,\pnt} = \pn{-1}^{\val_{\pnt}\pn{f} \val_{\pnt}\pn{\alp}} 
\frac{f^{\val_{\pnt}\pn{\alp}}}{\alp^{\val_{\pnt}\pn{f}}}(\pnt)$\; 
for $f \in \Gm\pn{\sK_{C,\pnt}}$, $\alp \in \sK_{C,\pnt}^*$. 
\medskip 

Due to associativity of the local symbol, 
i.e.\ $\pn{\psi \circ \phe, \lull}_{C,\pnt} = \psi \circ \pn{\phe,\lull}_{C,\pnt}$, 
we obtain for the infinitesimal part of $\Fmo{C}{\mdl}$: 

\noindent 
$\pn{\chr(\fld) = 0}$ 
\[ \xymatrix{ 
\pn{\Fmo{C}{\mdl}}_{\inf}\pn{R} \tms \ACHm{C}{\mdl} \ar[rr]^-{\pair{\llul,\lull}_{C,\mdl}} 
&& \Gm\pn{R} \\ 
\bigpn{\Gac\pn{R} \tens \sO_C\pn{\mdl - \mdl_{\red}}} 
\tms \bigcap_{\pnt} \sO_{C,\pnt}^* 
\ar[rr]_-{\id \tens \pn{\llul,\lull}_{C,\pnt}} 
\ar[u]^{\pn{\exp \circ \pn{\id \tens \res}} \tms \res} 
&& \Gac\pn{R} \tens \Ga\pn{\fld} \ar[u]_{\exp} 
} 
\] 
\hspace{9mm}
where \;$\pn{f,\alp}_{C,\pnt} = \Res_{\pnt}\bigpn{f \,\der \log\pn{\alp}}$\; 
for $f \in \Ga\pn{\sK_{C,\pnt}}$, $\alp \in \sK_{C,\pnt}^*$. 
\medskip 

\noindent 
$\pn{\chr(\fld) = p}$ 
\[ \xymatrix{ 
\pn{\Fmo{C}{\mdl}}_{\inf}\pn{R} \tms \ACHm{C}{\mdl} \ar[rr]^-{\pair{\llul,\lull}_{C,\mdl}} 
&& \Gm\pn{R} \\ 
\bigpn{\Wcfl{r}\pn{R} \tens \filf_{\mdl-\mdl_{\red}} \Witt_r\pn{\sK_{C,\pnt}}} 
\tms \bigcap_{\pnt} \sO_{C,\pnt}^* 
\ar[rr]_-{\id \tens \pn{\llul,\lull}_{C,\pnt}} 
\ar[u]^{\pn{\Expah \circ \pn{\id \tens \res}} \tms \res} 
&& \Wcfl{r}\pn{R} \tens \Witt_r\pn{\fld} \ar[u]_{\Expah} 
} 
\] 
\hspace{9mm} 
where 
\[ \pn{f,\alp}_{C,\pnt} = 
\Frob^{1-r} \Res_{\pnt}\bigpn{\Phi_{r-1}(f) \,\der \log\pn{\wt{\alp}}} 
\] 
for $f \in \Witt_r(\sK_{C,\pnt})$, $\alp \in \Kc_{C,\pnt}^* = k((t))^*$, 
$\wt{\alp} \in \Witt_r(k)[[t]][t^{-1}]$ a lift of $\alp$, 
\begin{align*} 
\Frob: & \;\Witt_r(k) \lra \Witt_r(k) \;, 
& 
\pn{w_0, \ldots, w_{r-1}} & \lmt \pn{w_0^p, \ldots, w_{r-1}^p} \;\laurink \\ 
\Phi_j: & \;\Witt_r(\sK_{C,\pnt}) \lra  \Witt_r(k)[[t]][t^{-1}] \;, 
& 
\pn{f_0, \ldots, f_{r-1}} & \lmt \sum_{0 \leq i \leq j} p^i \wt{f_i}^{p^{j-i}} \;\laurink \\ 
\Res: & \;\Omega^1_{\Witt_r(k)[[t]][t^{-1}]} \lra \Witt_r(k) \;, 
& 
\sum_i a_i \,t^i\,\der \log(t) & \lmt a_0 \;\laurink 
\end{align*} 
see \cite[No.~6.5]{KR10} 
for a description of the local symbol map to $\Witt_r$. 

For any perfect base field $k$ let $\clfld$ be an algebraic closure. 
If $\rmpC: C \dra G$ is a rational map to a $k$-group scheme, 
$\rmpC \tens \clfld: C \tens \clfld \dra G \tens \clfld$ the base changed map, 
$f \in \sK_C$ 
and $\pntt \in C$ a closed point, 
then 
\begin{eqnarray*} 
\pn{\rmpC, f}_{\pntt} 
& = & \sum_{C(\clfld) \ni \pnt \ra \pntt} \pn{\rmpC \tens \clfld, f \tens \clfld}_{\pnt} \\ 
& = & \Trace_{\pntt / \strl(\pntt)} \wt{\pn{\rmpC,f}_{\pntt}} 
\end{eqnarray*} 
with $\strl: C \lra \Spec(k)$ the structural morphism, 
$\Trace_{\Grnd / \Base}: G(\Grnd) \lra G(\Base)$ the trace map 
induced by a finite flat morphism $\Grnd \ra \Base$ 
($\seecite$\cite[XVII, (6.3.13.2)]{SGA4}) 
and $\wt{\pn{\rmpC,\lull}_{\pntt}}: \sK_C \lra G\bigpn{k(\pntt)}$ 
the map defined by the same formula as 
$\pn{\rmpC \tens \clfld, \lul \tens \clfld}_{\pnt}$ 
over an algebraic closure $\clfld$ of $k$, for $\pnt$ over $\pntt$. 
\ePnt 

\bPnt 
\label{Pairing_curve_lim}
Since the sheaf of relative Cartier divisors with support in $S$ 
can be expressed as an inductive limit of divisors with bounded pole order, 
we have \;$ \Divsfc{\Crv}{S,0} = 
\underset{\substack{D \\ \vrt{D} \subset S}} \varinjlim \,\Fmo{\Crv}{D}$, 
and let \;$\Rlnmc{\Crv}{S} = 
\underset{\substack{D \\ \vrt{D} \subset S}} \varprojlim \,\ACHm{\Crv}{D}$, 
where in both limits 
$D$ ranges over all effective divisors on $\Crv$ with support in $S$. 
Taking the limit of the pairings from Proposition \ref{duality_curve} 
yields a canonical perfect pairing  $\phantom{\wh{{D^0}^0}}$ 
\begin{align*} 
\bigpair{\llul,\lull}_{\Crv,S}: \; 
\Divsfc{\Crv}{S,0} \;\tms\; \Rlnmc{\Crv}{S} \; & \lra \; \Gm \laurin
\end{align*} 
In particular, 
this pairing yields an isomorphism 
\bDpl 
\Divsfc{\Crv}{S,0} \,\cong\,  
\Homfabk\Biggpn{\Bigpn{\prod_{\pnt \in S} \Oc_{\Crv,\pnt}^*} \Big/ k^*,\Gm} \laurin 
\eDpl 
\ePnt 

\vspace{\vs} 

Let $X$ be a suitable projective variety over a perfect field. 
Now we are looking at curves $C$ in $X$ 
that are not necessarily normal or irreducible. 

\bPnt 
\label{pairingFml}
Let $\sD \in \Divor{X}(R)$, where $R$ is a finite $k$-ring. 
As in \Point \ref{pairingFml_curve} 
let $\Gp\pn{\sD} \in \Extabk\pn{\Alba{X}, \Lin_R} \cong \Pic_{\Alba{X}}^0(R) 
\cong \Picor{X}(R)$ 
be the algebraic group corresponding to $\sO_{X \tens R}(\sD)$, 
where $\Lin_R := \WeilRes{R}{k} \Gm$ 
is the Weil restriction of $\GmS{R}$ from $R$ to $k$. 
The canonical 1-section of $\sO_{X \tens R}(\sD)$ induces a rational map 
$\rmp{\sD}: X \dra \Gp\pn{\sD}$, 
which is regular away from $\vrt{\sD}$. 
Let $C \subset X$ be a curve that meets $\vrt{\mdl}$ properly, 
and suppose $\sD \cut \Ct \in \Fmo{\Ct}{\mdl_{\Ct}}$. 
Then by \Point \ref{pairingFml_curve} we have 
for $f \in \sO_{C,\vrt{\mdl \cut C}}^{\,*}$: 
\begin{eqnarray*} 
\Bigpair{\sD \cut \Ct, \bigbt{\dv\pn{\wt{f}}}}_{\Ct,\mdl_{\Ct}} 
 & = & \prod_{\pnt \in \vrt{\mdl \cut \Ct}} \Bigpn{\rmp{\sD}\big|_{\Ct}, \wt{f}\,}_{\pnt} \\ 
 & = & \prod_{\pnt \in \vrt{\dv\pn{\wt{f}}}} 
            \Bigpn{\rmp{\sD}\big|_{\Ct}, \wt{f}\,}_{\pnt}^{-1} \\ 
 & = & \prod_{\pnt \in \vrt{\dv\pn{\wt{f}}}} 
            \rmp{\sD}\big|_{\Ct}(\pnt)^{-\val_{\pnt}\pn{\wt{f}}} \\ 
 & = & \prod_{\pntt \in \vrt{\dv\pn{f}}} 
            \rmp{\sD}\big|_{C}(\pntt)^{-\ord_{\pntt}\pn{f}} \\ 
 & = & \rmp{\sD} \bigpn{\dv\pn{f}_C}^{-1} 
\end{eqnarray*}
i.e.\ the expression depends only on the cycle 
$\dv\pn{f}_C \in \ZO{X}$ in $X$, 
not on the curve $C$. 
\ePnt 

\bLem 
\label{(D,div(f))}
Let $\pn{C,f}, \pn{Z,g} \in \ZLm{X}{0 \mdl}$ 
and suppose that these two 
pairs satisfy 
$\bt{\dv(f)_C} = \bt{\dv(g)_Z} \in \ACHm{X}{\mdl}$. 
Let $\sD \in \Divor{X}(R)$, where $R$ is a finite ring, 
with the property that $\sD \cut \Ct \in \Fmo{\Ct}{\mdl_{\Ct}}$ 
and $\sD \cut \Zt \in \Fmo{\Zt}{\mdl_{\Zt}}$. 
Then \;$\bigpair{\sD \cut \Ct, \dv\pn{\wt{f}}}_{\Ct,\mdl_{\Ct}}$ 
$= \bigpair{\sD \cut \Zt, \dv\pn{\wt{g}}}_{\Zt,\mdl_{\Zt}}$. 
\eLem

\bPf 
If $\bt{\dv(f)_C} = \bt{\dv(g)_Z} \in \ACHm{X}{\mdl}$, 
then there are $f' \in \sK_C^*$, $g' \in \sK_Z^*$ 
with $\bt{\dv(f')_C} = \bt{\dv(f)_C}$ resp.\ $\bt{\dv(g')_Z} = \bt{\dv(g)_Z}$ 
such that ${\dv(f')_C} = {\dv(g')_Z} \in \ZO{X}$. 
Then the assertion follows from \Point \ref{pairingFml}. 
\ePf 


\subsection{Representability of Chow Groups for Singular Curves}
\label{represent-Chow}

Let $k$ be an algebraically closed field. 

\bNot 
\label{group-algebraic}
We will say that an abstract group $\Gam$ is 
\emph{represented by a group scheme} 
(resp.\ \emph{algebraic group}) 
if there is a $k$-group scheme (resp.\ algebraic $k$-group) \,$G$ 
such that \,$\Gam = G\pn{k}$. 
(We will sometimes drop the ``represented by''.) 
\eNot 

This notation is justified by the fact that the functor 
\[ G \lmt G\pn{k} 
\] 
from the category of $k$-group schemes (resp.\ algebraic $k$-groups) 
to the category of abstract groups 
is an embedding, if $k$ is algebraically closed. 

\bThm 
\label{CHm-algGrp} 
Let $C$ be a projective curve over $k$, not necessarily irreducible or normal, 
and $\Ct$ its normalization. 
Let $\mdl$ be a modulus for $C$. 
Then \,$\CHmo{C}{\mdl}$ is represented by an algebraic group. 

Moreover, any morphism $\zeta: Z \ra C$ of projective curves over $k$ 
(with $\mdl$ as modulus) 
induces a homomorphism \,$\CHmo{Z}{\mdl} \lra \CHmo{C}{\mdl}$ 
which is represented by a homomorphism of algebraic groups. 
\eThm 

\bPf
By Proposition \ref{CHmSing-motivation} and taking degree 0 parts, we have 
\[ \CHmo{C}{\mdl} = \coker\Bigpn{\KZm{C}{\mdl} \overset{\gam}{\lra}  \CHmo{\Ct}{\mdl}} 
\laurin 
\] 
Here $\CHmo{\Ct}{\mdl}$ is an algebraic group due to \Point \ref{Curve-business}, 
whereas 
\[ \KZm{C}{\mdl_{\Cgst}} 
 =  \ker\Biggpn{\dsum_{\pnt \,\in\, \Ct \tms_C \Sing\pn{C,\mdl}} 
	\hspace{-5mm} \Zint \,\pnt 
	\lra \dsum_{\pntt \,\in\, \Sing\pn{C,\mdl}} 
	\hspace{-3mm} \Zint \,\pntt}
\]
and the natural map $\gam$ 
obviously define a group scheme and a homomorphism of those. 
The image of $\gam$ is a closed subgroup of an algebraic group, 
hence $\CHmo{C}{\mdl}$ arises as a quotient of algebraic groups. 
Since the category of algebraic groups is abelian, 
$\CHmo{C}{\mdl}$ is an algebraic group as well. 

Now assume that $\zeta: \Z \ra C$ is a morphism of projective curves. 
According to the description above, 
we obtain an induced map 
\[ \xymatrix{ 
	\CHmo{Z}{\mdl} \ar[r] & \CHmo{C}{\mdl} 
} 
\] 
as the cokernel of the map of two-term complexes 
\[ \xymatrix{ 
	\KZm{Z}{\mdl} \ar[r] \ar[d] & \KZm{C}{\mdl} \ar[d] \\ 
	\CHmo{\Zt}{\mdl} \ar[r] & \CHmo{\Ct}{\mdl} \laurin 
} 
\] 
The category of group schemes is a full subcategory of $\Abk$
(= category of sheaves of abelian groups). 
We have already seen 
that the groups $\CHmo{Y}{\mdl}$ and $\CHmo{X}{\mdl}$ 
are (in particular) group schemes, 
therefore it remains to show that $\CHmo{Z}{\mdl} \lra \CHmo{C}{\mdl}$ 
is a morphism in $\Abk$. 
The category $\Compl^{[-1,0]}\pn{\Abk}$ 
of two-term complexes in $\Abk$ is abelian. 
The lower row of the square is a homomorphism of group schemes 
due to the functoriality of Albanese varieties with modulus ($\see$\cite[Prop.\ 3.22]{Ru13}), whereas 
for the upper row this assertion is obvious due to their plain structure.  
Therefore both rows lie in $\Compl^{[-1,0]}\pn{\Abk}$, 
hence so does the quotient of these rows. 
\ePf


\subsection{Skeleton Relative Divisors} 
\label{skeleton}

Let $X$ be a projective variety, 
and assume that the base field $k$ is algebraically closed 
(unless stated otherwise). 
Let $\mdl$ be a modulus for $X$, 
as in Section \ref{ChowMod}. 
Let $C$ be a curve in $X$, not necessarily normal or irreducible, 
and $\Ct$ its normalization. 

\bPnt 
\label{CHm(C)_algGrp}
According to \Point \ref{Curve-business} 
the degree 0 part of the Chow group with modulus of the smooth curve $\Ct$ 
is an algebraic group: $\CHmo{\Ct}{\mdl} = \Jacm{\Ct}{\mdl}\pn{k}$. 
Then by Theorem \ref{CHm-algGrp} the same is true for $\CHmo{C}{\mdl}$, 
we denote the corresponding algebraic group by $\Jacm{C}{\mdl}$. 
\ePnt 

\bDef 
\label{Def_Fm_sing}
Let $\Lm{C}{\mdl}$ denote the affine part of $\Jacm{C}{\mdl}$. 
We define $\Emo{C}{\mdl_{\Cgst}}$ to be the Cartier dual of $\Lm{C}{\mdl}$: 
\[ \Emo{C}{\mdl} := \pn{\Lm{C}{\mdl}}^{\vee} 
\laurin 
\] 
\eDef 

\bPnt 
\label{dual_ACHm}
By definition, $\Lm{C}{\mdl}\pn{k}$ is the affine part of $\CHm{C}{\mdl_{\Cgst}}$, 
given by 
\[ \ACHm{C}{\mdl_{\Cgst}} 
:= \ker \Bigpn{\CHm{C}{\mdl_{\Cgst}} \lra \CHm{C}{0 \mdl_{\Cgst}}} 
\laurin 
\] 
Then the Cartier dual $\Emo{C}{\mdl}$ of $\Lm{C}{\mdl}$ 
is the annihilator of 
\[ \KACHm{C}{\mdl} 
:= \ker\Bigpn{\ACHm{\Ct}{\mdl} \lra \ACHm{C}{\mdl}} 
\] 
with respect to the pairing from Proposition \ref{duality_curve}. 
\ePnt 

\bRmk 
\label{Char_Em}
The description from \Point \ref{dual_ACHm} shows that 
it is possible to characterize 
the formal groups $\Emo{C}{\mdl_{\Cgst}}$ as follows: 
\[ \Emo{C}{\mdl_{\Cgst}} = \lrst{ \sD \in \Fmo{\Ct}{\mdl_{\Ct}} \left| 
\begin{array}{l} 
\bigpair{\sD,\kap}_{\Ct,\mdl_{\Ct}} = 0 \\ 
\forall \, \kap \in \ker\bigpn{\ACHm{\Ct}{\mdl_{\Ct}} \lra \ACHm{C}{\mdl_{\Cgst}}} 
\end{array}
\right.} 
\laurin 
\] 
In particular, if $C$ and $Z$ are curves in $X$, then 
\[ \Emo{C \cup Z}{\mdl_{\CcupZgst}} = 
\lrst{ 
\begin{array}{l} 
\hspace{0mm} \pn{\sD_C,\sD_Z} \in \\ 
\Emo{C}{\mdl_{\Cgst}} \tms \Emo{Z}{\mdl_{\Zgst}} 
\end{array}
\left|
\begin{array}{l} 
\bigpair{\sD_C,\alp}_{C,\mdl_{\Cgst}} = 
\bigpair{\sD_Z,\bet}_{Z,\mdl_{\Zgst}} \\ 
\forall \, \alp \in \ACHm{C}{\mdl_{\Cgst}}, \bet \in \ACHm{Z}{\mdl_{\Zgst}} \\ 
\textrm{with} \hspace{2mm}{\iot{C}}_*\alp = {\iot{Z}}_*\bet \in \ACHm{C \cup Z}{\mdl_{\CcupZgst}} 
\end{array}
\right.} 
\laurin 
\] 
\eRmk 


\bPrp 
\label{Functor_Em}
The affine algebraic group $\Lm{C}{\mdl}$ 
and the formal group $\Emo{C}{\mdl_{\Cgst}}$ from \Point \ref{dual_ACHm} 
are covariant resp.\ contravariant functorial in $C$: 
Let $\zeta: Z \lra C$ be a morphism of proper curves 
in $X$ that intersect $\mdl$ properly. 
Then $\zeta$ induces homomorphisms of algebraic groups resp. of formal groups 
\begin{align*} 
\zeta_*: \Lm{Z}{\mdl} & \lra\, \Lm{C}{\mdl} \\ 
              \Emo{Z}{\mdl_{\Zgst}} & \lla\, \Emo{C}{\mdl_{\Cgst}} : \zeta^* \laurin 
\end{align*} 
\ePrp 

\bPf 
Note that a morphism $\zeta: Z \lra C$ of curves in $X$ is necessarily injective. 
Then $\zeta$ 
yields a morphism 
$\wt{\zeta}: \Zt \lra \Ct$ between the normalizations of $Z$ and $C$. 
Due to the construction of $\Lm{C}{\mdl}$ 
the map $\zeta_*: \Lm{Z}{\mdl} \lra\, \Lm{C}{\mdl}$ 
is induced by the homomorphism 
$\wt{\zeta}_*: \Llm{0}{\Zt}{\mdl_{\Zt}} \lra\, \Llm{0}{\Ct}{\mdl_{\Ct}}$ 
via passing through the quotients, 
which in turn is the affine part of the functoriality map 
$\Alba{\wt{\zeta}}: \Albm{\Zt}{\mdl_{\Zt}} \lra \Albm{\Ct}{\mdl_{\Ct}}$. 
As $\Alba{\wt{\zeta}}$ is a homomorphism of algebraic groups, 
the same is true for $\zeta_*$. 
The homomorphism $\zeta^*: \Emo{C}{\mdl_{\Cgst}} \lra\, \Emo{Z}{\mdl_{\Zgst}}$ 
is obtained from $\zeta_*: \Lm{Z}{\mdl} \lra\, \Lm{C}{\mdl}$ 
via Cartier duality. 
\ePf 

\bDef
\label{Def_skeletonDiv}
Due to functoriality of $\Emo{C}{\mdl_{\Cgst}}$ 
we can form the projective limit $\varprojlim \Emo{C}{\mdl_{\Cgst}}$ 
ranging over all curves $C$ in $X$ intersecting $\vrt{\mdl_X}$ 
properly. 
An element of this projective limit is a compatible system 
of relative Cartier divisors on the various curves in $X$, 
which we will call a \emph{skeleton relative divisor} 
on $X$. 
\eDef 

\bLem 
\label{Fm-to-limEm}
Let $S$ be the support of a Cartier divisor $\geq \mdl_X$ on $X$. 
There is a canonical homomorphism 
\begin{eqnarray*}
\skelm{X}{\mdl,S}: \Fmor{X}{\mdl} & \lra & 
\varprojlim_{C \iscp S} \Emo{C}{\mdl_{\Cgst}} \\ 
                                                 \sD & \lmt & \pn{\sD \cut C}_C 
\end{eqnarray*}
given by pull-back of a relative Cartier divisor on \,$X$\, with support in \,$S$\, 
to \,$\Ct$\, for all curves \,$C \subset X$\, that meet \,$S$\, properly. 
Here we denote the image of \,$\sD \in \Fmor{X}{\mdl}$\, in \,$\Emo{C}{\mdl_{\Cgst}}$\, 
by \,$\sD \cut C$\, (instead of \,$\sD \cut \Ct$). 
\eLem 

\bPf 
According to Corollary \ref{Fm=Fm^CH} and Definition \ref{DefFch} we have 
\[  \Fmor{X}{\mdl} 
    \,=\, \bigcap_{C} \bigpn{ \llul\cut\Ct }^{-1} \Fmo{\Ct}{\mdl_{\Ct}} 
\] 
hence $\Fmor{X}{\mdl} \cut \Ct \subset \Fmo{\Ct}{\mdl_{\Ct}}$ 
for every curve $C$ in $X$ that meets $\vrt{\mdl}$ properly. 
The description from \Point \ref{pairingFml} 
of $\bigpair{\sD \cut \Ct, \alp}_{\Ct,\mdl_{\Ct}}$ 
for $\alp \in \ACHm{\Ct}{\mdl_{\Ct}}$ and $\sD \in \Fmor{X}{\mdl}$ 
says $\sD \cut \Ct$ annihilates 
\;$\ker\bigpn{\ACHm{\Ct}{\mdl_{\Ct}} \lra \ACHm{C}{\mdl_{\Cgst}}}$. 
Then the characterization of the group $\Emo{C}{\mdl_{\Cgst}}$ 
in Remark \ref{Char_Em} 
implies that in fact $\Fmor{X}{\mdl} \cut C \subset \Emo{C}{\mdl_{\Cgst}}$. 
Thus $\skelm{X}{\mdl,S}$ is a well-defined homomorphism. 
%
\ePf


\subsection{Duality in Higher Dimension} 
\label{duality_highDim}

\bPrp 
\label{Pairing_curve}
Let $C$ be a curve in $X$ which intersects $\vrt{\mdl_X}$ properly. 
In the following diagram, 
the upper row is a well-defined and perfect pairing 
(induced by the perfect pairing of the lower row 
from Proposition \ref{duality_curve}): 
\[ \xymatrix{ 
\Emo{C}{\mdl_{\Cgst}} \tms \ACHm{C}{\mdl_{\Cgst}} \ar[r]  
\ar@{}@<-15mm>[d]|(.34){\hspace{0mm} \cap} 
\ar@{}@<-15mm>[d]|(.46){\hspace{-1.70mm} |} 
\ar@{}@<-15mm>[d]|(.46){\hspace{+1.70mm} |} 
\ar@{}@<-15mm>[d]|(.62){\hspace{-1.70mm} |} 
\ar@{}@<-15mm>[d]|(.62){\hspace{+1.70mm} |} 
&  \Gm 
\phantom{\laurin} \\ 
\Fmo{\Ct}{\mdl_{\Ct}} \tms \ACHm{\Ct}{\mdl_{\Ct}} 
\ar[r] \ar@{->>}@<-6.0mm>[u] 
& \Gm \ar@<+1.8mm>@{=}[u] 
\laurin 
} 
\] 
\ePrp 

\bPf 
By construction ($\see$\Point \ref{dual_ACHm}), 
$\Emo{C}{\mdl_{\Cgst}}$ is exactly the annihilator of \\ 
$\ker\bigpn{\ACHm{\Ct}{\mdl_{\Ct}} \lra \ACHm{C}{\mdl_{\Cgst}}}$. 
\ePf 

\bLem 
\label{compSystem}
The pairings 
$\pair{\llul,\lull}_{C,\mdl_{\Cgst}}: \Emo{C}{\mdl_{\Cgst}} \tms \ACHm{C}{\mdl_{\Cgst}} \lra \Gm$, 
where $C$ ranges over curves in $X$ that intersect $\vrt{\mdl_X}$ properly, 
form a compatible system of perfect pairings: 
if $\zeta: Z \lra C$ is a $X$-morphism of curves in $X$ that intersect $\mdl_X$ properly, 
then for $\alp \in \ACHm{Z}{\mdl_{\Zgst}}$ and $\sD \in \Emo{C}{\mdl_{\Cgst}}$ it holds 
$\bigpair{\zeta^*\sD,\alp}_{Z,\mdl_{\Zgst}} = \bigpair{\sD,\zeta_*\alp}_{C,\mdl_{\Cgst}} \laurin 
$ 
\[ \xymatrix{ 
\Emo{C}{\mdl_{\Cgst}} \tms \ACHm{C}{\mdl_{\Cgst}} 
\ar[r] \ar@<-14.6mm>[d]_{\zeta^*} 
 & \Gm \\ 
\Emo{Z}{\mdl_{\Zgst}} \tms \ACHm{Z}{\mdl_{\Zgst}} \ar[r] \ar@<-5.8mm>[u]_{\zeta_*} 
 & \Gm \ar@<+1.2mm>@{=}[u] 
} 
\] 
\eLem 

\bPf 
A morphism of subvarieties of $X$ over $X$ is necessarily injective. 
Then the assertion follows directly 
from the definition of the pairings in Proposition \ref{Pairing_curve} 
and from the characterization of the formal groups $\Emo{C}{\mdl_{\Cgst}}$ 
in Remark \ref{Char_Em}. 
\ePf 

\bCor 
\label{Pairing_limit}
Taking the limit of the perfect pairings 
$\pair{\llul,\lull}_{C,\mdl_{\Cgst}}$ 
from Proposition \ref{Pairing_curve}, 
where $C$ ranges over all curves in $X$ that intersect $\vrt{\mdl_X}$ properly, 
yields a perfect pairing 
\[ 
\varprojlim_C \Emo{C}{\mdl_{\Cgst}} \tms \varinjlim_C \ACHm{C}{\mdl_{\Cgst}} \lra \Gm 
\laurin 
\] 
\eCor 

\bPrp 
\label{Pairing_X}
There is a canonical pairing 
\begin{align*} 
\pair{\llul,\lull}_{X,\mdl}: \Fmor{X}{\mdl} \;\tms\; \ACHm{X}{\mdl} \; & \lra \; \Gm \\ 
\Bigpair{\sD,\bigbt{\dv(f)_C}}_{X,\mdl} 
& = \Bigpair{\sD \cut C,\bigbt{\dv(f)}}_{C,\mdl_{\Cgst}} \\ 
\Bigpair{\sD,\bt{z}}_{X,\mdl} & = \; \rmp{\sD}\pn{z}^{-1} 
\end{align*} 
for $z \in \Rlnm{X}{0 \mdl} \subset \ZOm{X}{\mdl}$ 
and where $\rmp{\sD}$ is the rational map associated with $\sD$, 
i.e.\ induced by the 1-section of $\sO_{X \tens R}(\sD)$ 
(see \Point \ref{pairingFml}). 
\ePrp 

\bPf 
Since $\ACHm{X}{\mdl}  =  \varinjlim \ACHm{C}{\mdl_{\Cgst}}$ 
by Proposition \ref{CHm_indLim_algclosed} 
and due to the homomorphism 
\,$\skelm{X}{\mdl}: \Fmor{X}{\mdl} \lra \varprojlim_C \Emo{C}{\mdl_{\Cgst}}$ 
from Lemma \ref{Fm-to-limEm}, 
we can define the pairing via the following diagram 
\[ \xymatrix{ 
\Fmor{X}{\mdl} \hspace{2.5mm} \tms \hspace{2,5mm} \ACHm{X}{\mdl} 
\ar[r] \ar@<-17.5mm>[d]_{\skelm{X}{\mdl}} 
 & \Gm \ar@<-1.8mm>@{=}[d] \phantom{\laurin} \\ 
\underset{C} \varprojlim \,\Emo{C}{\mdl_{\Cgst}} \tms 
\underset{C} \varinjlim \ACHm{C}{\mdl_{\Cgst}} 
\ar[r] \ar@<-10mm>[u]^-{\wr} 
 & \Gm \laurin 
} 
\] 
The formula $\bigpair{\sD,\bt{z}}_{X,\mdl} = \rmp{\sD}\pn{z}^{-1}$ 
has been shown in \Point \ref{pairingFml}.  
\ePf 

\bRmk 
\label{restrict to curves}
Proposition \ref{Pairing_X} shows that 
the pairing $\pair{\llul,\lull}_{X,\mdl}$ 
can be defined by restriction to curves $C \subset X$, 
using the pairings $\pair{\llul,\lull}_{C,\mdl_{\Cgst}}$ 
from Proposition \ref{Pairing_curve}. 
Here 
the definition of $\bigpair{\sD,\alp}_{X,\mdl}$ for $\alp \in \ACHm{X}{\mdl}$ 
is independent of the choice of representatives $\pn{C_i,f_i} \in \Rlnpm{X}{0 \mdl}$ 
such that $\alp = \bigbt{\sum \dv\pn{f_i}_{C_i}}$, 
due to Lemma \ref{(D,div(f))}. 
\eRmk 


\bCor 
\label{Pairing_X_perfect}
The pairing 
$\pair{\llul,\lull}_{X,\mdl}: \Fmor{X}{\mdl} \tms \ACHm{X}{\mdl} \lra \Gm$ 
from Proposition \ref{Pairing_X} is perfect 
if and only if the skeleton map from Lemma \ref{Fm-to-limEm} 
\;$\skelm{X}{\mdl}: \Fmor{X}{\mdl} \lra \varprojlim_C \Emo{C}{\mdl_{\Cgst}}$\; 
is an isomorphism. 
\eCor

\subsection{Descent of the Base Field}
\label{sub: Descent}

Consider now a projective variety $X$ with modulus $\mdl$ 
over a finite field $\finfld$. 
Assume $X$ to be smooth outside of $\mdl_X$. 
Let $\fld$ be an algebraic closure of $\finfld$. 

\bPnt 
\label{descent_group}
Let $C$ be a curve in $X$ over $\finfld$. 
The Chow group of the base change of $C$ to $\fld$ 
is given by the $\fld$-valued points of the Jacobian of $C$: 
\[ \xymatrix{
	\hspace{-28mm}
	\Jacm{\Ct}{\mdl}(\fld) = \CHm{\Ct \tens \fld}{\mdl \tens \fld} 
	\ar@<-4.25mm>[d]^{\mod \KZm{C \tens \fld}{\mdl \tens \fld}} \\ 
	\hspace{13mm}
	\Jacm{C}{\mdl}(\fld) = \QHm{C \tens \fld}{\mdl \tens \fld} 
	= \CHm{C \tens \fld}{\mdl \tens \fld}. 
	}
\] 
Since $C$, $\mdl$ and hence $\Sing(C,\mdl)$ are already defined over $\finfld$, 
the finite group 
$\KZm{C \tens \fld}{\mdl \tens \fld} = \KZm{C}{\mdl}$ is constant over $\finfld$. 
Therefore $\Jacm{C}{\mdl}$ is the quotient of $\Jacm{\Ct}{\mdl}$ modulo $\KZm{C}{\mdl}$, i.e.\ we obtain 
\[ \xymatrix{
	\hspace{-10mm}
	\Jacm{\Ct}{\mdl}(\finfld) = \CHm{\Ct}{\mdl}  
	\ar@<1.7mm>[d]^{\mod \KZm{C}{\mdl}} \\ 
	\hspace{-10mm}
	\Jacm{C}{\mdl}(\finfld) = \QHm{C}{\mdl}
	}
\] 
and accordingly for the affine parts we get 
\[ \Lm{C}{\mdl}(\finfld) = \AQHm{C}{\mdl}
\laurin  
\]
\ePnt 

\bPnt 
\label{descent_direct-sys}
Let $\anyfld$ be an extension field of $\finfld$. 
Consider the directed system of curves in $X$ over $\anyfld$ 
from Proposition \ref{CHm_indLim}: 
\[ \fC(\anyfld) 
	= \bigst{ \Crv \;\big|\; \Crv \subset X \tens \anyfld \textrm{ curve}, 
						\;\Crv \iscp \vrt{\mdl_X}
}
\laurin 
\]
For any curve $\Crv \subset X \tens \fld$ let $\fld_{\Crv}$ 
be the field of definition for $\Crv$. 
As $\fld_{\Crv}$ is a finite extension of $\finfld$, 
\[ \Curv := \bigcup_{\gam \in \Gal(\fld_{\Crv} | \finfld)} \gam \,\Crv 
\]
is an algebraic curve in $X$ over $\finfld$ which $\Crv$ embeds into.  
Thus we can replace the directed system $\fC(k)$ by $\fC(\finfld)$ 
without changing any limit taken over this directed system. 
This implies that $\underset{C} \varprojlim \,\Emo{C}{\mdl_{\Cgst}}$ 
and $\underset{C} \varinjlim \ACHm{C}{\mdl_{\Cgst}}$ over $\fld$ 
are compatible with Galois descent. 
The pairing from Proposition \ref{Pairing_X} over $\finfld$ then becomes 
\[ \xymatrix{ 
	\Fmor{X}{\mdl} \hspace{2.5mm} \tms \hspace{2,5mm} \ACHm{X}{\mdl} 
	\ar[r] \ar@<-17.5mm>[d]_{\skelm{X}{\mdl}} 
	& \Gm \ar@<-1.8mm>@{=}[d] \phantom{\laurin} \\ 
	\underset{C} \varprojlim \,\Emo{C}{\mdl_{\Cgst}} \tms 
	\underset{C} \varinjlim \AQHm{C}{\mdl_{\Cgst}} 
	\ar[r] \ar@<-10mm>[u]^-{\wr} 
	& \Gm \laurin 
} 
\] 
\ePnt

\newpage 

\section{Skeleton Theorem}
\label{sec: Skeleton}

The topic of this section is the relation between 
relative Cartier divisors on a projective variety $X$ 
and compatible systems of relative Cartier divisors 
on the various curves in $X$, so called ``skeleton divisors'' 
($\see$Definition \ref{skeleton}). 
The goal is to establish a 1-1 correspondence between those. 
This ``Skeleton Theorem'' will be a key tool 
for all what follows in this note, 
but the relevance of this issue 
might go beyond the scope of this paper. 

This section is very technical and cumbersome, 
but the proof of the Skeleton Theorem (Theorem \ref{SkelThm}) 
on page \pageref{Idea} 
gives an overview about what we are doing.

\subsection{Rigidity of Relative Divisors} 
\label{rigidity}

Let $C$ be a smooth proper irreducible curve over a field $k$, 
the exponential characteristic of which we denote by $p$. 
Let $\mdl = \sum_{\pnt\in S} n_{\pnt} \bt{\pnt}$ be an effective divisor on $C$. 
Let $\sD \in \Fmo{C}{\mdl}\pn{R} \subset \Divf_{C}^{0}\pn{R}$ 
be a relative Cartier divisor on $C$ in the Cartier dual $\Fmo{C}{\mdl}$ 
of $\Llm{0}{C}{\mdl}$. 
Lemma \ref{Basis-determines_0} resp.\ \ref{Basis-determines_p} 
will show that $\sD$ is uniquely determined by the values 
$\bigpair{\sD, \alp_{n,\nu,\sig,\mu}}_{C,\mdl}$ 
of the Cartier duality pairing 
on certain $\alp_{n,\nu,\sig,\mu} \in \ACHmp{C}{\mdl}$, 
whenever the set $\st{\alp_{n,\nu,\sig,\mu}}_{n,\nu,\sig,\mu}$ 
forms a \emph{basis} resp.\ \emph{pro-basis of $\ACHmp{C}{\mdl}$}; 
this notion is introduced in 
Definition \ref{basis-ACHm} resp.\ \ref{pro-basis-ACHm}. 

\bExm 
\label{Lambda-ring}
Let $\sO = \resfld[[t]]$ be the ring of formal power series over a ring $\resfld$, 
and $K = \resfld[[t]][t^{-1}]$ the ring of formal Laurent series. 
Define 
\[ \Lam_n := \ker\Biggpn{\WR{\Gm}{\frac{\resfld[[t]]}{(t^n)}} \lra \Gm} 
\] 
i.e.\ $\Lam_n$ is the smooth connected unipotent algebraic $\resfld$-group 
that is characterized by 
$\Lam_n\pn{\resfld} = U_K^{(1)} \big/ U_K^{(n)}$, 
where $U_K^{(m)} = 1 + t^m \resfld[[t]]$ is the $m^{\textrm{th}}$ higher unit group. 
This algebraic group carries the structure of an algebraic ring scheme \linebreak 
$\pn{\Lam_n, +_{\Lam_n}, \cdot_{\Lam_n}}$, 
which is induced by the ring structure of 
the scheme of big Witt vectors $\Lam := \Lam_{\infty}$: 
addition on $\Lam$ is given by multiplication of formal power series, 
while multiplication on $\Lam$ is given by 
\begin{eqnarray*}
\lull \cdot_{\Lam} \lul : \Lam \tms \Lam & \lra & \Lam \\ 
\prod_{i \geq 1} \pn{1 - a_i t^i} \,,\, \prod_{j \geq 1} \pn{1 - b_j t^j} & \lmt & 
\prod_{i,j \geq 1} \Bigpn{1 - a_i^{\frac{j}{(i,j)}} b_j^{\frac{i}{(i,j)}} t^{\frac{ij}{(i,j)}} }^{(i,j)} 
\end{eqnarray*} 
where $(i,j) := \gcd(i,j)$ denotes the greatest common divisor of $i$ and $j$, 
cf.\ \cite[I, \S1, Prop.\ (1.1)]{Bl_algK}. 
\eExm 
 
\bDef 
\label{basis-Lam}
A family $\st{\lam_{\nu}}_{\substack{0 < \nu < n \\ (\nu,p) = 1}} \subset \Lam_n$ 
is called a \emph{basis of $\Lam_n$} \\
if there are $g_{\nu} \in \resfld[[t]] / (t^n)$ s.t.\ 

\noindent 
$\pn{\chr(\resfld) = 0}$ the exponential map \;$\exp: \Gac \lra \Gm$\; 
yields an isomorphism 
\[ \dsum_{0 < \nu < n} \Ga  \iso  \Lam_n \;, 
\hspace{10mm} 
\sum a_{\nu}  \lmt  \prod \exp\pn{g_{\nu} a_{\nu}} 
\] 
\phantom{\pn{$\chr(\resfld) = p}$} 
and \;$\lam_{\nu} = \exp\pn{g_{\nu}}$\; for all $\nu$. 

\noindent 
$\pn{\chr(\resfld) = p}$ the Artin-Hasse exponential \;$\Expah: \Wittc \lra \Gm$\; 
yields an \\ 
\phantom{\pn{$\chr(\resfld) = p}$} isomorphism 
\[ \dsum_{\substack{0 < \nu < n \\ (\nu,p) = 1}} \Witt_{r(\nu,n)} 
 \iso  \Lam_n \;, 
\hspace{10mm} 
\sum w_{\nu}  \lmt  \prod \Expah\pn{g_{\nu} \cdot w_{\nu}} 
\] 
\phantom{\pn{$\chr(\resfld) = p}$} 
where \;$r(\nu,n) := \min\st{i \;|\; \nu p^i \geq n}$\; \\ 
\phantom{\pn{$\chr(\resfld) = p}$} 
and \;$\lam_{\nu} = \Expah\pn{g_{\nu}}$\; for all $\nu$. 
\eDef

\bRmk 
\label{basis-Criterion}
A family $\st{g_{\nu}}_{\substack{0 < \nu < n \\ (\nu,p) = 1}}$ 
gives a basis of $\Lam_n$ 
if \;$\ord\pn{g_{\nu}} = \nu$\; for all $\nu$, 
if $\resfld$ is a field, 
$\seecite$\cite[V, No.~15, Prop.~8 and No.~16, Prop.~9]{S59}. 
\eRmk 

\bRmk 
\label{dual-decomp}
Via duality, any choice of basis of $\Lam_n$ 
determines a splitting of the Cartier dual of $\Lam_n$: 
\begin{tabular}[t]{ll} 
$\Lam_n^{\vee} \cong \dsum_{0 < \nu < n} \Gac$ & $\textrm{if } \chr(\resfld) = 0$, \\ 
$\Lam_n^{\vee} \cong \dsum_{\substack{0 < \nu < n \\ (\nu,p) = 1}} \; \Wcfl{r(\nu,n)}$ 
 & $\textrm{if } \chr(\resfld) = p$. 
\end{tabular} 
\eRmk

\bDef 
\label{ACHm'}
We use the notation fixed at the beginning of Subsection \ref{rigidity}. 
Let \,$\strl: C \lra \Spec k$\, be the structural morphism. 
Set 
\[ \ACHmp{C}{\mdl} 
:= \strl_*\sO_{\mdl}^* 
 = \frac{\bigcap_{\pnt \in \vrt{\mdl}} \sO_{C,\pnt}^*} 
       {\bigcap_{\pnt \in \vrt{\mdl}} \pn{1 + \fm_{C,\pnt}^{n_{\pnt}}}} 
\] 
and let 
\[ \Lmp{C}{\mdl} 
:= \WR{\Gm}{\strl_*\sO_{\mdl}} 
\] 
be the affine algebraic $k$-group with 
\[ \Lmp{C}{\mdl}\pn{k} = \ACHmp{C}{\mdl} \laurin 
\] 
\eDef 

\bPrp 
\label{duality'}
The group $\Lmp{C}{\mdl}$ is Cartier dual to $\Fm{C}{\mdl}$, 
and this duality 
is expressed by the perfect pairing 
\begin{align*} 
\Fm{\Crv}{\mdl_{\Crvgst}} \;\tms\; \ACHmp{\Crv}{\mdl_{\Crvgst}} \; & \lra \; \Gm \\ 
\bigpair{\sD,\bt{f}}_{\Crv,\mdl} & = \,\prod_{\pnt \in \vrt{\mdl}} \,\pn{\sD,f}_{\pnt} 
\end{align*} 
inspired by the pairing from Definition \ref{Pairing_curve,S}. 
\ePrp 

\bPf 
Let $\extfld := \H^0\pn{C,\sO_C}$. 
According to the dual exact sequences 
\[ 0 \lra \Fmo{C}{\mdl} \lra \Fm{C}{\mdl} \xra{\deg} \WeilRes{\extfld}{\fld} \Zint \lra 0 
\] 
and 
\[ 1 \lla \ACHm{C}{\mdl} \lla \ACHmp{C}{\mdl} \lla \extfld^* \lla 1 
\] 
the perfectness of the pairings induced by  $\pair{\llul,\lull}_{C,\vrt{\mdl}}$ 
on $\Fmo{C}{\mdl} \tms \ACHm{C}{\mdl}$ and on 
$\WeilRes{\extfld}{\fld} \Zint \tms \extfld^*$ 
extends to $\Fm{C}{\mdl} \tms \ACHmp{C}{\mdl}$. 
\ePf 

\bPnt 
\label{decomposition ULT}
Consider the following decomposition of $\ACHmp{C}{\mdl}$: 
\[ \xymatrix{ 
& & \frac{\bigcap_{\pnt} \sO_{C,\pnt}^*} 
       {\bigcap_{\pnt} \pn{1 + \fm_{C,\pnt}^{n_{\pnt}}}} 
       \ar[d]^{\wr} 
& & \\ 
1 \ar[r] & 
\prod_{\pnt} \frac{1 + \fm_{C,\pnt}}{1 + \fm_{C,\pnt}^{n_{\pnt}}}   \ar[r] & 
\prod_{\pnt}  \frac{\sO_{C,\pnt}^*}{1 + \fm_{C,\pnt}^{n_{\pnt}}}   \ar[r] & 
\prod_{\pnt} \fld(\pnt)^*    \ar[r]  & 1} 
\]        
If the residue fields $\fld(\pnt)$ are separable over $\fld$ for all $\pnt \in \vrt{\mdl}$, 
the lower row of this diagram 
corresponds to a decomposition of $\Lmp{C}{\mdl}$ 
into a unipotent part $\Ump{C}{\mdl}$ 
and a torus part $\Tmp{C}{\mdl}$: 
\bMyEqn \label{ULT} \mytag{ULT} 
0 \lra \Ump{C}{\mdl} \lra \Lmp{C}{\mdl} \lra \Tmp{C}{\mdl} \lra 0 \laurin 
\eMyEqn
\ePnt 

\bDef 
\label{n-part}
Let $\Sm := \vrt{\mdl}$ be the support of 
$\mdl = \sum_{\pnt \in \Sm} n_{\pnt} \bt{\pnt}$, 
\[ \Sm(n) := \st{\pnt \in \Sm \;|\; n_{\pnt} = n} 
\] 
and 
\[ N(\mdl) = \st{n \in \Nat \;|\; \Sm(n) \neq \varnothing} 
\] 
hence $\mdl = \sum_{n \in N(\mdl)} \sum_{\pntt \in \Sm(n)} n \,\bt{\pntt}$. 
Set 
\[ \resfld(n) := \prod_{\pntt \in \Sm(n)} \fld(\pntt) 
\] 
the product of the residue fields of points in $\Sm(n)$, 
and 
\[ \sm(n) := \dim_k \resfld(n) 
\laurin 
\] 
Define the algebraic $\fld(\pntt)$-group $\Lam_{\pntt,n_{\pntt}}$ by 
\[ \Lam_{\pntt,n_{\pntt}}\bigpn{\fld(\pntt)} 
:= \frac{1 + \fm_{C,\pntt}}{1 + \fm_{C,\pntt}^{n_{\pntt}}} 
\] 
and the algebraic $\resfld(n)$-group $\Lam_{n,\resfld(n)}$ by 
\[ \Lam_{n,\resfld(n)} 
= \prod_{\pntt \in \Sm(n)} \Lam_{\pntt,n_{\pntt}} 
\] 
for $n \in N(\mdl)$. 
Let $\pjn_{n}$ 
be the composition of the isomorphism obtained 
from the Approximation Lemma ($\see$\cite[I, \S3]{S66}) 
and the projection to the $n$-component: 
\[ \pjn_{n}: \; \frac{\bigcap_{\pnt \in \vrt{\mdl}} \sO_{C,\pnt}^*}
                                 {\bigcap_{\pnt \in \vrt{\mdl}} \pn{1 + \fm_{C,\pnt}^{n_{\pnt}}}} 
   \; \underset{\sim} {\xra{\mathrm{apr}}} \; 
   \prod_{\pnt \in \vrt{\mdl}} \frac{\sO_{C,\pnt}^*}{1 + \fm_{C,\pnt}^{n_{\pnt}}} 
   \; \overset{\pr_{n}} \llsur \; 
   \prod_{\pntt \in \Sm(n)} \frac{\sO_{C,\pntt}^*}{1 + \fm_{C,\pntt}^{n_{\pntt}}} 
   \laurin 
\] 
\eDef

\bDef 
\label{basis-ACHm}
Assume the base field $k$ is algebraically closed of characteristic $0$. 
We call a family 
\[ \bigst{\alp_{n,\nu,\sig}}_{n,\nu,\sig} 
= \bigcup_{\substack{n \in N(\mdl) \\ 1 \leq \sig \leq \sm(n)}} 
   \bigst{\alp_{n,0,\sig}} 
   \cup \bigst{\alp_{n,\nu,\sig}}_
            {1 \leq \nu < n} 
\hspace{3mm} \subset \; 
\frac{\bigcap_{\pnt \in \vrt{\mdl}} \sO_{C,\pnt}^*}
        {\bigcap_{\pnt \in \vrt{\mdl}} \pn{1 + \fm_{C,\pnt}^{n_{\pnt}}}} 
\] 
a \emph{basis of $\ACHmp{C}{\mdl}$} 
if 
\begin{align*}
\textrm{(B1)} \hspace{5mm} & \pjn_{n'}\pn{\alp_{n,\nu,\sig}} = 1 
  \hspace{19mm} \forall \,n' \neq n, \; \forall \, 0 \leq \nu < n, \; 
  \forall \,1 \leq \sig \leq \sm(n) 
  \laurin \\ 
\textrm{(B2)} \hspace{5mm} & 
  \textrm{The system } \hspace{5mm} 
  \prod_{\pntt \in \Sm(n)} \alp_{n,0,\sig}(\pntt)^{\nbr_{\pntt}} = 1, 
  \hspace{10mm} \sig = 1, \ldots, \sm(n) \\ 
 & \textrm{has only the trivial solution } \; 
 \pn{\nbr_{\pntt}}_{\pntt} = 0 \in \Zint^{\Sm(n)} 
 \hspace{10mm} \forall \, n \in N(\mdl) 
 \laurin \\ 
\textrm{(B3)} \hspace{5mm} &  
  \textrm{There is a basis }\, 
  \bigst{\exp\pn{g_{n,\nu}}}_{1 \leq \nu < n} 
  \,\textrm{ of }\, \Lam_{n,\resfld(n)} \hspace{10mm} \forall \,n \in N(\mdl) \\ 
  & \textrm{and a $k$-basis } \; \st{c_{n,\nu,\sig}}_{1 \leq \sig \leq \sm(n)} 
  \,\textrm{ of }\, \resfld(n) 
  \hspace{15mm} 
  \forall\,n, \;\forall\,0 < \nu < n \\ 
  & \textrm{such that } \;
  \pjn_{n}\pn{\alp_{n,\nu,\sig}} = 
  \exp\pn{g_{n,\nu} \,c_{n,\nu,\sig}} 
  \hspace{16mm} \forall \,n, \,\nu \geq 1, \,\sig 
  \laurin 
\end{align*}
\eDef 

\bDef 
\label{pro-basis-ACHm}
Assume the base field $k$ is algebraically closed of characteristic $p > 0$. 
We call a family 
\[ \bigst{\alp_{n,\nu,\sig,\mu}}_{n,\nu,\sig,\mu} 
= \bigcup_{\substack{n \in N(\mdl) \\ 1 \leq \sig \leq \sm(n) \\ \mu \geq 1}}
   \bigst{\alp_{n,0,\sig,\mu}} 
   \cup \bigst{\alp_{n,\nu,\sig,\mu}}_
   {\substack{1 \leq \nu < n \\ (\nu,p) = 1}} 
\hspace{0mm} \subset \; \frac{\bigcap_{\pnt \in \vrt{\mdl}} \sO_{C,\pnt}^*}
        {\bigcap_{\pnt \in \vrt{\mdl}} \pn{1 + \fm_{C,\pnt}^{n_{\pnt}}}} 
\] 
a \emph{pro-basis of $\ACHmp{C}{\mdl}$} 
if 
\begin{align*}
\textrm{(PB1)} \hspace{3mm} & \pjn_{n'}\pn{\alp_{n,\nu,\sig,\mu}} = 1 
  \hspace{9mm} 
  \forall \,n' \neq n, \;\forall \,0 \leq \nu < n,\;\forall \,1\leq \sig \leq \sm(n),\; 
  \forall\,\mu \geq 1 
  \laurin \\ 
\textrm{(PB2)} \hspace{3mm} & 
  \textrm{The system } \hspace{0mm} 
  \prod_{\pntt \in \Sm(n)} \alp_{n,0,\sig,\mu}(\pntt)^{\nbr_{\pntt}} = 1, 
  \hspace{4mm} \sig = 1, \ldots, \sm(n) \hspace{5mm} \textrm{has only the } \\ 
  & \textrm{trivial solution } 
  \,\pn{\nbr_{\pntt}}_{\pntt} = 0 \, 
  \textrm{ in } \st{-p^{\mu},\ldots,p^{\mu}}^{\Sm(n)} 
  \hspace{3mm} 
  \forall \, n \in N(\mdl),\;\forall \mu \geq 1 
  \laurin \\ 
\textrm{(PB3)} \hspace{3mm} &  
  \textrm{There is a basis }\, 
  \bigst{\Expah\pn{g_{n,\nu}}}_{\substack{0 < \nu < n \\ (\nu,p) = 1}} 
  \,\textrm{ of }\, \Lam_{n,\resfld(n)} \hspace{10mm} \forall \,n \in N(\mdl), \\ 
  & \textrm{a $k$-basis } \; \st{c_{n,\nu,\sig}}_{1 \leq \sig \leq \sm(n)} 
  \,\textrm{ of } \, \resfld(n) \hspace{13mm} 
  \forall \,n \in N(\mdl), \;\forall\,0 < \nu < n \\ 
  & \textrm{and a family } \; \st{c_{\mu}}_{\mu \geq 1} \subset k 
  \;\textrm{ with }\; c_{\mu} \neq c_{\mu'} 
  \hspace{11mm} \forall \,\mu \neq \mu' \geq 1 \\ 
  & \textrm{such that } \;
  \pjn_{n}\pn{\alp_{n,\nu,\sig,\mu}} = 
  \Expah\pn{g_{n,\nu} \,c_{n,\nu,\sig} \,c_{\mu}} 
  \hspace{10mm} \forall \,n, \,\nu \geq 1, \,\sig, \,\mu 
  \laurin 
\end{align*}
\eDef 

\bDef 
\label{pro-basis-abuse}
By abuse of notation we will call a family 
\[ \bigst{\tha_{n,\nu,\sig,\mu}}_{n,\nu,\sig,\mu} 
\subset \bigcap_{\pnt \in \vrt{\mdl}} \sO_{C,\pnt}^* 
\subset \sK_C^* 
\] 
a \emph{(pro-)basis of $\ACHmp{C}{\mdl}$}, 
if the family of the corresponding classes in 
$\frac{\bigcap_{\pnt \in \vrt{\mdl}} \sO_{C,\pnt}^*}
        {\bigcap_{\pnt \in \vrt{\mdl}} \pn{1 + \fm_{C,\pnt}^{n_{\pnt}}}} 
 = \ACHmp{C}{\mdl}$ is a (pro-)basis of $\ACHmp{C}{\mdl}$. 
\eDef 

\bDef 
\label{pro-basis-C-in-X}
Let $\pn{X, \mdl}$ be is as in the beginning of Section \ref{AbelJacobiMap}, 
and let $C$ be a curve in $X$, not necessarily normal or irreducible. 
We call a family 
\[ \bigst{\tha_{n,\nu,\sig,\mu}}_{n,\nu,\sig,\mu} 
\subset \sK_C^* 
\] 
a \emph{(pro-)basis of} $\ACHmp{C}{\mdl}$ 
\footnote{ Here we define the symbolic expression 
``(pro-)basis of $\ACHmp{C}{\mdl}$'', 
while we do not care whether the expression ``$\ACHmp{C}{\mdl}$'' 
is defined or not.} 
if 
\[ \bigst{\wt{\tha}_{n,\nu,\sig,\mu}}_{n,\nu,\sig,\mu} 
\subset \sK_{\Ct}^* 
\,\textrm{ is a (pro-)basis of } \ACHmp{\Ct}{\mdl_{\Ct}} 
\laurin 
\] 
\eDef 

\bPrp 
\label{basis-exists}
Assume $k$ is algebraically closed. 
If $\chr(k) = 0$, a basis of $\ACHmp{C}{\mdl}$ exists, 
and  if $\chr(k) = p > 0$, a pro-basis of $\ACHmp{C}{\mdl}$ exists. 
\ePrp 

\bPf 
For every $n \in N\pn{\mdl}$ we construct the set 
$\bigst{\tha_{n,\nu,\sig}}_{0 \leq \nu < n, \,1 \leq \sig \leq \sm(n)}$ 
respectively 
$\bigst{\tha_{n,\nu,\sig,\mu}}_
{\substack{0 \leq \nu < n \\ (\nu_{\pnt},p) = 1}, \,1 \leq \sig \leq \sm(n), \,\mu \geq 1}$ 
as follows. 

Note that $\#\Sm(n) = \sm(n)$ since the base field $k$ is algebraically closed. 
Hence we can write $\Sm(n) = \st{\,\pnt_{\sig} \;|\; 1 \leq \sig \leq \sm(n)}$ 
and $\resfld(n) = \prod_{1 \leq \sig \leq \sm(n)} \fld(\pnt_{\sig})$. 
For every $1 \leq \sig \leq \sm(n)$, 
if $\chr(k) = 0$ 
choose $a_{\sig} \in \fld(\pnt_{\sig})^*$ not a root of unity (e.g.\ $a_{\pnt} = 2$); 
in the case of $\chr(k) = p$ 
let $a_{\sig,\mu} \in \fld(\pnt_{\sig})^*$ be a generator of the cyclic group 
$\bF_{p^{\mu+1}}^*$ for every $\mu \geq 1$. 
According to the Approximation Lemma ($\see$\cite[I, \S3]{S66}) there exist 
$\tha_{n,0,\sig}$ resp.\ 
$\tha_{n,0,\sig,\mu} \in \bigcap_{\pnt \in \vrt{\mdl}} \sO_{C,\pnt}^*$ 
such that $\tha_{n,0,\sig}(\pnt_{\sig}) = a_{\sig}$ 
resp.\ $\tha_{n,0,\sig,\mu}(\pnt_{\sig}) = a_{\sig,\mu}$ 
and $\tha_{n,0,\sig}(\pnt_{\sig'}) = 1$ 
resp.\ $\tha_{n,0,\sig,\mu}(\pnt_{\sig'}) = 1$ 
for all $\sig' \neq \sig$, 
and $\pjn_{n'}\pn{\tha_{n,0,\sig}} = 1$ 
resp.\ $\pjn_{n'}\pn{\tha_{n,0,\sig,\mu}} = 1$ 
for all $n' \neq n$. 
Obviously the set 
$\st{\tha_{n,0,\sig} \;|\; 1 \leq \sig \leq \sm(n)}$ satisfies condition (B2), 
resp.\ $\st{\tha_{n,0,\sig,\mu} \;|\; 1 \leq \sig \leq \sm(n), \,\mu \geq 1}$ 
satisfies condition (PB2). 

Let $\st{\lam_{\pnt,\nu}}_{0 < \nu < n}$ 
$= \st{\exp\pn{g_{\pnt,\nu}}}_{0 < \nu < n}$ 
be a basis of 
$\Lam_{\pnt,n} = \frac{1 + \fm_{C,\pnt}}{1 + \fm_{C,\pnt}^{n}}$ 
if $\chr(k) = 0$, 
respectively $\st{\lam_{\pnt,\nu}}_{\substack{0 < \nu < n \\ (\nu,p) = 1}}$ 
$= \st{\Expah\pn{g_{\pnt,\nu}}}_{\substack{0 < \nu < n \\ (\nu,p) = 1}}$ 
if $\chr(k) = p$, 
for every $\pnt \in \Sm(n)$; 
such a basis always exists by Remark \ref{basis-Criterion}. 
Let $g_{n,\nu} = \pn{g_{\pnt,\nu}}_{\pnt \in \Sm(n)}$ 
and $\lam_{n,\nu} = \pn{\lam_{\pnt,\nu}}_{\pnt \in \Sm(n)}$, 
so  $\st{\lam_{n,\nu}}_{\substack{0 < \nu < n \\ (\nu,p) = 1}}$ 
yields a basis of $\Lam_{n,\resfld(n)} \cong \prod_{\pnt \in \Sm(n)} \Lam_{\pnt,n}$. 
Set $c_{n,\nu,\sig} := \pn{\del_{\sig \sig'}}_{1 \leq \sig' \leq \sm(n)} \in \resfld(n)$ 
for all $\nu$, $\sig$, 
where $\del_{\sig \sig'}$ is the Kronecker delta, 
and let $c_{\mu}$ be a generator of the cyclic group $\bF_{p^{\mu+1}}^*$. 
By the Approximation Lemma ($\see$\cite[I, \S3]{S66}) there are functions 
$\tha_{n,\nu,\sig}$ resp.\ 
$\tha_{n,\nu,\sig,\mu} \in \bigcap_{\pnt \in \vrt{\mdl}} \sO_{C,\pnt}^*$ 
such that 
$\pjn_{n}\pn{\tha_{n,\nu,\sig}} = \exp\pn{g_{n,\nu} c_{n,\nu,\sig}}$ 
for all $\nu$, $\sig$, 
and with $\pjn_{n'}\pn{\tha_{n,\nu}} = 1$ 
for all $n' \neq n$ and all $\nu$, $\sig$, 
resp.\ 
$\pjn_{n}\pn{\tha_{n,\nu,\sig,\mu}} = \Expah\pn{g_{n,\nu} c_{n,\nu,\sig} c_{\mu}}$ 
for all $\nu$, $\sig$, $\mu$, 
and with $\pjn_{n'}\pn{\tha_{n,\nu}} = 1$ 
for all $n' \neq n$ and all $\nu$, $\sig$, $\mu$. 
Then conditions (B1), (B2), (B3) of Definition \ref{basis-ACHm} 
resp. (PB1), (PB2), (PB3) of Definition \ref{pro-basis-ACHm} 
are satisfied by construction. 
\ePf 

\bLem 
\label{Basis-determines_0}
Assume $k$ is algebraically closed of characteristic $0$. 
Let \linebreak 
$\st{\tha_{n,\nu,\sig}}_{n,\nu,\sig}$ 
be a basis of $\ACHmp{C}{\mdl}$. 
Let $R$ be a local finite $\fld$-ring. 
Let 
$\st{t_{n,\sig}}_{n \in N\pn{\mdl}, 1 \leq \sig \leq \sm(n)}$ 
be a family in 
$\Trz_R\pn{\fld} = \Gm\pn{R_{\red}}$ 
such that the system of equations 
$\prod_{\pntt \in \Sm(n)} \tha_{n,0,\sig}(\pntt)^{\nbr_{\pntt}} = t_{n,\sig}$, 
$\sig = 1, \ldots, \sm(n)$ 
is solvable in $\Zint^{\Sm(n)}$ 
for all $n \in N\pn{\mdl}$, 
and let 
$\st{u_{n,\nu,\sig}}_{n \in N\pn{\mdl}, \,0 < \nu < n, \,1 \leq \sig \leq \sm(n)}$ 
be a family in 
$\Upf_R\pn{\fld} = \ker\bigpn{\Gm\pn{R} \ra \Gm\pn{R_{\red}}}$. 

Then there is a unique relative Cartier divisor 
$\sD = \sD_{\et} + \sD_{\inf} \in \Fm{C}{\mdl}\pn{R}$ 
such that 
\;$\bigpair{\sD_{\et},\tha_{n,0,\sig}}_{C,\mdl} = t_{n,\sig}$\; 
for all $n \in N\pn{\mdl}$, $1 \leq \sig \leq \sm(n)$, 
and 
\;$\bigpair{\sD_{\inf},\tha_{n,\nu,\sig}}_{C,\mdl} = u_{n,\nu,\sig}$\; 
for all $n \in N\pn{\mdl}$, $0 < \nu < n$, $1 \leq \sig \leq \sm(n)$. 
\eLem 

\bPf 
By assumption the base field $k$ is in particular perfect, 
so the formal group $\Fm{C}{\mdl}$ splits canonically into 
\'etale part and infinitesimal part: 
\[ \Fm{C}{\mdl} = \pn{\Fm{C}{\mdl}}_{\et} \tms \pn{\Fm{C}{\mdl}}_{\inf} 
\] 
and we may treat these parts separately. 

\medskip 
\textbf{Step 1:} The \'etale part of $\Fm{C}{\mdl}$ is given by 
\[ \pn{\Fm{C}{\mdl}}_{\et}\pn{R} 
\,\cong\, \dsum_{\pnt \in \vrt{\mdl}} \Zint 
\,=\, \dsum_{n \in N(\mdl)} \dsum_{\pntt \in \Sm(n)} \Zint 
\,\laurin 
\] 
As by assumption the system of equations 
\[ \bigpair{\sD_{\et},\tha_{n,0,\sig}}_{C,\mdl} 
= \prod_{\pntt \in \Sm(n)} \tha_{n,0,\sig}(\pntt)^{\nbr_{\pntt}} 
= t_{n,\sig}
\hspace{15mm} \sig = 1, \ldots, \sm(n) 
\]  
has a solution 
$\pn{\nbr_{\pntt}}_{\pntt \in \Sm(n)} \in \Zint^{\Sm(n)}$ 
for every $n \in N(\mdl)$, 
there is an \'etale relative Cartier divisor 
$\sD_{\et} \cong \pn{\nbr_{\pnt}}_{\pnt \in \vrt{\mdl}}$ 
solving these equations, 
and due to Definition \ref{basis-ACHm} (B2) 
a divisor with this property is uniquely determined. 

\medskip 
\textbf{Step 2:} By \Point \ref{decomposition ULT} 
and due to the bases $\st{\exp\pn{g_{n,\nu}}}_{0 < \nu < n}$ 
of $\Lam_{n,\resfld(n)}$ 
from Definition \ref{basis-ACHm} (B3) 
for every $n \in N(\mdl)$, 
the unipotent part of $\Lmp{C}{\mdl}$ has a decomposition 
\[ \Ump{C}{\mdl} 
= \prod_{n \in N(\mdl)} \Lam_{n,\resfld(n)} 
\cong \dsum_{n \in N(\mdl)} \dsum_{0 < \nu < n} \GaS{\resfld(n)} 
\laurink 
\] 
where $\Lam_{\pnt,n_{\pnt}}(k) = \frac{1 + \fm_{C,\pnt}}{1 + \fm_{C,\pnt}^{n_{\pnt}}}$ 
and $\Lam_{n,\resfld(n)} = \prod_{\pntt \in \Sm(n)} \Lam_{\pntt,n_{\pntt}}$ 
($\see$Definition \ref{n-part}). 
By duality we obtain a decomposition 
of the infinitesimal part of $\Fm{C}{\mdl}$ 
\[ \pn{\Fm{C}{\mdl}}_{\inf} 
= \prod_{n \in N\pn{\mdl}} \,\pn{\Fm{C}{\mdl}}_{\inf}^{n} 
\cong \dsum_{n \in N(\mdl)} \dsum_{0 < \nu < n} \GacS{\resfld(n)} 
\laurin 
\] 
According to these decompositions, 
an infinitesimal divisor $\sD_{\inf} \in \pn{\Fm{C}{\mdl}}_{\inf}(R)$ 
corresponds to a tuple 
$\pn{a^{n,\nu}}_{n,\nu} \in 
\dsum_{n \in N(\mdl)} \dsum_{0 < \nu < n} \GacS{\resfld(n)}(R_{\resfld(n)})$ 
with $a^{n,\nu} = \pn{a^{n,\nu}_{\pntt}}_{\pntt} 
\in \prod_{\pntt \in \Sm(n)} \fld(\pntt) \tens_{\fld} R 
= \resfld(n) \tens_{\fld} R$; 
by (B3) for each $n' \in N\pn{\mdl}$, $0 < \nu' < n'$ 
the basis element $\tha_{n',\nu',\sig'} \in \Ump{C}{\mdl}(k)$ 
corresponds to the tuple 
$\pn{\del_{n' n} \,\del_{\nu' \nu} \,c_{n',\nu',\sig'}}_{n,\nu} \in 
\dsum_{n \in N(\mdl)} \dsum_{0 < \nu < n,} \GaS{\resfld(n)}(k_{\resfld(n)})$ 
with $c_{n,\nu,\sig} = \pn{c_{n,\nu,\sig}^{\,\pntt}}_{\pntt}$ 
$\in \prod_{\pntt \in \Sm(n)} \fld(\pntt) = \resfld(n)$ 
and where $\del_{\alp \bet}$ denotes the Kronecker delta. 
Then 
\[ \bigpair{\sD_{\inf},\tha_{n,\nu,\sig}}_{C,\mdl} 
= \prod_{\pntt \in \Sm(n)} \exp\bigpn{a^{n,\nu}_{\pntt} \,c_{n,\nu,\sig}^{\,\pntt}} 
= \exp\biggpn{\sum_{\pntt \in \Sm(n)} a^{n,\nu}_{\pntt} \,c_{n,\nu,\sig}^{\,\pntt}} 
\] 
As \;$\exp: \Nil(R) \lra 1 + \Nil(R) = \Upf_R(k)$\, is invertible, 
we can set 
\[ v_{n,\nu,\sig} := \exp^{-1}\pn{u_{n,\nu,\sig}} 
\laurin 
\] 
Then for every $n, \nu$ the conditions 
\;$\bigpair{\sD_{\inf},\tha_{n,\nu,\sig}}_{C,\mdl} \,=\, u_{n,\nu,\sig} 
\hspace{10mm} \forall \,\sig$ \\ 
are expressed by the matrix equation 
\[ \pn{a^{n,\nu}_{\pntt}}_{\pntt \in \Sm(n)} \, 
\bigpn{c_{n,\nu,\sig}^{\,\pntt}}_{\substack{\pntt \in \Sm(n) \\ 1 \leq \sig \leq \sm(n)}} 
= \pn{v_{n,\nu,\sig}}_{1 \leq \sig \leq \sm(n)} 
\laurin 
\] 
Since the matrix 
$\pn{c_{n,\nu,\sig}^{\,\pntt}}_{\pntt,\sig}$ 
is invertible by (B3), 
these conditions determine uniquely an infinitesimal relative Cartier divisor 
$\sD_{\inf} \cong \pn{a_{n,\nu}}_{n,\nu}$. 

This yields a unique relative Cartier divisor 
$\sD = (\sD_{\et},\sD_{\inf}) \in \Fm{C}{\mdl}(R)$ 
satisfying the required properties. 
\ePf 


\bLem 
\label{Basis-determines_p}
Assume that $k$ is algebraically closed of characteristic $p > 0$. 
Let 
$\st{\tha_{n,\nu,\sig,\mu}}_{n,\nu,\sig,\mu}$ 
be a pro-basis of $\ACHmp{C}{\mdl}$. 
Let $R$ be a finite $\fld$-ring. 
\linebreak 
Consider families 
\,$\st{t_{n,\sig,\mu}}_{n \in N\pn{\mdl}, 1 \leq \sig \leq \sm(n), 1 \leq \mu}$\, 
in \,$\Trz_R\pn{\fld} := \Gm\pn{R_{\red}}$ 
and \linebreak 
\,$\st{u_{n,\nu,\sig,\mu}}_
{n \in N\pn{\mdl}, 0 < \nu < n, 1 \leq \sig \leq \sm(n), 1 \leq \mu}$\, 
in \,$\Upf_R\pn{\fld} := \ker\bigpn{\Gm\pn{R} \ra \Gm\pn{R_{\red}}}$. \\ 
For an indeterminate relative Cartier divisor $\sD \in \Fm{C}{\mdl}\pn{R}$ 
with \'etale part $\sD_{\et}$ and infinitesimal part $\sD_{\inf}$ 
consider the system of equations 
\begin{align} 
& & & \fral \,n \in N\pn{\mdl}, \;\;\fral \,1 \leq \sig \leq \sm(n) \notag \\ 
\bigpair{\sD_{\et},\tha_{n,0,\sig,\mu}}_{C,\mdl} &\,=\, t_{n,\sig,\mu} 
 & & \fral \,\mu \gg 1 \label{Sy} \mytag{Sy} 
 \\ 
\bigpair{\sD_{\inf},\tha_{n,\nu,\sig,\mu}}_{C,\mdl} &\,=\, u_{n,\nu,\sig,\mu} 
& & \fral \,\mu \geq 1, 
\;\;\fral \,0 < \nu < n, \,\pn{\nu,p} = 1 
\laurin \notag
\end{align} 
If the system (\ref{Sy}) is solvable, 
then the solution $\sD$ is unique 
and can be found via an explicit procedure, 
described in Algorithm \ref{algorithm}. 
\eLem 

\bPf 
Since $\Fm{C}{\mdl}$ commutes with finite products 
and a finite ring is the product of finitely many local rings 
we may assume that $R$ is local. 
As in the proof of Lemma \ref{Basis-determines_0} 
we may treat \'etale part and infinitesimal part of $\Fm{C}{\mdl}$ separately. 

\medskip 
\textbf{Step 1:} The \'etale part of $\Fm{C}{\mdl}$ is given by 
\bMyEqn \label{FZ} \mytag{FZ} 
\pn{\Fm{C}{\mdl}}_{\et}\pn{R} 
\,\cong\, \dsum_{\pnt \in \vrt{\mdl}} \Zint 
\,=\, \dsum_{n \in N(\mdl)} \dsum_{\pntt \in \Sm(n)} \Zint 
\,\laurin 
\eMyEqn  
According to this decomposition, 
any \'etale relative divisor 
$\sD_{\et} \in \pn{\Fm{C}{\mdl}}_{\et}(R)$ 
corresponds to a tuple 
\bMyEqn \label{Ln} \mytag{Ln} 
\pn{\nbr_{n,\pntt}}_{n,\pntt} \in \dsum_{n \in N(\mdl)} \dsum_{\pntt \in \Sm(n)} \Zint 
\,\laurin 
\eMyEqn 
Assume that the system above has a solution $\sD \in \Fm{C}{\mdl}$, 
then the system 
\begin{align} 
\label{DL} \mytag{DL} 
& \bigpair{\sD_{\et},\tha_{n,0,\sig,\mu}}_{C,\mdl} 
= \prod_{\pntt \in \Sm(n)} \tha_{n,0,\sig,\nmb}\pn{\pntt}^{\nbr_{n,\pntt}} 
= t_{n,\sig,\mu} \\ 
& \hspace{0mm} 
\forall \,1 \leq \sig \leq \sm(n) \notag 
\end{align} 
has a solution 
$\pn{\nbr_{n,\pntt}}_{\pntt \in \Sm(n)} \in \Zint^{\Sm(n)}$ 
for every $n \in N\pn{\mdl}$, 
valid for sufficiently large $\mu$. 
Then there is an $m \in \Nat$ with $\vrt{\nbr_{n,\pntt}} \leq p^{\nmb}$ 
for all $n, \pntt$. 
By Definition \ref{pro-basis-ACHm} (PB2) 
these conditions determine the divisor 
$\sD_{\et} \cong \pn{\nbr_{n,\pntt}}_{n,\pntt}$ uniquely. 

\medskip 
\textbf{Step 2:} By \Point \ref{decomposition ULT} 
and due to the bases 
$\st{\Expah\pn{g_{n,\nu}}}_{\substack{0 < \nu < n \\ (\nu,p) = 1}}$ 
of $\Lam_{n,\resfld(n)}$ 
from Definition \ref{pro-basis-ACHm} (PB3) 
for every $n \in N(\mdl)$, 
the unipotent part of $\Lmp{C}{\mdl}$ has a decomposition 
\bMyEqn \label{UW} \mytag{UW} 
\Ump{C}{\mdl} 
= \prod_{n \in N\pn{\mdl}} \Lam_{n,\resfld(n)} 
\cong \dsum_{n \in N\pn{\mdl}} \dsum_{\substack{0 < \nu < n \\ (\nu,p) = 1}} 
\Witt_{r\pn{\nu,n},\resfld(n)} 
\laurin 
\eMyEqn
By duality we obtain a decomposition 
of the infinitesimal part of $\Fm{C}{\mdl}$ 
\bMyEqn \label{FW} \mytag{FW} 
\pn{\Fm{C}{\mdl}}_{\inf} 
= \prod_{n \in N\pn{\mdl}} \pn{\Fm{C}{\mdl}}_{\inf}^{n} 
\cong \dsum_{n \in N\pn{\mdl}} \dsum_{\substack{0 < \nu < n \\ (\nu,p) = 1}} 
\Wcfl{r\pn{\nu,n}}_{\resfld(n)} 
\laurin 
\eMyEqn 
According to these decompositions, 
any infinitesimal relative divisor 
$\sD_{\inf} \in \pn{\Fm{C}{\mdl}}_{\inf}(R)$ 
corresponds to a tuple 
\bMyEqn \label{Wn} \mytag{Wn} 
\pn{w^{n,\nu}}_{n,\nu} 
\in \dsum_{n \in N\pn{\mdl}} \dsum_{\substack{0 < \nu < n \\ (\nu,p) = 1}} 
\Wcfl{r\pn{\nu,n}}_{\resfld(n)}\bigpn{R_{\resfld(n)}} 
\eMyEqn 
with 
\begin{align} 
\label{Wi} \mytag{Wi} 
w^{n,\nu} = \pn{w^{n,\nu}_{\,\pntt}}_{\pntt \in \Sm(n)} 
& \in \dsum_{\pntt \in \Sm(n)} \Wcfl{r\pn{\nu,n}}(R) \\ 
w^{n,\nu}_{\,\pntt} = \pn{w^{n,\nu}_{i,\pntt}}_{i \geq 0}
& \in \Wcfl{r\pn{\nu,n}}(R) \notag
\laurink 
\end{align} 
and by (PB3) for each $n \in N\pn{\mdl}$, $0 < \nu < n$, 
$1 \leq \sig \leq \sm(n)$, $\mu \geq 1$ 
the basis element $\tha_{n,\nu,\sig,\mu} \in \Ump{C}{\mdl}(k)$ 
corresponds to the tuple 
\[ \mytag{Th} \label{Th}
\pn{\del_{n' n} \,\del_{\nu' \nu} \,c_{n,\nu,\sig} \,c_{\mu}}_{n',\nu'} 
\in \dsum_{n' \in N\pn{\mdl}} \dsum_{\substack{0 < \nu' < n' \\ (\nu',p) = 1}} 
\Witt_{r\pn{\nu',n'},\resfld(n')}\pn{k_{\resfld(n')}} 
\] 
with 
\[ \mytag{CC} \label{CC} 
c_{n,\nu,\sig} = \pn{c_{n,\nu,\sig}^{\,\pntt}}_{\pntt \in \Sm(n)} 
\in \prod_{\pntt \in \Sm(n)} \fld(\pntt) = \resfld(n), 
\hspace{7mm} c_{\mu} \in k 
\] 
and where $\del_{\alp \bet}$ denotes the Kronecker delta. 
Then, using Cartier duality of Witt vectors ($\seecite$\cite[Exm.~1.11]{Ru13}) 
\begin{align} 
\label{DE} \mytag{DE} 
\bigpair{\sD_{\inf},\tha_{n,\nu,\sig,\mu}}_{C,\mdl} 
& =   \prod_{\pntt \in \Sm(n)} 
         \Expah\Bigpn{w^{n,\nu}_{\,\pntt} 
         \cdot \bigbt{c_{n,\nu,\sig}^{\,\pntt}} 
         \cdot \bigbt{c_{\mu}} } \\ 
& =\, \Expah\Bigpn{\sum_{\pntt \in \Sm(n)} w^{n,\nu}_{\,\pntt} 
         \cdot \bigbt{c_{n,\nu,\sig}^{\,\pntt}} 
         \cdot \bigbt{c_{\mu}} } \notag \\ 
& =\, \Expaht\Bigpn{\sum_{\pntt \in \Sm(n)} w^{n,\nu}_{\,\pntt} 
         \cdot \bigbt{c_{n,\nu,\sig}^{\,\pntt}}, 
         \,c_{\mu} } \notag 
\end{align} 
where $\bt{c} \in \Witt(k)$ denotes the Teichm\"uller representative of $c \in k$. \\ 
If we find elements \;$u^{n,\nu,\sig}\pn{t} \in 1 + t R[[t]] = \Lam\pn{R}$\; 
for all \,$n, \nu, \sig$\, such that 
\[ \label{UC} \mytag{UC} 
u^{n,\nu,\sig} \pn{c_{\mu}} \,=\, u_{n,\nu,\sig,\mu} 
\hspace{16mm} \forall \,\mu \geq 1 
\laurink 
\] 
then the conditions 
\;$\bigpair{\sD_{\inf},\tha_{n,\nu,\sig,\mu}}_{C,\mdl} \,=\, u_{n,\nu,\sig,\mu} 
\hspace{13mm} \forall \,\mu \geq 1$ \\ 
are expressed by 
\bMyEqn \label{WU} \mytag{WU} 
   \Expaht\biggpn{\sum_{\pntt \in \Sm(n)} w^{n,\nu}_{\,\pntt} 
   \cdot \bigbt{c_{n,\nu,\sig}^{\,\pntt}}, t } 
   \,=\, u^{n,\nu,\sig} \pn{t}
\laurin 
\eMyEqn 
Assume that there is a solution 
$\pn{w^{n,\nu}}_{n,\nu} 
\in \dsum_{n,\nu} \Wcfl{r\pn{\nu,n}}_{\resfld(n)}\bigpn{R_{\resfld(n)}}$. 
Taking into account that 
\bMyEqn \label{TR} \mytag{TR} 
w^{n,\nu}_{i,\pntt} \in \Nil(R) 
\hspace{6mm} \textrm{and} \hspace{6mm} 
\exists \,\nmb \in \Nat \hspace{3mm} \forall \,i \geq \nmb: 
\hspace{3mm} w^{n,\nu}_{i,\pntt} = 0 
\eMyEqn 
and writing for every \,$n, \nu, \sig$ 
\[ u^{n,\nu,\sig} \pn{t} \,=\, \sum_{i \geq 0} r_i^{n,\nu,\sig} \, t^i 
\hspace{5mm} \in\, 1 + t \, R[[t]] 
\laurink 
\] 
there is hence $\oma \in \Nat$ such that 
\,$r^{n,\nu,\sig}_{i} = 0$\, for $i \geq \oma$. 
Then for every $n, \nu, \sig$ 
the coefficients \,$r_i^{n,\nu,\sig} \in R$\, are uniquely determined 
by the matrix equation 
\[ \bigpn{r_{i}^{n,\nu,\sig}}_{0 \leq i < \oma} \, 
\Bigpn{\pn{c_{\mu}}^{i}}_{\substack{0 \leq i < \oma \\ 1 \leq \mu \leq \oma}} 
\,=\, \bigpn{u_{n,\nu,\sig,\mu}}_{1 \leq \mu \leq \oma} 
\laurin 
\] 
Here the matrix $\bigpn{\pn{c_{\mu}}^{i}}_{i,\oma}$ is invertible, 
since by (PB3) the Vandermonde determinant is 
\[ \det \Bigpn{\pn{c_{\mu}}^{i}}_{\substack{0 \leq i < \oma \\ 1 \leq \mu \leq \oma}} 
\,=\, \prod_{\mu < \mu'} \pn{c_{\mu} - c_{\mu'}} 
\,\neq\, 0 
\laurin 
\] 
Moreover, we have an isomorphism 
\bMyEqn \label{WL} \mytag{WL} 
   \Pi\eps: \prod_{\substack{\nmr \geq 1 \\ \pn{\nmr,p}=1}} \Witt \iso \Lam \;, 
   \hspace{12mm} 
   \pn{v_{\nmr}}_{\nmr} \lmt \prod_{\nmr} \Expaht\bigpn{v_{\nmr}, t^{\nmr}} 
\eMyEqn 
defined over $\Zint_{(p)}$ ($\seecite$\cite[III, No.~1, Prop.\ on p.~53]{D}). 
By (\ref{WU}), the power series 
$u^{n,\nu,\sig}\pn{t}$ lies in the image of 
\,$\Expaht\pn{\llul,t}: \Wittc \lra \Lam$, 
thus only the first component of the inverse image of $u^{n,\nu,\sig}\pn{t}$ 
under the isomorphism \,$\Pi\eps$\, is non-trivial, for all $n, \nu, \sig$. 
Hence there is a unique \,$v^{n,\nu,\sig} \in \Wittc\pn{R}$\, such that 
\bMyEqn \label{VU} \mytag{VU} 
\Expaht\bigpn{v^{n,\nu,\sig}, t} \,=\, u^{n,\nu,\sig}\pn{t} 
\laurin 
\eMyEqn 
As $\Expaht\pn{\llul, t}$ is injective, 
(\ref{WU}) and (\ref{VU}) imply that for every $n, \nu$ 
the Witt vector \,$w^{n,\nu}$\, 
is given by the system of equations 
\[ \sum_{\pntt \in \Sm(n)} w^{n,\nu}_{\,\pntt} \cdot \bigbt{c_{n,\nu,\sig}^{\,\pntt}} 
\,=\, v^{n,\nu,\sig} 
\hspace{16mm} \forall \,1 \leq \sig \leq \sm(n) 
\laurink 
\] 
i.e.\ by the matrix equation 
\[ \label{WV} \mytag{WV} 
\bigpn{w^{n,\nu}_{\,\pntt}}_{\pntt \in \Sm(n)} \, 
\Bigpn{\bigbt{c_{n,\nu,\sig}^{\,\pntt}}}_
{\substack{\pntt \in \Sm(n) \\ 1 \leq \sig \leq \sm(n)}} 
\,=\, \bigpn{v^{n,\nu,\sig}}_{1 \leq \sig \leq \sm(n)} 
\laurin 
\] 
Since the matrix 
$\bigpn{\bt{c_{n,\nu,\sig}^{\,\pntt}}}_{\pntt,\sig}$ 
is invertible by (PB3), 
these conditions 
determine uniquely an infinitesimal relative Cartier divisor 
$\sD_{\inf} \cong \pn{w^{n,\nu}}_{n,\nu}$.  

This yields a unique relative Cartier divisor 
$\sD = (\sD_{\et},\sD_{\inf}) \in \Fm{C}{\mdl}(R)$ 
satisfying the required conditions. 
An explicit procedure for finding the solution is described 
in Algorithm \ref{algorithm} below.  
\ePf 

\bAlg 
\label{algorithm}
The following algorithm 
produces a sequence 
$\pn{\sD_{\nmb}}_{\nmb \geq 1}$, 
with $\sD_{\nmb} = (\sD_{\nmb}^{\,\et},\sD_{\nmb}^{\,\inf})$, 
which can be thought of as an approximation to a solution 
of the system (\ref{Sy}) of Lemma \ref{Basis-determines_p}. 
For every $\nmb \geq 1$ consider the following finite subsystem (\ref{Sy(m)}) 
of (\ref{Sy}): 
\begin{align} 
\bigpair{\sD_{\et},\tha_{n,0,\sig,\nmb}}_{C,\mdl} &\,=\, t_{n,\sig,\nmb} 
& & \fral \,n \in N\pn{\mdl}, \;\;\fral \,1 \leq \sig \leq \sm(n) \notag \\ 
\label{Sy(m)} \mytag{Sy(m)} 
\bigpair{\sD_{\inf},\tha_{n,\nu,\sig,\mu}}_{C,\mdl} &\,=\, u_{n,\nu,\sig,\mu} 
& & \fral \,n \in N\pn{\mdl}, \;\;\fral \,1 \leq \sig \leq \sm(n) \\ 
& & & \fral 
\begin{array}{c} 
0 < \nu < n,  \notag \\ 
\pn{\nu,p} = 1, 
\end{array} 
\;\;\fral \;1 \leq \mu \leq \nmb 
\laurin \notag 
\end{align} 

\textbf{Step 1:} \'Etale part. \\ 
Using the decomposition (\ref{FZ}) of $\pn{\Fm{C}{\mdl}}_{\et}$ 
and the representation (\ref{Ln}) 
of an \'etale divisor $\sD_{\,\et}$ as a tuple of integers 
$\pn{\nbr_{n,\pntt}}_{n,\pntt}$, 
then due to (\ref{DL}) 
the \'etale part of (\ref{Sy(m)}) 
is expressed by the following system of equations: 
\begin{align} 
\label{Et(m)} \mytag{Et(m)} 
& \prod_{\pntt \in \Sm(n)} \tha_{n,0,\sig,\nmb}\pn{\pntt}^{\nbr_{n,\pntt}} 
\,=\, t_{n,\sig,\nmb} \\ 
& \hspace{0mm} \fral\,n \in N\pn{\mdl}, 
\;\;\fral \,1 \leq \sig \leq \sm(n) \notag 
\end{align} 
If the system (\ref{Et(m)}) has a solution $\pn{{}^{\nmb} \nbr_{n,\pntt}}_{n,\pntt}$ 
in $\st{-p^{\nmb},\ldots,p^{\nmb}}^{N(\mdl) \tms \Sm(n)}$ 
(since this is a finite set, it can be checked), 
the solution is unique by (PB2) of Definition \ref{pro-basis-ACHm}. 
Then by definition $\sD_{\nmb}^{\,\et}$ 
is the \'etale divisor corresponding to $\pn{{}^{\nmb} \nbr_{n,\pntt}}_{n,\pntt}$ 
with respect to the decomposition (\ref{FZ}). 
Otherwise set \,$\sD_{\nmb}^{\,\et} = 0$. 

\medskip 
\textbf{Step 2:} Infinitesimal part. \\ 
Using the decomposition (\ref{FW}) of $\pn{\Fm{C}{\mdl}}_{\inf}$ 
and the representations (\ref{Wn}) and (\ref{Wi}) 
of an infinitesimal divisor $\sD_{\inf}$ as a tuple 
$\pn{w^{n,\nu}}_{n,\nu}$, 
then due to (\ref{DE}) 
the infinitesimal part of (\ref{Sy(m)}) 
is expressed by the following system of equations: 
\begin{align} 
\label{In(m)} \mytag{In(m)} 
& \Expah\biggpn{ 
\sum_{\pntt \in \Sm(n)} w^{n,\nu}_{\,\pntt} 
	\cdot \bigbt{c_{n,\nu,\sig}^{\,\pntt}} \cdot \bigbt{c_{\mu}} }
\,=\, u_{n,\nu,\sig,\mu} \\ 
& \hspace{0mm} \fral\,n \in N\pn{\mdl}, 
\;\;\fral 
\begin{array}{c} 
0 < \nu < n,  \\ 
\pn{\nu,p} = 1, 
\end{array} 
\;\;\fral \,1 \leq \mu \leq \nmb, 
\;\;\fral \,1 \leq \sig \leq \sm(n) \notag 
\end{align} 
The system (\ref{In(m)}) does not always admit a solution 
$\pn{{}^{\nmb} w^{n,\nu}}_{n,\nu}$. 
But if we replace the indeterminate Witt vectors 
\,$w^{n,\nu} \in \Witt_{\resfld(n)}\bigpn{R_{\resfld(n)}}$ 
by families of Witt vectors 
\[ \bigpn{{}_{\nmr} w^{n,\nu}}_{\substack{\nmr \geq 1 \\ (\nmr,p) = 1}} 
\in \prod_{\substack{\nmr \geq 1 \\ \pn{\nmr,p}=1}} \Witt_{\resfld(n)}\bigpn{R_{\resfld(n)}} 
\laurink 
\] 
we can consider (\ref{In(m)}) as a special case of the following system 
\begin{align} 
\label{In'(m)} \mytag{In'(m)} 
& \Expah\biggpn{\sum_{\substack{\nmr \geq 1 \\ (\nmr,p) = 1}} 
\sum_{\pntt \in \Sm(n)} {}_{\nmr} w^{n,\nu}_{\,\pntt} 
	\cdot \bigbt{c_{n,\nu,\sig}^{\,\pntt}} \cdot \bigbt{\pn{c_{\mu}}^{\nmr}} }
\,=\, u_{n,\nu,\sig,\mu} \\ 
& \hspace{0mm} \fral\,n \in N\pn{\mdl}, 
\;\fral 
\begin{array}{c} 
0 < \nu < n,  \\ 
\pn{\nu,p} = 1, 
\end{array} 
\;\fral \,1 \leq \mu \leq \nmb, 
\;\fral \,1 \leq \sig \leq \sm(n) \notag 
\end{align} 
or written as a system of matrix equations 
\footnote{ $\Expah$ applied to a matrix is defined componentwise: 
$\Expah\bigpn{\pn{a_{ij}}_{ij}} := \bigpn{\Expah\pn{a_{ij}}}_{ij}$.} 
\begin{align*} 
& \Expah\biggpn{ \Bigpn{\bigbt{c_{n,\nu,\sig}^{\,\pntt}}}_{\sig,\pntt} 
\cdot \Bigpn{{}_{\nmr} w^{n,\nu}_{\,\pntt}}_{\pntt,\nmr} 
\cdot \Bigpn{\bigbt{\pn{c_{\mu}}^{\nmr}}}_{\nmr,\mu} } 
\,=\, \Bigpn{u_{n,\nu,\sig,\mu}}_{\sig,\mu} \\ 
& \hspace{0mm} \fral\,n \in N\pn{\mdl}, 
\;\fral \,0 < \nu < n,  \,\pn{\nu,p} = 1 
\end{align*} 
which is always solvable, due to the isomorphism \,$\Pi\eps$\, from (\ref{WL}). 

\noindent 
A solution $\bigpn{\pn{{}^{\nmb}_{\;\nmr} w^{n,\nu}}_{l}}_{n,\nu}$ of (\ref{In'(m)}) 
is found as follows: 
\begin{align} 
\label{Alg} \mytag{ALG} 
\bigpn{{}^{\nmb} r_{i}^{n,\nu,\sig}}_{0 \leq i < \nmb} 
  & \,:=\, \bigpn{u_{n,\nu,\sig,\mu}}_{1 \leq \mu \leq \nmb} 
  \Bigpn{\pn{c_{\mu}}^{i}}_{\substack{0 \leq i < \nmb \\ 1 \leq \mu \leq \nmb}}^{-1} 
  && \forall \,n, \nu, \sig \\ 
{}^{\nmb} u^{n,\nu,\sig} \pn{t} 
  & \,:=\, \sum_{0 \leq i < m} {}^{\nmb} r_i^{n,\nu,\sig} \, t^i 
  && \forall \,n, \nu, \sig \notag \\ 
\bigpn{{}^{\nmb}_{\;\nmr} v^{n,\nu,\sig}} 
  _{\substack{\nmr \geq 1 \\ \pn{\nmr,p}=1}} 
  & \,:=\, \Pi\eps^{-1} \bigpn{{}^{\nmb} u^{n,\nu,\sig}\pn{t}} 
  && \forall \,n, \nu, \sig \notag \\ 
\bigpn{{}^{\nmb}_{\;\nmr} w^{n,\nu}_{\,\pntt}}_{\pntt \in \Sm(n)}
  & \,:=\, \bigpn{{}^{\nmb}_{\;\nmr} v^{n,\nu,\sig}}_{1 \leq \sig \leq \sm(n)} 
  \,\Bigpn{\bigbt{c_{n,\nu,\sig}^{\,\pntt}}}
  _{\substack{\pntt \in \Sm(n) \\ 1 \leq \sig \leq \sm(n)}}^{-1} 
  && \forall \,n, \nu, \pntt, \nmr \notag 
\end{align} 
Then by definition $\sD^{\,\inf}_{\nmb}$ is the tuple 
$\pn{{}^{\nmb}_{\;1} w^{n,\nu}}_{n,\nu}$ (with $\nmr = 1$). 
This tuple corresponds to a relative Cartier divisor in $\Divf_C^{S,0}\pn{R}$ 
(also denoted by $\sD^{\,\inf}_{\nmb}$): 
if \,${}^{\nmb}_{\;1} w^{n,\nu}_{\,\pntt} \in \Wcfl{r\pn{\nu,n}}(R)$\, 
for all $n, \nu, \pntt$, 
then $\sD^{\,\inf}_{\nmb} \in \pn{\Fm{C}{\mdl}}_{\inf}\pn{R}$ 
due to the decomposition (\ref{FW}); 
for the general case cf.\ \Point \ref{Pairing_curve_lim}.  
Finally $\sD^{\,\inf}_{\nmb}$ is a solution of (\ref{In(m)}) 
if in addition \,${}^{\nmb}_{\;\nmr} w^{n,\nu} = 0$\, 
for all $\nmr \geq 2$ and all $n, \nu$. 
\eAlg 

\bThm 
\label{System-Solve}
The system {\rm (\ref{Sy})} from Lemma \ref{Basis-determines_p} 
is solvable if and only if 
the output $\pn{\sD_{\nmb}}_{\nmb \geq 1}$ of Algorithm \ref{algorithm} 
satisfies the following conditions: 
\begin{align*} 
{\rm (SC1)} \hspace{6mm}  & \exists \,\alp \geq 1 \;\; \forall \,m \geq \alp: 
\hspace{1mm} \sD_m \,\textrm{ is a solution of {\rm (\ref{Sy(m)})}} \laurink \\ 
{\rm (SC2)} \hspace{6mm}  & \pn{\sD_{\nmb}}_{\nmb \geq 1} \;\textrm{ becomes stationary} \laurin 
\end{align*} 
In this case the solution of {\rm (\ref{Sy})} is given by the limit of 
$\pn{\sD_{\nmb}}_{\nmb \geq 1}$. 
\eThm 

\bPf 
$(\Lra)$ results from Lemma \ref{Basis-determines_p} and its proof; 
the deeper reason is the structure of $\Fm{C}{\mdl}$, cf.\ (\ref{TR}). 

$(\Lla)$ Conditions (SC1) and (SC2) imply that (\ref{Sy}) 
is equivalent to a solvable finite subsystem (\ref{Sy(m)}). 
\ePf 

\vspace{\vs} 

The next question is: 
If (\ref{Sy}) is solvable 
and $\pn{\sD_{\nmb}}_{\nmb}$ the output of Algorithm \ref{algorithm}, 
after how many elements does this sequence become stationary? 

\bDef 
\label{Def-bounds}
Let $R$ be a finite $k$-ring, and set 
\[ \nil(R) := \min \st{\,n \;|\; r^n = 0 \,,\hspace{3mm} \forall \,r \in \Nil(R)} 
\laurin 
\] 
For $w = \pn{w_i}_{i \geq 0} \in \Wittc\pn{R}$ define 
\[ \lgth(w) := \min \st{\,l \;|\; w_i = 0 \,,\hspace{3mm} \forall \,i \geq l} 
\laurin 
\] 
There is 
an increasing 
function \,$\dime: \Nat \ra \Nat$, 
depending only on \,$\nil(R)$, s.t.\ 
\bMyEqn \label{LD} \mytag{LD} 
\forall \,w \in \Wittc(R): \hspace{5mm} 
\lgth(w) \leq l 
\hspace{5mm} \Lra \hspace{5mm} 
\deg \,\Expaht\pn{w,t} \leq \dime(l) 
\eMyEqn 
where \,$\Expaht\pn{w,t} \in R[t]$. 

In the situation of the Rigidity Lemma \ref{Basis-determines_p}, 
if \,$\sD_{\inf} \in \pn{\Fm{C}{\mdl}}_{\inf}(R)$\, corresponds to 
\,$\pn{w^{n,\nu}}_{n,\nu} 
\in \dsum_{n,\nu} \;\Wcfl{r(n,\nu)}_{\resfld(n)}\pn{R_{\resfld(n)}}$\, 
with respect to (\ref{FW}), 
then define 
\[ \lgth\pn{\sD_{\inf}} 
:= \max \bigst{\lgth\pn{w^{n,\nu}} \;\big|\; n, \nu} 
\laurin 
\] 
(A priori this definition depends on the choice of a basis of $\Ump{C}{\mdl}$). \\ 
If \,$\sD_{\et} \in \pn{\Fm{C}{\mdl}}_{\et}(R)$\, corresponds to 
\,$\pn{\nbr_{n,\pntt}}_{n,\pntt} \in \dsum_{n,\pntt} \Zint$\, 
with respect to (\ref{FZ}), 
then define 
\[ \lgth\pn{\sD_{\et}} 
:= \max \bigst{\log_{p} \pn{\vrt{\nbr_{n,\pntt}}} \;\big|\; n,\pntt} 
\laurin 
\] 
For \,$\sD \in \Fm{C}{\mdl}\pn{R}$\, define 
\[ \lgth\pn{\sD} 
:= \max \bigst{\lgth\pn{\sD_{\et}}, \lgth\pn{\sD_{\inf}}} 
\laurin 
\] 
\eDef 

\bLem 
\label{bounds}
Let $R$ be a finite $k$-ring. 
There exists an increasing function \,$\mty: \Nat \ra \Nat$, 
depending only on \,$\nil(R)$, 
such that for all systems {\rm (\ref{Sy})}: 

If \,$\sD$\, is a solution of system {\rm (\ref{Sy})}, 
then the output $\pn{\sD_{\nmb}}_{\nmb \geq 1}$ of Algorithm \ref{algorithm} 
becomes stationary after \,$\mty\bigpn{\lgth\pn{\sD}}$\, elements: 
\bMyEqn \label{LM} \mytag{LM} 
\lgth(\sD) \leq l 
\hspace{5mm} \Lra \hspace{5mm} 
\forall \,\nmb \geq \mty(l): \hspace{2mm} \sD_{\nmb} = \sD 
\laurin 
\eMyEqn 
\eLem  

\bPf 
If $\sD_{\et} \cong \pn{\nbr_{n,\pntt}}_{n,\pntt}$ 
is the solution of $(\textrm{\ref{Sy}})_{\et}$, 
then Algorithm \ref{algorithm} applied to $(\textrm{\ref{Sy(m)}})_{\et}$ 
finds this solution 
if \,$p^{m} \geq \vrt{\nbr_{n,\pntt}}$\, for all $n,\pntt$, 
i.e.\ if \,$m \geq \lgth\pn{\sD_{\et}}$. 

The solution $\sD_{\inf}$ of $(\textrm{\ref{Sy}})_{\inf}$ 
is determined by the equations (\ref{WU}). 
Then by (\ref{LD}) the degree of the $u^{n,\nu,\sig}(t)$ 
is bounded by a function in \,$\lgth(\sD_{\inf})$. 
On the other hand, looking at the procedure (\ref{Alg}), 
if \,$\nmb > \deg u^{n,\nu,\sig}(t)$\, for all ${n,\nu,\sig}$, 
then ${}^{\nmb} u^{n,\nu,\sig}(t) = u^{n,\nu,\sig}(t)$, 
since a polynomial of degree $< \nmb$ 
is determined by its values at $\nmb$ distinct points. 
It follows $\sD_{\nmb}^{\,\inf} = \sD_{\inf}$. 

Thus there is an increasing function \,$\mty: \Nat \ra \Nat$, 
depending only on the increasing function 
$\dime: \Nat \ra \Nat$ from (\ref{LD}), 
with the required property. 
\ePf 

\bCor 
\label{finite subsystem}
Suppose $\sD$ is a solution of {\rm (\ref{Sy})}. 
Let $M \subset \Nat\setminus\st{0}$ be a subset 
and {\rm (Sy($M$))} be the subsystem of {\rm (\ref{Sy})} 
for $\mu \in M$. 
Then Algorithm \ref{algorithm} applied to {\rm (Sy($M$))} 
produces the solution \,$\sD$\, 
if \;$\# M \geq \mty\bigpn{\lgth(\sD)}$. 
\eCor 

\bRmk 
\label{non-irreducible curves}
It is a straightforward exercise to extend the results of this Subsection \ref{rigidity} 
to the case of smooth proper non-irreducible curves, 
and hence to the case of arbitrary curves in $X$, 
due to Definitions \ref{pro-basis-C-in-X} and \ref{Pairing_curve}. 
(Smooth irreducible curves would be sufficient for this work, 
but it keeps formulations simpler if we may consider arbitrary curves in $X$.) 
\eRmk


\subsection{Skeleton Theorem} 
\label{skeleton theorem}

Let $X$ be a suitable projective variety of dimension $d$ over $\fld$, 
a finite field or an algebraic closure of a finite field of characteristic $p$, 
or an algebraically closed field  of characteristic $0$. 
Let $\mdl$ be a modulus for $X$, as in Section \ref{ChowMod}. 
Let $S$ be the support of an effective Cartier divisor $\geq \mdl_X$ on $X$.

\bThm[Skeleton Theorem] 
\label{SkelThm}
The canonical homomorphism \\ 
\;$ \skelm{X}{\mdl,S}: \Fmor{X}{\mdl} \lra \varprojlim \Emo{C}{\mdl_{\Cgst}} 
$\; 
from Lemma \ref{Fm-to-limEm} is 
an isomorphism. 
\eThm 

\bPf
\label{Idea}
The following proof makes use of a number of lemmata 
that will be provided afterwards in the rest of this section. 
By means of a Galois descent, 
which is possible due to \Point \ref{descent_direct-sys}, 
it is sufficient to work over an algebraic closure $k$ of the base field. 

Injectivity of \,$\skelm{X}{\mdl,S}$\, 
follows from the fact that over an algebraically closed base field  
a relative divisor on $X$ with support in $\vrt{\mdl_X}$ 
is detected on a dense subset of $\vrt{\mdl_X}$, 
and the curves $C$ in $X$ meeting $\vrt{\mdl_X}$ properly 
obviously cover $\vrt{\mdl_X}$. 

It remains to show surjectivity. 
Let \,$\pn{\sD_C}_C \in \varprojlim \Emo{C}{\mdl_{\Cgst}}$\, 
be a compatible system of relative Cartier divisors on curves in $X$. 
The task is to construct a relative Cartier divisor 
$\sD_X \in \Divf_X^{\vrt{\mdl \cut X}}$ 
such that 
\bMyEqn \label{SP} \mytag{SP} 
\sD_X \cut C = \sD_C \hspace{20mm} \forall \,C \subset X, \; C \iscp \mdl 
\laurin 
\eMyEqn  
Then $\sD_X \in \Fmor{X}{\mdl}$ by Corollary \ref{Fm=Fm^CH} 
and Definition \ref{DefFch}. 
%
To this end, we proceed as follows: 

\textbf{Step 1:} 
Consider relative curves \,$\rC \ra \Base$\, over irreducible varieties $\Base$, 
such that the fibres $\rC(b)$ over $b \in \Base$ 
are curves in $X$ that intersect $\vrt{\mdl_X}$ properly. 
Let $\mfn: \rC \ra X$ be the canonical map. 
There exists such a family $\rC \ra \Base$ 
of curves in $X$ that generates $\ACHm{X}{\mdl}$: 
\[ \ACHm{X}{\mdl}  =  
\sum_{b \in \Base} {\iot{\rC(b)}}_* \ACHm{\rC(b)}{\mdl_{\rCbgst}} 
\] 
($\see$Lemma \ref{ACH generated by B}). 
By duality, the compatible system 
$\pn{\sD_C}_C \in \varprojlim \Emo{C}{\mdl_{\Cgst}}$ 
is then completely determined by the subsystem $\pn{\sD_{\rC(b)}}_{b \in \Base}$ 
($\see$Lemma \ref{compSys-determined}). 
In order that $\sD_X \in \Divf_X^{\vrt{\mdl}}$ satisfies (\ref{SP}) 
it is thus sufficient that it satisfies 
\bMyEqn \label{SPB} \mytag{SPB} 
\sD_X \cut \rC(b) = \sD_{\rC(b)} \hspace{20mm} \forall \,b \in \Base 
\laurin 
\eMyEqn 

\textbf{Step 2:} 
We construct a relative Cartier divisor $\sD_{\rC} \in \Divf_{\rC}^{\vrt{\mdl \cut \rC}}$ 
mapping to $\pn{\sD_{\rC(b)}}_{b \in \Base}$ 
($\see$Lemma \ref{locRelDiv}): 
We find an open dense $U \subset \Base$ 
and a relative divisor $\sD_{\rC \tms_{\Base} U} 
\in \Divf_{\rC \tms_{\Base} U}^{\vrt{\mdl \cut \rC} \tms_{\Base} U}$ 
on $\rC \tms_{\Base} U$ 
mapping to $\pn{\sD_{\rC(b)}}_{b \in U}$. 
Then we show that $\sD_{\rC \tms_{\Base} U}$ 
extends to a relative divisor $\sD_{\rC}$ on $\rC$. 
For the construction of $\sD_{\rC \tms_{\Base} U}$ 
we consider the generic fibre $\rC(\gnc)$ of $\mfn$ 
(here $\gnc$ is the generic point of $\Base$). 
Let $\ol{\fld(\gnc)}$ be an algebraic closure of the function field $\fld(\gnc)$ 
of $\Base$, 
denote by $\gnccl = \bigpn{\Spec \ol{\fld(\gnc)} \lra \Base}$ 
the corresponding geometric point 
and let $\rC(\gnccl) = \rC(\gnc) \tens_{\fld(\gnc)} \ol{\fld(\gnc)}$ 
be the corresponding geometric fibre. 
On $\rC(\gnccl)$ 
we apply rigidity from Subsection \ref{rigidity} and a reciprocity trick 
(explained in the following): 

Let \,$\idx$\, be a multi-index, 
$\idx = \pn{n,\nu,\sig}$ for $\chr(k) = 0$ and 
$\idx = \pn{n,\nu,\sig, \mu}$ for $\chr(k) = p > 0$. 
Let $\st{\tha_{\idx}}_{\idx} \subset \sK_X^*$ 
be a family of functions on $X$ that yields 
a (pro-)basis of $\ACHmp{\rC(\gnccl)}{\mdl_{\rCgncclgst}}$, 
and define 
$\Tha_{\idx} := \dv\pn{\tha_{\idx}}_X$. 
By induction on the dimension of $X$ 
we may assume the Skeleton Theorem to hold 
for suitable varieties of dimension $< \dim X$. 
Hence there exist relative divisors $\sD_{\Tha_{\idx}}$ 
satisfying (\ref{SP}) mutatis mutandis. 
(The precise construction of $\sD_{\Tha_{\idx}}$ will be done in 
the proof of Lemma \ref{locRelDiv}, Step 1, and Claim \ref{compSystemHi}.) 
By rigidity, any element of $\Divf_{\rC(\gnccl)}^{\vrt{\mdl \cut \rC(\gnccl)}}$ 
is uniquely determined by its values on the basis 
$\st{\tha_{\idx}}_{\idx}$ 
with respect to the pairing $\pair{\llul,\lull}_{\rC(\gnccl),\mdl_{\rCgncclgst}}$. 
For every $b \in \Base$ 
the compatible system 
$\pn{\sD_C}_C \in \varprojlim \Emo{C}{\mdl_{\Cgst}}$ 
satisfies 
\,$ \bigpair{\sD_{\rC(b)}, \Tha_{\idx} \cut \rC(b)}_{\rC(b),\mdl} 
= \bigpair{\sD_{\Tha_{\idx}}, 
   \Tha_{\idx} \cut \rC(b)}_{\vrt{\Tha_{\idx}},\mdl} 
$\, 
by reciprocity. 
Then we show that 
there exists a unique relative Cartier divisor $\sD_{\rC\pn{\gnccl}}$ 
on $\rC\pn{\gnccl}$ subject to the conditions 
\[ \bigpair{\sD_{\rC\pn{\gnccl}}, 
   \Tha_{\idx} \cut \rC(\gnccl)}_{\rC(\gnccl),\mdl} 
= \bigpair{\sD_{\Tha_{\idx}}, 
   \Tha_{\idx} \cut \rC\pn{\gnccl}}_{\vrt{\Tha_{\idx}},\mdl} 
\] 
where the RHS 
is the image in $\GmS{k(\gnccl)}$ 
of the unique $R_{k(\gnc)}$-valued point of 
$\GmS{k(\gnc)}$ that specializes to 
$\bigpair{\sD_{\Tha_{\idx}}, 
  \Tha_{\idx} \cut \rC(b)}_{\vrt{\Tha_{\idx}},\mdl} 
  \in \GmS{k(b)}\pn{R_{k(b)}}$ 
for all $b \in U$, 
and $U \subset \Base$ is dense open with 
$\mfn^*\tha_{\idx} 
\in \sO_{\rC \tms_{\Base} U, \vrt{\mdl \cut \rC} \tms_{\Base} U}^*$ 
for all $n,\nu,\sig$ and all $\mu \leq m$, $m$ sufficiently large 
(for $\chr(k) = p > 0$). 
By construction $\sD_{\rC\pn{\gnccl}}$ comes from 
a relative divisor $\sD_{\rC \tms_{\Base} U}$ 
on $\rC \tms_{\Base} U$, 
and $\sD_{\rC \tms_{\Base} U}$ 
extends to a relative divisor $\sD_{\rC}$ on $\rC$ 
by Lemma \ref{locRelDiv}, Step 3. 

\textbf{Step 3:} 
We choose a subscheme $\SubBase \subset \Base$ such that 
$\rC_{\SubBase} := \rC \tms_{\Base} \SubBase$ 
is suitable of dimension $\dim X$ 
and the canonical map \,$\mfn_{\SubBase}: \rC_{\SubBase} \lra X$\, 
as well as the restriction \,$\mfn_{\SubBase}^{-1}\vrt{\mdl} \lra \vrt{\mdl}$\, 
are birational. 
Then the induced map 
\[ \mfn_{\SubBase}^*: \Divor{X} \iso \Divor{\rC_{\SubBase}}
\] 
is an isomorphism, 
since $\Divf_X$ sits inside an exact sequence 
\[ 0 \lra \ul{k}^* \lra \ul{\sK}_X^* \lra \Divor{X} \lra \Picorf{X} \lra 0 
\] 
and all group-sheaves in this sequence (except $\Divor{X}$ a priori) 
are birational invariants among suitable proper schemes $X$ 
($\see$Remark \ref{birInvariant}), 
thus the same is true for $\Divf_X$. 
Then the relative divisor $\sD_{\rC} \cut \rC_{\SubBase}$ 
yields via $\pn{\mfn_{\SubBase}^*}^{-1}$ 
a unique element $\sD_X \in \Divf_X^{\vrt{\mdl}}$. 
The relative divisor $\sD_X$ on $X$ obtained in this way 
is independent of the choice of the subscheme $\SubBase$, 
i.e.\ independent of the choice of a rational section 
$X \dra \rC$ of the map $\mfn: \rC \ra X$ 
($\see$Lemma \ref{Secn-independent}). 
Then 
$ \sD_{\rC} = \mfn^* \sD_X 
\laurink 
$ 
hence $\sD_X$ satisfies (\ref{SPB}): 
$\sD_{\rC}$ has this property by construction, 
and the inclusion $\rC(b) \inj X$ factors through $\mfn: \rC \ra X$ 
for every $b \in \Base$. 
\ePf 

\vspace{\vs} 
The rest of this section is devoted to prove the lemmata 
that were used in the proof of Theorem \ref{SkelThm}. 

In the following we will work over an algebraic closure of the base field, 
unless stated otherwise. 
Let $\Base$ be an irreducible space of curves in $X$, 
denote by $\gnc$ the generic point of $\Base$. 
Let $\strl: \rC \ra \Base$ be the universal curve over $\Base$ 
and \,$\mfn: \rC \ra X$\, the canonical map, 
i.e.\ $\rC$ is given by the incidence variety 
\[ \xymatrix{ 
& \, \rC \ar[dl]_-{\mfn} \ar[dr]^-{\strl} & & 
\hspace{-22mm} = \bigst{\pn{x,b} \in X \tms \Base \;\big|\; x \in b} \\ 
X & & \Base \laurin & \\ 
} 
\] 


\bPnt 
\label{Base conditions}
We will assume that $\Base$ has the following properties: \\ 
\begin{tabular}{rl} 
(T0) & The generic fibre $\rC(\gnc)$ is smooth. 
(Possible since $X$ is normal.) \\ 
(T1) & $\rC(b)$ meets $\mdl$ properly \hspace{4mm} $\forall \,b \in \Base$. \\ 
(T2) & Every truncated curve 
$T \in \Trc{X}{\mdl}{\deg X \cdot \pn{\deg\mdl}^{\dim X -1}}$ \\ 
 & (see Not.~\ref{truncated curve}) 
 is a subspace of $\rC(b)$ for some element $b \in \Base$. 
\end{tabular} 
\ePnt 

\bLem 
\label{Base property}
If $\sL$ is an ample line bundle on $X$, 
then for some $n \in \Nat$ there is an open dense subscheme $\Base$ 
of the $\pn{\dim X - 1}$-fold complete linear system $\vrt{\sL^{\tens n}}^{d-1}$, 
defined over the base field, 
satisfying conditions {\rm (T0), (T1), (T2)} of \ref{Base conditions}. 
\eLem 

\bPf 
A very ample line bundle $\sA$ on $X$ 
defines an embedding $\iot{\sA}$ into the projective space 
$P := \Prj\bigpn{\H^0\pn{X,\sA}^{\vee}}$ 
such that $\sA = \iot{\sA}^*\sO_P(1)$. 
Viewing $X$ as a subspace of $P$, 
every curve $C \in \vrt{\sA}^{d-1} = \Prj\bigpn{\H^0\pn{X,\sA}}^{d-1}$ 
arises as a multiple hyperplane section 
$H_1 \cap \ldots \cap H_{d-1} \cap X$ for some $H_i \in \vrt{\sO_P(1)}$. 
By Bertini's Theorem ($\seecite$\cite[II, Thm.~8.18]{H}), 
for hyperplanes $H \in \vrt{\sO_P(1)}$ it is an open dense property 
that $H \cap X$ is smooth. 
Thus by induction any dense open subscheme $B$ of $\vrt{\sA}^{d-1}$ 
satisfies (T0). 

Since $\sL$ is ample, 
some power of $\sL$ is very ample 
and thus defines an embedding of $X$ 
into a projective space $\Prj^m$. 
Replacing $\sL$ by this power of $\sL$ 
we may assume that $\sL$ comes from $\sO_{\Prj^m}(1)$. 
The approximation property (T2) for finite subschemes of bounded length 
by sections of $\sO_{\Prj^m}(n)$ is then evident for $n \gg 0$. 

Property (T1), i.e.\ meeting $\mdl$ properly, is obviously an open property 
among the elements $C \in \vrt{\sL^{\tens n}}$. 
The subspace $\Base \subset \vrt{\sL^{\tens n}}^{d-1}$ of those elements 
satisfying (T1) is non-empty 
and $\vrt{\sL^{\tens n}}^{d-1} = \Prj\bigpn{\H^0\pn{X,\sL^{\tens n}}}^{d-1}$ 
is irreducible, 
so $\Base$ is dense in $\vrt{\sL^{\tens n}}^{d-1}$. 

Now assume $\Base' \subset \vrt{\sL^{\tens n}}^{d-1}$ is an open dense subscheme 
satisfying conditions (T0), (T1), (T2). 
As $\Base'$ is defined over a finite extension of the base field, 
and by the arguments above, 
every Galois conjugate of $\Base'$ will satisfy (T0), (T1), (T2) as well. 
Therefore, defining $\Base$ as the (finite) union of all Galois conjugates of $\Base'$ 
yields a subscheme defined over the base field. 
\ePf 


\bLem 
\label{ACH generated by B}
The affine Chow group with modulus $\ACHm{X}{\mdl}$ 
is generated by the groups $\ACHm{\rC(b)}{\mdl_{\rCbgst}}$ 
for $b \in \Base$: 
\[ \ACHm{X}{\mdl}  =  
\sum_{b \in \Base} {\iot{\rC(b)}}_* \ACHm{\rC(b)}{\mdl_{\rCbgst}} 
\laurin 
\] 
\eLem 

\bPf 
Property (T1) of $\Base$ implies that the expression 
$\ACHm{\rC(b)}{\mdl_{\rCbgst}}$ is well-defined for every $b \in \Base$. 
By property (T2) of $\Base$ 
we find for any $T \in \Trc{X}{\mdl}{\deg\mdl}$ 
a curve $b \in \Base$ such that $T \subset \rC(b)$ and hence 
\[ {\iot{C_{T}}}_* \ACHm{C_{T}}{\mdl_{\CTgst}} 
\,\subset\, {\iot{\rC(b)}}_* \ACHm{\rC(b)}{\mdl_{\rCbgst}} 
\] 
for any curve $C_{T}$ in $X$ with $C_{T} \tms_X \bt{\mdl} = T$. 
Then the assertion follows from Theorem \ref{ACH generated by truncations}. 
\ePf 

\bLem
\label{compSys-determined}
Any skeleton divisor $\pn{\sD_C}_C \in \varprojlim \Emo{C}{\mdl_{\Cgst}}(R)$, 
where $R$ is a finite $k$-ring, 
is uniquely determined by the subsystem $\pn{\sD_{\rC(b)}}_{b \in \Base}$. 
\eLem 

\bPf 
For any curve $C \subset X$ we have the duality pairing 
$\pair{\llul,\lull}_{C,\mdl_{\Cgst}}: \Emo{C}{\mdl_{\Cgst}} \tms \ACHm{C}{\mdl_{\Cgst}} \lra \Gm$, 
which is perfect. 
Therefore a relative Cartier divisor $\sD_C \in \Emo{C}{\mdl_{\Cgst}}$ 
is uniquely determined by the values 
\[ \bigst{\bigpair{\sD_C,\alp}_{C,\mdl_{\Cgst}} \;\big|\; \alp \in \ACHm{C}{\mdl_{\Cgst}}} 
\laurink 
\] 
cf.\ Lemma \ref{Basis-determines_p}. 
By Lemma \ref{ACH generated by B} 
for every $\alp \in \ACHm{C}{\mdl_{\Cgst}}$ we find curves $b_1,\ldots, b_r \in \Base$ 
and elements $\bet_i \in \ACHm{\rC(b_i)}{\mdl_{\rCbigst}}$, $i = 1,\ldots, r$ 
such that 
${\iot{C}}_* \alp = \sum_{i = 1}^{r} {\iot{\rC(b_i)}}_* \bet_i \in \ACHm{X}{\mdl}$. 
Then 
\[ \bigpair{\sD_C,\alp}_{C,\mdl_{\Cgst}} 
= \sum_{i = 1}^{r} \bigpair{\sD_{\rC(b_i)},\bet_i}_{\rC(b_i),\mdl_{\rCbigst}} 
\] 
due to the compatibility condition of 
$\Emo{C \cup \bigcup_i \rC(b_i)}{\mdl_{C \cup \bigcup_i \rC(b_i)}}$ 
from Remark \ref{Char_Em} 
and the fact that 
$\pn{\sD_C}_C \in \varprojlim \Emo{C}{\mdl_{\Cgst}}(R)$ 
forms a compatible system. 
\ePf 

\bLem
\label{locRelDiv}
For every skeleton divisor $\pn{\sD_C}_C \in \varprojlim \Emo{C}{\mdl_{\Cgst}}(R)$, 
where $R$ is a finite $k$-ring, 
there is a unique relative Cartier divisor 
\,$\sD_{\rC} \in \Divf_{\rC}^{\vrt{\mdl \cut \rC}}(R)$ 
on \,$\rC \ra \Base$ with the property 
\,$\sD_{\rC} \cut \rC(b) = \sD_{\rC(b)}$\, for all $b \in \Base$. 
\eLem 

\bPf 
Let $\pn{\sD_C}_C \in \varprojlim \Emo{C}{\mdl_{\Cgst}}(R)$ 
be a skeleton relative divisor on $X$. 
We seek a relative Cartier divisor 
$\sD_{\rC} \in \Divf_{\rC/\Base}^{\vrt{\mdl \cut \rC}}$, 
such that $\sD_{\rC} \cut \rC(b) = \sD_{\rC(b)}$ 
for all $b \in \Base$. 
Pulling-back to the generic fibre $\rC(\gnc)$, 
where $\gnc$ is the generic point of $\Base$, 
yields a map 
$\Divf_{\rC}^{\vrt{\mdl \cut \rC}} \lra \Divf_{\rC(\gnc)}^{\vrt{\mdl \cut \rC(\gnc)}}$. 
We change the base $\fld(\gnc)$ to an algebraic closure $\ol{\fld(\gnc)}$: 
let \,$\gnccl = \bigpn{\Spec \ol{\fld(\gnc)} \ra \Base}$, 
and \,$\rC(\gnccl) = \rC(\gnc) \tens_{\fld(\gnc)} \ol{\fld(\gnc)}$. 
\[ \xymatrix{ 
\rC(\gnccl) \ar[d] \ar[r] \ar@{}[dr]|{\square} & 
\rC(\gnc) \ar[d] \ar[r] \ar@{}[dr]|{\square} & 
\rC \ar[d] \ar[r]^-{\mfn} & X \\ 
\gnccl \ar[r] & \gnc \ar[r] & \Base & 
}
\] 
Then we can apply rigidity from Subsection \ref{rigidity}: 
a relative Cartier divisor $\sD_{\Crv} \in \Fmo{C}{\mdl_{\Cgst}}$ on a curve $\Crv$ 
over an algebraically closed field 
can be recovered from the values 
$\bigpair{\sD_{\Crv}, \alp_{n,\nu,\sig,\mu}^{\Crv}}_{\Crv,\mdl_{\Crvgst}}$ 
of the pairing $\pair{\llul,\lull}_{\Crv,\mdl_{\Crvgst}}$ 
on a pro-basis $\st{\alp_{n,\nu,\sig,\mu}^{\Crv}}_{n,\nu,\sig,\mu}$ 
of $\ACHmp{\Crv}{\mdl_{\Crvgst}}$. 
This will be done in Step 2. 

Theorem \ref{SkelThm} is trivially true if $X$ is a curve. 
By induction on the dimension of $X$, 
we may assume that Theorem \ref{SkelThm} is true 
for suitable projective varieties of dimension $< d := \dim X$. 

\medskip 
\textbf{Step 1:} Using Reciprocity. \\ 
Let $\tha \in \sK_X^*$, 
set $\Tha := \dv\pn{\tha} = \sum_i n_i \,H_i$. 
Let $\sig_i: \Ht_i \ra H_i$ be suitabilizations of the $H_i$, 
let $\vrt{\That} := \coprod_i \Ht_i$ 
and $\That = \sum_i n_i \,\Ht_i$, 
considered as a divisor on $\vrt{\That} \tms \Afn^1$.  
Then for every $i$ there is a relative Cartier divisor 
$\sD_{\Ht_i}$ on $\Ht_i$ 
such that \,$\sD_{\Ht_i} \cut \Ct = \sD_C$\, 
for all curves $C$ in $H_i$ 
(note that $\sD_C \in \Emo{C}{\mdl_{\Cgst}} \subset \Emo{\Ct}{\mdl}$, 
and as $\Ht_i \ra H_i$ is birational 
we have a factorization $\Ct \ra C_{\Ht_i} \ra C$, 
where $C_{\Ht_i}$ denotes the proper transform of $C$ in $\Ht_i$). 
For $\dim X = 2$ this is trivial since the $H_i$ are curves. 
For $\dim X > 2$ this holds by induction hypothesis 
(we apply Theorem \ref{SkelThm} to the case 
$X = \Ht_i$ and $S = \vrt{\mdl_{\Ht_i} + E}$, 
where $E$ is the exceptional divisor of $\sig_i: \Ht_i \ra H_i$). 

\bClm 
\label{compSystemHi}
The relative divisors $\sD_{\Ht_i}$ 
form a compatible system, i.e.\ 
\[ \bigpair{\sD_{\Ht_i},\alpt}_{\Ht_i,\mdl} 
= \bigpair{\sD_{\Ht_j},\betat}_{\Ht_j,\mdl} 
\] 
for all \,$\alpt \in \ACHm{\Ht_i}{\mdl}$, 
$\betat \in \ACHm{\Ht_j}{\mdl}$ such that 
$\pn{\iot{\Ht_i}}_*\alpt 
= \pn{\iot{\Ht_j}}_*\betat \in \ACHm{X}{\mdl}$. 
This compatible system will be denoted by $\sD_{\That}$. 
\eClm

\bBf 
Let $\alp = {\sig_i}_* \alpt \in \ACHm{H_i}{\mdl}$, 
$\beta = {\sig_j}_* \betat \in \ACHm{H_j}{\mdl}$, 
and let $C \subset H_i$ and $Z \subset H_j$ be curves such that 
$\alp \in \ACHm{C}{\mdl}$ and $\beta \in \ACHm{Z}{\mdl}$. 
Denote by $\nu^*\alp \in \ACHm{\Ct}{\mdl}$ a lift of $\alp$, 
by $\nu^*\beta \in \ACHm{\Zt}{\mdl}$ a lift of $\beta$. 
By definition of the pairing in higher dimension 
($\see$Proposition \ref{Pairing_X}) 
we have 
\begin{align*} 
\bigpair{\sD_{\Ht_i}, \alpt}_{\Ht_i,\mdl} 
 = \bigpair{\sD_{\Ht_i} \cut C_{\Ht_i}, \alpt}_{C_{\Ht_i},\mdl} 
& = \bigpair{\sD_{\Ht_i} \cut \Ct, \nu^*\alp}_{\Ct,\mdl} 
\laurin 
\end{align*} 
Due to \,$\sD_{\Ht_i} \cut \Ct = \sD_C$\, 
and by definition of the pairing for curves 
($\see$Proposition \ref{Pairing_curve}) 
it holds 
\begin{align*} 
\bigpair{\sD_{\Ht_i} \cut \Ct, \nu^*\alp}_{\Ct,\mdl} 
& = \bigpair{\sD_C, \alp}_{C,\mdl} 
\laurin 
\end{align*} 
Since ${\iot{C}}_*\alp = {\iot{Z}}_*\bet \in \ACHm{X}{\mdl}$ 
we obtain by compatibility of skeleton divisors ($\see$Remark \ref{Char_Em}) 
\begin{align*} 
\bigpair{\sD_C, \alp}_{C,\mdl} 
& = \bigpair{\sD_Z, \beta}_{Z,\mdl} 
\laurin 
\end{align*} 
Finally 
\begin{align*} 
\bigpair{\sD_Z, \beta}_{Z,\mdl} 
& = \bigpair{\sD_{\Ht_j}, \betat}_{\Ht_j,\mdl} 
\end{align*} 
by the same arguments as for 
\,$\pair{\sD_C, \alp}_{C,\mdl} = \pair{\sD_{\Ht_i}, \alpt}_{\Ht_i,\mdl}$\,
above. 
\eBf 

\medskip 

By Lemma \ref{Base property} we may assume that the curves $\rC(b)$ 
for $b \in \Base$ are complete intersection curves. 
According to Lemma \ref{tameLift} and induction, 
for any $b \in \Base$ such that the curve $\rC(b)$ meets $\Tha$ properly 
we find $\gam_1,\ldots,\gam_{d-1} \in \sK_X^*$ 
with $\Gam_i := \dv\pn{\gam_i}_X$ 
such that $\Gam_1 \cdots \Gam_{d-1} \supset \rC(b)$ 
and \;$\dv\pn{\tha}_C = \dv\pn{\tha}_{\Gam_1 \cdots \Gam_{d-1}}$, 
and there is a reciprocity 
\bMyEqn 
\label{Rec} \mytag{Rec} 
\Tha \cut \rC(b) 
= \dv\pn{\tha}_{\rC(b)} 
= \dv\pn{\tha}_{\Gam_1 \cdots \Gam_{d-1}} 
= \dv\pn{\gam_i}_{\Tha \cdot \Gam_1 \cdots \overset{\vee} \Gam_i \cdots 
   \Gam_{d-1}} 
\laurin 
\eMyEqn 
If we denote \,$\rZ(b) := \Tha \cdot \Gam_1 \cdots \Gam_{d-2}$, 
then due to compatibility of skeleton divisors ($\see$Remark \ref{Char_Em}) 
and the formula for the pairing in terms of a rational map associated 
with a relative divisor ($\see$Proposition \ref{Pairing_X}), 
we have 
\begin{align} 
\label{PU} \mytag{PU} 
\bigpair{\sD_{\rC(b)}, \Tha \cut \rC(b)}_{\rC(b),\mdl} 
& =\, \bigpair{\sD_{\rC(b)}, \rZ(b) \cut \rC(b)}_{\rC(b),\mdl} \\ 
& =\, \bigpair{\sD_{\rZ(b)}, \rZ(b) \cut \rC(b)}_{\vrt{\rZ(b)},\mdl} \notag \\ 
& =\, \bigpair{\sD_{\rZ(b)}, \wt{\rZ(b)} \cut \rC(b)}_{\vrt{\wt{\rZ(b)}},\mdl} \notag \\ 
& =\, \bigpair{\sD_{\rZ(b)}, \rZ(b)_{\That} \cut \rC(b)}_{\vrt{\rZ(b)_{\That}},\mdl} \notag \\ 
& =\, \bigpair{\sD_{\That}, \That \cut \rC(b)}_{\vrt{\That},\mdl} \notag \\ 
& =\, \rmp{\sD_{\That}} \bigpn{\That \cut \rC(b)}^{-1} \notag \\ 
& =\, \rmp{\sD_{\That}}_{\Base} \bigpn{\That \cut \rC(b)}^{-1} \notag \\ 
& =\, u_{\Tha} \pn{b} 
\hspace{22mm} \in\, \GmS{\fld(b)}\bigpn{R_{\fld(b)}} \notag 
\end{align} 
where 
\[ \rmp{\sD_{\That}}_{\Base}: \vrt{\mfn^*\That} \dra \Gp\pn{\sD_{\That}}_{\Base} 
\] 
is the map obtained from the map 
$\rmp{\sD_{\That}}: \vrt{\That} \dra 
  \Gp\pn{\sD_{\Tha}}$ 
from \Point \ref{pairingFml}   
and the universal property of the fibre product 
\[ \xymatrix{ 
\vrt{\mfn^*\That} \ar[r] \ar@{-->}[ddr]^(0.36){\rmp{\sD_{\That}}_{\Base}} \ar[d] 
& \vrt{\That} \ar@{-->}[ddr]^{\rmp{\sD_{\That}}} \ar[d] & \\ 
\rC \ar[r]|!{[u];[dr]}\hole_(0.36){\mfn} \ar@/_0.0pc/[ddr] & X \ar[ddr]|!{[d];[dr]}\hole \\ 
& \Gp\pn{\sD_{\Tha}}_{\Base} \ar[r] \ar[d] \ar@{}[dr]|{\square} 
& \Gp\pn{\sD_{\Tha}} \ar[d] \\ 
& \Base \ar[r] & \Spec k 
}
\] 
and 
\[ u_{\Tha} := \rmp{\sD_{\That}}_{\Base} 
      \bigpn{\mfn^*\That \tms_{\Base} \Opn_{\Tha} }^{-1} 
\hspace{3mm} 
\in \GmS{\Base} \bigpn{R_{\,\Opn_{\Tha}}} 
\laurink 
\] 
for some open dense $\Opn_{\Tha} \subset \Base$ 
(here we use the notation from 
\Point \ref{trace_algGrps}). 

\medskip 

\textbf{Step 2:} Using Rigidity. \\ 
Let $\Pnt$ be either $\gnccl$ or a closed point $b \in \Base$. 
Let \,$\st{\tha_{n,\nu,\sig,\mu}^{\;\bgst}}_{n,\nu,\sig,\mu} \subset \sK_X^*$ 
be a family of functions on $X$, 
such that $\bigst{\tha_{n,\nu,\sig,\mu}^{\;\bgst}|_{\rC(\Pnt)}}_{n,\nu,\sig,\mu}$ 
is a pro-basis of $\ACHmp{\rC(\Pnt)}{\mdl_{\rCPntgst}}$. 
The existence of such a family of functions on $X$ is shown 
for $\Pnt = \gnccl$\, in Lemma \ref{local basis over B}, 
and for $\Pnt = b$ 
it follows from Proposition \ref{basis-exists} 
and Remark \ref{non-irreducible curves}. 
Set 
\[ \Tha_{n,\nu,\sig,\mu}^{\;\bgst} := \dv\bigpn{\tha_{n,\nu,\sig,\mu}^{\;\bgst}}_X 
\;\laurin 
\] 
For any curve \,$\iot{C}: C \ra X$\, denote 
\begin{align*} 
\alp_{n,\nu,\sig,\mu}^{C} 
& 
    := \bigbt{{\iot{C}}^*\Tha_{n,\nu,\sig,\mu}^{\;\bgst}} 
    \in \ACHmp{C}{\mdl_{\Cgst}} \\ 
u_{n,\nu,\sig,\mu}^{\;\bgst} 
& := \rmp{\sD_{\That_{n,\nu,\sig,\mu}^{\;\bgst}}}_{\Base}
   \bigpn{\mfn^*\That_{n,\nu,\sig,\mu}^{\;\bgst} \tms_{\Base} \Opn_{n,\nu,\sig,\mu}}^{-1} 
   \in \GmS{\Base} \pn{R_{\,\Opn_{n,\nu,\sig,\mu}}} 
\end{align*} 
for some open dense $\Opn_{n,\nu,\sig,\mu} \subset \Base$. 
Consider the system of equations 
\begin{align} 
\label{SYS} \mytag{SYS} 
\bigpair{\sD, \alp_{n,\nu,\sig,\mu}^{\rC(\Pnt)}}_{\rC(\Pnt),\mdl} 
\,=\;\, & u_{n,\nu,\sig,\mu}^{\;\bgst}\pn{\Pnt} \\ 
\fral \, n \in N\bigpn{\mdl_{\rCPntgst}},  
\;\; \fral 
\begin{array}{c} 
0 \leq \nu < n,  \\ 
\pn{\nu,p} = 1, 
\end{array} 
\;\; & \fral \, 1 \leq \sig \leq \sm(n),  
\;\; \fral \, \mu \geq 1 \notag 
\end{align} 
where $\sD \in \Fm{\rC(\Pnt)}{\mdl_{\rCPntgst}}\pn{R}$ 
is an indeterminate divisor, 
and the modified system 
\begin{align} 
& & & \fral \,n \in N\bigpn{\mdl_{\rCPntgst}}, \;\;\fral \,1 \leq \sig \leq \sm(n) \notag \\ 
\label{SY} \mytag{SY} 
\bigpair{\sD_{\et}, \alp_{n,0,\sig,\mu}^{\rC(\Pnt)}}_{\rC(\Pnt),\mdl} 
   &\,=\, t_{n,\sig,\mu}^{\;\bgst}(\Pnt) 
   & & \fral \,\mu \gg 1 \\ 
\bigpair{\sD_{\inf}, \alp_{n,\nu,\sig,\mu}^{\rC(\Pnt)}}_{\rC(\Pnt),\mdl} 
   &\,=\, u_{n,\nu,\sig,\mu}^{\;\bgst}(\Pnt) 
   & & \fral \,\mu \geq 1, 
   \;\;\fral \,0 < \nu < n, \,\pn{\nu,p} = 1 \notag 
\end{align} 
where 
\[ t_{n,\sig,\mu}^{\;\bgst} := \bigpn{\Lin_R \ra \Trz_R} \bigpn{u_{n,0,\sig,\mu}^{\;\bgst}} 
\laurink 
\] 
in particular 
$t_{n,\sig,\mu}^{\;\bgst}(\Pnt)$ is the image of $u_{n,0,\sig,\mu}^{\;\bgst}(\Pnt)$ 
under the functoriality map \\ 
$\Gm\bigpn{R \tens k(\Pnt)} \lra \Gm\bigpn{R_{\red} \tens k(\Pnt)}$. 

\bClm 
\label{SYS-implies-SY}
Any solution of (\ref{SYS}) is a solution of (\ref{SY}). 
\eClm 

\bBf 
Write $\resfld := k(\Pnt)$, 
$\pair{\llul,\lull} := \pair{\llul,\lull}_{\rC(\Pnt),\mdl}$ 
and 
$\fmlG := \Fm{\rC(\Pnt)}{\mdl_{\rCPntgst}}$ 
for the formal $\resfld$-group Cartier dual to 
$L := \Lmp{\rC(\Pnt)}{\mdl_{\rCPntgst}}$, 
$L(\resfld) = \ACHmp{\rC(\Pnt)}{\mdl_{\rCPntgst}}$. 

As \,$\sD_{\et} := \pn{\fmlG \ra \fmlG_{\et}}\pn{\sD}$ 
is by definition the image of \,$\sD \in \fmlG\pn{R_{\resfld}}$ 
under the functoriality map 
\,$\fmlG\pn{R \tens \resfld} \lra \fmlG\pn{R_{\red} \tens \resfld}$, 
it holds for all $\alp \in L(\resfld)$: 
\[ \bigpair{\sD_{\et}, \alp} 
\,=\, \bigpn{\Lin_R \ra \Trz_R} \bigpair{\sD, \alp} 
\laurin 
\] 
This shows that (\ref{SYS}) implies the \'etale part of (\ref{SY}). 

By construction $\tha_{n,\nu,\sig,\mu}^{\;\bgst}|_{\rC(\Pnt)}$ 
lies in the unipotent part $U$ of $L$ for $\nu > 0$. 
The image of any $\yps \in U(\resfld)$ 
under $\pair{\sD,\lull}: L(\resfld) \lra \Lin_R\pn{\resfld}$ 
lies in the unipotent part \,$\Upf_R = \ker\bigpn{\Lin_R \ra \Trz_R}$ of \,$\Lin_R$, 
so $\bigpair{\sD_{\et},\yps} = \bigpn{\Lin_R \ra \Trz_R} \bigpair{\sD,\yps}= 0$. 
Then 
\[ \bigpair{\sD, \yps} 
= \bigpair{\sD_{\et} + \sD_{\inf}, \yps} 
= \bigpair{\sD_{\inf}, \yps} 
\laurin 
\] 
This shows that (\ref{SYS}) implies also the infinitesimal part of (\ref{SY}). 
\eBf 

\medskip 

\noindent 
If $\chr(k) = p > 0$ 
let 
\begin{align*} 
\bigpn{\sD_{\rC(\Pnt),\nmb}}_{\nmb \geq 1} 
& \,:=\, \textrm{output of Algorithm \ref{algorithm} 
applied to (\ref{SY})} 
\laurin 
\end{align*} 
If $\chr(k) = 0$ 
let $\bigpn{\sD_{\rC(\Pnt),\nmb}}_{\nmb \geq 1}$ 
be the constant sequence 
that represents the unique solution of (\ref{SY}), 
which exists according to Lemma \ref{Basis-determines_0}. 

\bClm 
\label{Dmb stationary}
For every closed point \,$\bcp \in \Base$\, 
the sequence \,$\pn{\sD_{\rC(\bcp),\nmb}}_{\nmb \geq 1}$ 
becomes stationary and the limit is \,$\sD_{\rC(\bcp)}$: 
\[ 
\exists \,\mty(\bcp) \geq 1 \quad \forall \,\nmb \geq \mty(\bcp): 
\hspace{8mm} 
\sD_{\rC(\bcp),\nmb} = \sD_{\rC(\bcp)} 
\laurin 
\] 
\eClm 

\bBf 
If $\chr(k) = 0$ the sequence is stationary 
and coincides with the unique solution of (\ref{SY}) 
by definition. 
For $\chr(k)  > 0$, 
according to the Rigidity Lemma \ref{Basis-determines_p} 
the system (\ref{SY}) with \,$\Pnt = \bcp \in \Base$\, 
has as only candidate for a solution 
the limit of \,$\pn{\sD_{\rC(\bcp),\nmb}}_{\nmb \geq 1}$, 
if the limit exists. 
On the other hand, 
by construction (\ref{PU}) 
the divisor $\sD_{\rC(\bcp)}$ is a solution of (\ref{SYS}) for $\Pnt = \bcp$. 
By Claim \ref{SYS-implies-SY} then 
$\sD_{\rC(\bcp)}$ is a solution of (\ref{SY}) for $\Pnt = \bcp$ 
and so coincides with the limit of $\pn{\sD_{\rC(\bcp),\nmb}}_{\nmb \geq 1}$, 
in particular $\pn{\sD_{\rC(\bcp),\nmb}}_{\nmb \geq 1}$ becomes stationary 
($\see$Theorem \ref{System-Solve}), 
in either case. 
\eBf 


\bClm 
\label{BaseChange-for-Alg}
For $\Pnt = \gnccl$ and for every $\nmb \geq 1$ 
the relative Cartier divisor $\sD_{\rC(\gnccl),\nmb}$ 
descends to 
a relative Cartier divisor \,$\sD_{\rC(\gnc),\nmb}$ on $\rC(\gnc)$. 
For every \,$\nmb \geq 1$\, 
there exists an open dense subset 
$U_{\nmb}^{\bgst} \subset \Base$ 
such that $\sD_{\rC(\gnc),\nmb}$ yields a relative divisor on 
$\rC \tms_{\Base} U_{\nmb}^{\bgst}$, 
and for every $b \in U_{\nmb}^{\bgst}$ it holds 
\[ \sD_{\rC(b),\nmb} \,=\, \sD_{\rC(\gnc),\nmb} \cut \rC(b) 
\laurin 
\] 
(For the construction of the beginning part 
$\pn{\sD_{\rC(b),\mu}}_{1 \leq \mu \leq \nmb}$ 
only the part \\ 
$\st{\tha_{n,\nu,\sig,\mu}^{\;\bgst} \;|\; n,\nu,\sig,\mu \leq \nmb}$ 
of the pro-basis is involved.) 
\eClm 

\bBf 
While the functions $\tha_{n,\nu,\sig,\mu}^{\;\bgst}$ 
and the relative Cartier divisors $\sD_{\Tha_{n,\nu,\sig,\mu}^{\;\bgst}}$ 
are defined over $X$, 
the objects $\alp_{n,\nu,\sig,\mu}^{\rC(\gnc)}$, 
$u_{n,\nu,\sig,\mu}^{\;\bgst}(\gnc)$, $t_{n,\sig,\mu}^{\;\bgst}(\gnc)$ 
are defined over $\fld(\gnc)$, 
and $\alp_{n,\nu,\sig,\mu}^{\rC(\gnccl)}$, 
$u_{n,\nu,\sig,\mu}^{\;\bgst}(\gnccl)$, $t_{n,\sig,\mu}^{\;\bgst}(\gnccl)$ 
are their pull-backs to $\gnccl$. 

For the infinitesimal part of $\sD_{\rC(\gnccl),\nmb}$ 
one sees directly from procedure (\ref{Alg}) 
that the output of Algorithm \ref{algorithm} 
is defined over $k(\gnc)$. 
More precisely, 
$\sD_{\rC(\gnccl),\nmb}^{\,\inf}$ is defined over 
the open dense subset $U^{\bgst}_{\nmb} \subset \Base$ 
of all those $b \in \Base$ such that 
$\st{\tha_{n,\nu,\sig,\mu}^{\;\bgst}|_{\rC(b)}\;;\; n,\nu,\sig,\mu \leq \nmb}$ 
forms the beginning part of a pro-basis of $\ACHmp{\rC(b)}{\mdl}$ 
(cf.\ Lemma \ref{local basis over B}), 
and is compatible with specialization on $U^{\bgst}_{\nmb} \subset \Base$. 

For the \'etale part of $\sD_{\rC(\gnccl),\nmb}$ 
we show that $\sD_{\rC(\gnccl),\nmb}^{\,\et}$ 
is defined over $k(\gnc)$, 
hence over some open dense of $\Base$, 
and its multiplicities coincide with those of $\sD_{\rC(b),\nmb}^{\,\et}$ 
for general \,$b \in \Base$. 
The family of \'etale relative divisors 
$\bigst{\sD_{\rC(b)}^{\,\et}}_{b \in \Base}$ 
has constant multiplicities 
\;$\pn{\nbr_{\pnt_{b}}}_{\pnt_{b}} 
\in \dsum_{\pnt_{b} \in \vrt{\mdl_{\rC(b)}}} \Zint$\; 
among general \,$b \in \Base$: \\  
we have 
\,$
\bigpair{\sD_{\rC(b)}^{\,\et}, \alp^{\rC(b)}_{n,0,\sig,\nmb}}_{\rC(b),\mdl_{\rC(b)}} 
= \prod_{\pntt_{b} \in \Sm_{b}(n)} \tha_{n,0,\sig,\nmb}\pn{\pntt_{b}}^{\nbr_{\pntt_{b}}} 
= t_{n,\sig,\nmb}(b) 
\laurink 
$\, 
and $\tha_{n,0,\sig,\nmb}\pn{\pntt_{b}}$ and $t_{n,\sig,\nmb}(b)$ 
are both continuous functions in $b \in \Base$. 
Moreover, let $\rho: \EBase \ra \Base$ be a generically finite covering 
with $\EBase$ irreducible, 
and let 
$\rho_{\rC}: \rC_{\EBase} \ra \rC$ 
be the corresponding base changed map. 
Then the system 
\begin{align} 
\label{Etmp} \mytag{Et(\nmb,\Pnt)} 
& \prod_{\pntt_{\Pnt} \in \Sm_{\Pnt}(n)} 
\rho_{\rC}^*\,\tha_{n,0,\sig,\nmb}\pn{\pntt_{\Pnt}}^{\nbr_{n,\pntt_{\Pnt}}} 
= \rho_{\rC}^*\,t_{n,\sig,\nmb}(\Pnt) \\ 
& \hspace{30mm} 
\sig = 1, \ldots, \sm_{\Pnt}(n) \notag 
\end{align} 
($\see$Notation \ref{Modl-invariants}) 
has for $\Pnt = \ebase \in \EBase$ over $b \in \Base$ 
the unique solution 
\bMyEqn \label{LH} \mytag{LH} 
\pn{\nbr_{n,\pntt_{\ebase}}}_{\pntt_{\ebase} \in \Sm_{\ebase}(n)} 
= \pn{\nbr_{n,\pntt_{b}}}_{\pntt_{b} \in \Sm_{b}(n)} 
\hspace{15mm} \forall \,\ebase \ra b 
\eMyEqn 
in \;$\st{-p^{\nmb},\ldots,p^{\nmb}}^{\Sm_{\ebase}(n)}$\; 
if $\nmb$ is sufficiently large. 
Thus the solutions $\pn{\nbr_{\pnt_{\ebase}}}_{\pnt_{\ebase}}$ 
of (\ref{Etmp}) for $\Pnt = \ebase$ 
coincide with the solution for $\rho(h) \in \Base$, 
hence are constant among general $\ebase \in \EBase$. 
Let $\egnc$ be the generic point of $\EBase$. 
If the function field $k(\egnc)$ of $\EBase$ 
contains the residue fields $k(\pnt_{\gnc})$ 
for all $\pnt_{\gnc} \in \vrt{\mdl_{\rC(\gnc)}}$, 
there is a 1-1 correspondence 
$\vrt{\mdl_{\rC(\egnc)}} \llra \vrt{\mdl_{\rC(\ebase)}}$ 
between the irreducible components of $\mdl_{\rC_{\EBase}}$ 
and those of $\mdl_{\rC(\ebase)}$ 
for general $\ebase \in \EBase$. 
Then the system (\ref{Etmp}) for $\Pnt = \egnc$ 
has a unique solution 
in $\st{-p^{\nmb},\ldots,p^{\nmb}}^{\Sm_{\egnc}(n)}$ 
(for $\nmb$ sufficiently large) 
that corresponds to the one for general $b \in \Base$. 
Since $k(\gnccl)$ is the union of those $k(\egnc)$, 
the same holds for $\Pnt = \gnccl$. 
If $\nmb$ is not sufficiently large, 
we have $\sD_{\rC(b),\nmb}^{\,\et} = 0$ 
and $\sD_{\rC(\gnccl),\nmb}^{\,\et} = 0$, 
so in this case the assertion is fulfilled trivially. 
Due to (\ref{LH}), 
the solution 
$\pn{\nbr_{n,\pntt_{\gnccl}}}_{\pntt_{\gnccl}}$ 
of (\ref{Etmp}) for $\Pnt = \gnccl$\, satisfies 
\bMyEqn \label{LE} \mytag{LE} 
\hspace{0mm} \nbr_{n,\pntt_{\gnccl}} = \nbr_{n,\pntt'_{\gnccl}} 
\hspace{15mm} \forall \,\pntt_{\gnccl},\pntt'_{\gnccl} \ra \pntt_{\gnc} 
\laurink 
\eMyEqn 
and therefore the solution 
$\pn{\nbr_{n,\pntt_{b}}}_{\pntt_{b}}$ 
of (\ref{Etmp}) for $\Pnt = b \in \Base$ in general position satisfies 
\bMyEqn \label{LB} \mytag{LB} 
\hspace{0mm} \nbr_{n,\pntt_{b}} = \nbr_{n,\pntt'_{b}} 
\hspace{15mm} \forall \,\pntt_{b},\pntt'_{b} \ra \ol{\st{\pntt_{\gnc}}} 
\laurin 
\eMyEqn 
In other words, 
the \'etale relative divisor $\sD_{\rC(\gnccl),\nmb}^{\,\et}$ 
corresponding to 
$\pn{\nbr_{\pnt_{\gnccl}}}_{\pnt_{\gnccl}} 
\in \dsum_{\pnt_{\gnccl} \in \vrt{\mdl_{\rC(\gnccl)}}} \Zint$\, 
with respect to (\ref{FZ}) 
has the same multiplicity at all points of $\rC(\gnccl)$ 
lying over the same irreducible component of $\vrt{\mdl}$, 
which coincides with 
the multiplicity of the \'etale divisor $\sD_{\rC(b),\nmb}^{\,\et}$ 
corresponding to 
$\pn{\nbr_{\pnt_{b}}}_{\pnt_{b}} 
\in \dsum_{\pnt_{b} \in \vrt{\mdl_{\rC(b)}}} \Zint$\, 
at all points of $\rC(b)$ 
lying in this component of $\vrt{\mdl}$. 
This shows that $\sD_{\rC(\gnccl),\nmb}^{\,\et}$ is defined over $k(\gnc)$ 
and is compatible with specialization. 
\eBf 

\medskip 

The sequence $\pn{\sD_{\rC(\gnc),\nmb}}_{\nmb \geq 1}$ 
is constant by definition for $\chr(k) = 0$, 
and it becomes stationary 
according to Lemma \ref{pro-divisor terminates} 
for $\chr(k) = p > 0$. 
Hence it yields a relative Cartier divisor 
\,$\sD_{\rC(\gnc)} 
\in \Divf_{\rC(\gnc)}^{\vrt{\mdl \cut \rC(\gnc)}} \pn{R_{k(\gnc)}}$, 
defined over 
$U := U^{\bgst}_{\mty\pn{\gnc}}$, 
where \,$\mty\pn{\gnc}$\, is the smallest \,$\nmb \geq 1$\, 
such that \,$\sD_{\rC(\gnc),\nmb} = \sD_{\rC(\gnc)}$. 
Thus we obtain a relative Cartier divisor 
\,$\sD_{\rC \tms_{\Base} U} \in 
\Divf_{\rC \tms_{\Base} U}^{\vrt{\mdl \cut \rC} \tms_{\Base} U} \pn{R_U} 
$. 

For every $b \in U$ 
the beginning part 
$\pn{\sD_{\rC(b),\nmb}}_{1 \leq \nmb \leq \mty(\gnc)}$ 
is obtained from 
$\pn{\sD_{\rC(\gnc),\nmb}}_{1 \leq \nmb \leq \mty(\gnc)}$ 
by base change ($\see$Claim \ref{BaseChange-for-Alg}). 
By Claim \ref{Dmb stationary} 
the limit of $\pn{\sD_{\rC(b),\nmb}}_{\nmb \geq 1}$ exists 
and is $\sD_{\rC(b)}$. 
By Corollary \ref{finite subsystem} 
and since \,$\lgth\bigpn{\sD_{\rC(b)}} \leq \lgth\bigpn{\sD_{\rC(\gnc)}}$ 
we have \,$\sD_{\rC(b),\mty(\gnc)} = \sD_{\rC(b)}$. 
Inserting the limits into Claim \ref{BaseChange-for-Alg} we obtain 
\bMyEqn \label{SPC} \mytag{SPC} 
\sD_{\rC \tms_{\Base} U} \cut \rC(b) = \sD_{\rC(b)} 
\hspace{20mm} \forall \,b \in U 
\laurin 
\eMyEqn 

\medskip 

\textbf{Step 3:} From Local to Global. \\ 
For arbitrary $\bspp \in \Base$ 
and any $\alp \in \ACHm{\rC(\bspp)}{\mdl_{\rCbsppgst}}$ 
let $\tha \in \sK_X^*$ be a function on $X$ with $\Tha := \dv\pn{\tha}_X$ 
such that $\Tha \cut \rC(\bspp) = \alp$. 
Let $\Opn_{\Tha}$ be the set of all those $b \in \Base$ 
such that $\rC(b)$ meets $\Tha$ properly. 
Then due to (\ref{SPC}) and (\ref{PU}) it follows for all 
$b \in U^{\bgst} \cap \Opn_{\Tha}$ 
\[ \bigpair{\sD_{\rC \tms_{\Base} U} \cut \rC(b), \Tha \cut \rC(b)}_{\rC(b),\mdl} 
=\, \bigpair{\sD_{\rC(b)}, \Tha \cut \rC(b)}_{\rC(b),\mdl} 
=\, u_{\Tha} \pn{b} 
\laurin 
\] 
Since the algebraic function \,$u_{\Tha}$\, is regular at $b = \bspp$, 
the domain of definition of the expression 
\,$\bigpair{\sD_{\rC \tms_{\Base} U} \cut \rC(b), \Tha \cut \rC(b)}_{\rC(b),\mdl}$\, 
extends to the point \,$b = \bspp$, 
and then it holds 
\[ \bigpair{\sD_{\rC \tms_{\Base} U} \cut \rC(\bspp), \alp}_{\rC(\bspp),\mdl} 
= \bigpair{\sD_{\rC(\bspp)}, \alp}_{\rC(\bspp),\mdl} 
\laurin 
\] 
Therefore  
\[ \bigpair{\sD_{\rC \tms_{\Base} U} \cut \rC(\bspp), \lull}_{\rC(\bspp),\mdl} 
= \bigpair{\sD_{\rC(\bspp)}, \lull}_{\rC(\bspp),\mdl} 
: \ACHm{\rC(\bspp)}{\mdl_{\rCbsppgst}} \lra \Gm\pn{R} 
\] 
gives a well-defined element of 
\;$\Homabk\pn{\Lm{\rC(\bspp)}{\mdl_{\rCbsppgst}},\Lin_R}
= \Fmo{\rC(\bspp)}{\mdl_{\rCbsppgst}} \pn{R}$. 
This shows that $\sD_{\rC \tms_{\Base} U}$ specializes also to 
a divisor on $\rC(\bspp)$, 
and 
\,$\sD_{\rC \tms_{\Base} U} \cut \rC(\bspp) = \sD_{\rC(\bspp)}$. 
Thus $\sD_{\rC \tms_{\Base} U}$ extends to a relative divisor $\sD_{\rC}$ 
defined on all of $\rC$ satisfying 
\[ \sD_{\rC} \cut \rC(b) = \sD_{\rC(b)} 
\hspace{20mm} \forall \,b \in \Base 
\laurin 
\] 
This concludes the proof of Lemma \ref{locRelDiv}. 
\ePf 

\bNot  
\label{trace_algGrps}
Let $\rmpC: \rC \dra \sG$ be a rational map of schemes over $\Base$, 
where $\sG$ is a $\Base$-group scheme. 
Let $\fz$ be a relative 0-cycle on $\rC$ over $\Base$, 
i.e.\ a finite formal sum $\fz = \sum n_{\fx} \,\fx$, $n_{\fx} \in \Zint$, 
$n_{\fx} = 0$ for almost all $\fx$, 
where $\fx$ ranges over all integral subschemes $\fx$ of $\rC$ 
that are finite and flat over $\Base$, 
and suppose $\rmpC$ is defined on all $\fx$ with $n_{\fx} \neq 0$. 
Then we denote 
\[ \rmpC\pn{\fz} 
:= \sum n_{\fx} \Trace_{\fx/\Base}\bigpn{\rmpC\pn{\fx}} 
\in \sG\pn{\Base} 
\] 
where $\Trace_{\fx/\Base}: \sG\pn{\fx} \lra \sG\pn{\Base}$ 
is the trace map induced by the finite flat morphism $\fx \ra \Base$, 
see \cite[XVII, (6.3.13.2)]{SGA4}. 
\eNot  

\bNot 
\label{Modl-invariants}
For any morphism $\Pnt \ra \Base$ 
we define $\Sm_{\Pnt}$, $\Sm_{\Pnt}(n)$, $\resfld_{\Pnt}(n)$, $\sm_{\Pnt}(n)$ 
analogous to $\Sm$, $\Sm(n)$, $\resfld(n)$, $\sm(n)$ 
from Definition \ref{n-part}, 
with $C$ replaced by $\rC \tms_{\Base} \Pnt$ 
and $\mdl$ replaced by $\mdl_{\rC \tms_{\Base} \Pnt}$. 
\eNot 


\bLem 
\label{local basis over B}
For a general \,$\bsp \in \Base$\, 
there are functions \,$\tha_{n,\nu,\sig,\mu}^{\;\bgst} \in \sK_X^*$ on $X$ 
for $n \in N(\mdl)$, $0 \leq \nu < n$, $1 \leq \sig \leq \sm(n)$, $\mu \geq 1$, 
such that $\st{\mfn^*\tha_{n,\nu,\sig,\mu}^{\;\bgst}}_{n,\nu,\sig,\mu}$ 
is a pro-basis of $\ACHmp{\rC(\gnccl)}{\mdl_{\rCgncclgst}}$, 
and $\st{\tha_{n,\nu,\sig,\mu}^{\;\bgst}|_{\rC(\bsp)}}_{n,\nu,\sig,\mu}$ 
is a pro-basis of $\ACHmp{\rC(\bsp)}{\mdl_{\rCbspgst}}$. 
For any integer $\nmb \geq 1$, the locus 
\[ U^{\bgst}_{\nmb}
:= \lrst{b \in \Base \left| 
\begin{array}{l} 
\textrm{the set } 
\st{\tha_{n,\nu,\sig,\mu}^{\;\bgst}|_{\rC(b)}\;;\; n,\nu,\sig,\mu \leq \nmb} 
\textrm{ can be} \\ 
\textrm{extended to a pro-basis of } 
\ACHmp{\rC(b)}{\mdl_{\rCbgst}} 
\end{array} 
\right.} 
\] 
is open and dense in $\Base$. 
\eLem 

\bPf 
Write $\mfn^*\mdl = \mdl_{\rC} = \sum_{n \in N(\mdl_{\rC})} n \,\Es_{n}$, 
where the $\Es_{n}$ are reduced effective divisors on $\rC$, finite over $\Base$, 
with ideal sheaf $\sI_{\rC,\Es_{n}}$. 
Then we have $\Es_{n}(\Pnt) = \Sm_{\Pnt}(n)$ for any $\Pnt \ra \Base$. 
Let $\Opnn \subset \Base$ be the open dense subset  
where $\strl|_{\mfn^*\mdl}: \mfn^*\mdl \lra \Base$ is \'etale. 
For the algebraic $\Base$-group 
\[ \Lmp{\rC}{\mdl_{\rC}} 
:= \WRS{\GmS{\Base}}{\pi_*\sO_{\mdl_{\rC}}}{\sO_{\Base}} 
\laurink 
\] 
we have the following relation between the fibres over 
$\gnccl \ra \Base$ and $b \in \Opnn$: 
\[ \Lmp{\rC(\gnccl)}{\mdl_{\rCgncclgst}} 
\cong \Lmp{\rC(b)}{\mdl_{\rCbgst}} \tens_{\fld} \fld(\gnccl) 
\laurin 
\] 
For any point $\Pnt \ra \Base$ with $\rC(\Pnt)$ smooth 
we have (cf.\ Definition \ref{ACHm'}): 
\begin{align*} 
\ACHmp{\rC(\Pnt)}{\mdl_{\rCPntgst}} 
& = \Lmp{\rC(\Pnt)}{\mdl_{\rCPntgst}}\bigpn{\fld(\Pnt)} \\ 
& = \strl_*\sO_{\mdl_{\rCPnt}}^* \\ 
& = \prod_{n \in N(\mdl_{\rC})} \frac{\strl_*\sO_{\rC(\Pnt),\Es_{n}(\Pnt)}^*}
{1+ \strl_*\sI_{\rC(\Pnt),\Es_{n}(\Pnt)}^{\,n}} 
\;\laurin 
\end{align*} 
Let \,$\pjn_{n}: \pi_*\sO_{\mdl_{\rCPnt}}^* \lsur 
\frac{\strl_*\sO_{\rC(\Pnt),\Es_{n}(\Pnt)}^*}
{1+ \strl_*\sI_{\rC(\Pnt),\Es_{n}(\Pnt)}^{n}}$\, 
be the projection to the $n$-component, 
$\resfld_{\Pnt}(n) = \prod_{\pntt_{\Pnt} \in \Es_{n}\pn{\Pnt}} \fld\pn{\pntt_{\Pnt}}$, 
$\Lam_{n,\resfld_{\Pnt}(n)} 
= \prod_{\pntt_{\Pnt} \in \Es_{n}\pn{\Pnt}} \Lam_{\pntt_{\Pnt},n}$, 
$\Lam_{\pntt_{\Pnt},n} 
= \frac{1 + \fm_{\rC(\Pnt),\pntt_{\Pnt}}}{1 + \fm_{\rC(\Pnt),\pntt_{\Pnt}}^{\,n}}$. 
Let $g_{n} \in \sK_X$ such that 
$\mfn^*g_{n} \in \sI_{\rC,\Es_{n},\bsp} \setminus \sI_{\rC,\Es_{n},\bsp}^2$, 
where 
$\sI_{\rC,\Es_{n},\bsp} = \sI_{\rC,\Es_{n}} \tens_{\sO_{\Base}} \sO_{\Base,\bsp}$ 
for the given $\bsp \in \Opnn$, 
and with $\mfn^*g_{n}(\bsp) 
\in \fm_{\rC(\bsp),\pntt_{\bsp}} \setminus \fm_{\rC(\bsp),\pntt_{\bsp}}^2$ 
for all $\pntt_{\bsp} \in \Es_{n}(\bsp)$. 
If $\lam_{n,\nu} = \Expah\pn{\mfn^*g_{n}^{\,\nu}}$, then 
$\bigst{\lam_{n,\nu}(\Pnt)}_{\substack{0<\nu<n \\ (\nu,p)=1}}$ 
is a basis of $\Lam_{n,\resfld_{\Pnt}(n)}$ 
for $\Pnt = \gnccl$ and $\Pnt = \bsp \in \Opnn$. 
Since $\Es_{n} = \mfn^*\Es_{n,X}$ 
for some reduced effective divisors $\Es_{n,X}$ on $X$, 
and the ideal sheaves $\sI_{X,\Es_{n,X}}$ are relatively prime on a dense open 
subset of $X$, 
we find $\tha_{n,\nu} \in \sK_X^*$ 
with $\pjn_{n}\bigpn{\mfn^*\tha_{n,\nu}(\gnccl)} = \lam_{n,\nu}(\gnccl)$ 
and $\pjn_{n'}\bigpn{\mfn^*\tha_{n,\nu}(\gnccl)} = 1$ for all $n' \neq n$, 
by the Chinese Remainder Theorem. 

For $\nu \geq 1$ we construct $\tha_{n,\nu,\sig,\mu} \in \sK_X^*$ 
in the same way as $\tha_{n,\nu}$, 
but with $g_{n}^{\,\nu}$ replaced by $g_{n}^{\,\nu} \,c_{n,\nu,\sig} \,c_\mu$, 
where $c_{n,\nu,\sig} \in \sK_X$ satisfies 
$\mfn^*c_{n,\nu,\sig}\pn{\pntt_{\bsp,\sig'}} = \del_{\sig\sig'}$ 
for $\pntt_{\bsp,\sig'} \in 
\Es_{n}(\bsp) = \st{\pntt_{\bsp,\sig} \;|\; 1 \leq \sig \leq \sm_{\bsp}(n)}$, 
$\sm_{\bsp}(n) = \dim_{k}\resfld_{\bsp}(n) = \# \Es_{n}(\bsp)$, 
and $c_{\mu}$ a generator of the group of units of $\bF_{q^{\mu+1}}$. 

For $\nu = 0$ we construct $\tha_{n,0,\sig,\mu} \in \sK_X^*$ by the condition 
\[ \Bigpn{\pjn_{n'}\bigpn{\mfn^*\tha_{n,0,\sig,\mu}(\bsp)}}\pn{\pntt_{\bsp,\sig'}} = 
\left\{ 
\begin{array}{ll} 
c_{\mu} & \textrm{ if } n' = n, \; \sig' = \sig \\ 
1 & \textrm{ otherwise.} 
\end{array} 
\right. 
\] 
One checks that the family $\st{\tha_{n,\nu,\sig,\mu}}_{n,\nu,\sig,\mu}$ 
yields bases of $\ACHmp{\rC(\gnccl)}{\mdl_{\rC(\gnccl)}}$ 
and of $\ACHmp{\rC(\bsp)}{\mdl}$. 

For any $\nmb \geq 1$ 
let $\Opnn_{\nmb}^{\bgst} \subset \Base$ be the open dense subset 
of those $b \in \Base$ such that $\rC(b)$ meets 
$\Tha_{n,\nu,\sig,\mu}^{\;\bgst} = \dv\pn{\tha_{n,\nu,\sig,\mu}^{\;\bgst}}$ 
properly for all $\,n,\nu,\sig$ and $\mu \leq \nmb$: 
\[ \Opnn_{\nmb}^{\bgst} 
:= \bigst{b \in \Base \;\big|\;\; \rC(b) \iscp \Tha_{n,\nu,\sig,\mu}^{\;\bgst} 
\hspace{4mm} \forall \,n,\nu,\sig, \forall \,\mu \leq \nmb}
\laurin 
\] 
For $b \in \Opnn_{\nmb}^{\bgst}$ 
the set $\st{\tha_{n,\nu,\sig,\mu}^{\;\bgst}|_{\rC(b)}\;;\; n,\nu,\sig,\mu \leq \nmb}$ 
satisfies condition (PB1) from Definition \ref{pro-basis-ACHm} 
by construction; 
thus it can be extended to a pro-basis of $\ACHmp{\rC(b)}{\mdl_{\rCbgst}}$ 
if and only if 
\begin{align*}
\textrm{(PB2)} \hspace{3mm} & 
  \textrm{The system } \hspace{-2mm} 
  \prod_{\pntt_{b} \in \Es_{n}(b)} \tha_{n,0,\sig,\mu}\pn{\pntt_{b}}^{\nbr_{\pntt_{b}}} 
  = 1, 
  \hspace{2mm} \sig = 1, \ldots, s_{b}(n) \hspace{2.6mm} \textrm{has only the } \\ 
  & \textrm{trivial solution } 
  \,\pn{\nbr_{\pntt_{b}}}_{\pntt_{b}} = 0 \, 
  \textrm{ in } \st{-p^{\mu},\ldots,p^{\mu}}^{\Es_{n}(b)} 
  \hspace{5mm} 
  \forall \, n,\;\forall \,1 \leq \mu \leq \nmb \\ 
\textrm{(PB3)} \hspace{3mm} &  
  g_{n}|_{\rC(b)} \in \fm_{\rC(b),\pntt_{b}} \setminus \fm_{\rC(b),\pntt_{b}}^{2} 
  \hspace{43mm} \forall\,n, \;\forall \,\pntt_{b} \in \Es_{n}(b) \\ 
  & \bigpn{c_{n,\nu,\sig}\pn{\pntt_{b}}}_
  {\substack{\pntt_{b} \in \Es_{n}(b) \\ 1 \leq \sig \leq \sm_{b}(n)} } 
  \in \Gl_{\sm_{b}(n)}\pn{k} 
  \hspace{30.5mm} \forall\,n, \;\forall \,0 < \nu < n 
  \laurin 
\end{align*} 
Considering 
$\bigpn{\tha_{n,0,\sig,\mu}\pn{\pntt_{b,\sig'}}}_{\sig,\sig'} 
\in \Mat_{\sm_{b}(n) \tms \sm_{b}(n)}\pn{k} \cong k^{\sm_{b}(n)^2}$ 
as a point in  
\,$\Spec k\bt{X_{\sig \sig'} \;|\; 1 \leq \sig,\sig' \leq \sm_{b}(n)}$, 
conditions (PB2) and (PB3) are equivalent to 
\begin{align*}
\textrm{(PB2')} \hspace{3mm} & 
  \bigpn{\tha_{n,0,\sig,\mu}\pn{\pntt_{b,\sig'}}}_{\sig,\sig'} 
  \notin \bigcup_{\pn{\nbr_{\sig'}}_{\sig'} \neq 0} 
  \Zero\biggpn{\prod_{\sig'} X_{\sig \sig'}^{\,\nbr_{\sig'}} - 1 \;\bigg|\; 
  1 \leq \sig \leq \sm_{b}(n)} \\ 
  & 
  \textrm{where }\; 1 \leq \sig,\sig' \leq \sm_{b}(n), \quad 
  \pn{\nbr_{\sig'}}_{\sig'} \in \st{-p^{\mu},\ldots,p^{\mu}}^{\sm_{b}(n)} \\ 
  & \textrm{and }\; \Es_{n}(b) = \st{\pntt_{b,\sig'} \;|\; 1 \leq \sig' \leq \sm_{b}(n)}
  \hspace{18mm}  \forall \, n,\;\forall \,1 \leq \mu \leq \nmb \\ 
\textrm{(PB3')} \hspace{3mm} &  
  g_{n}|_{\rC(b)} \notin \fm_{\rC(b),\pntt_{b}}^{2} 
  \hspace{54mm} \forall\,n, \;\forall \,\pntt_{b} \in \Es_{n}(b) \\ 
  & \det \bigpn{c_{n,\nu,\sig}\pn{\pntt_{b}}}_
  {\substack{\pntt_{b} \in \Es_{n}(b) \\ 1 \leq \sig \leq \sm_{b}(n)} } 
  \neq 0 
  \hspace{32mm} \forall\,n, \;\forall \,0 < \nu < n 
  \laurink 
\end{align*} 
since \,$g_{n}|_{\rC(b)} \in \fm_{\rC(b),\pntt_{b}}$\, by construction 
of $\st{\tha_{n,\nu,\sig,\mu}}_{n,\nu,\sig,\mu}$. 
Then 
\[ U_{\nmb}^{\bgst} 
= \bigst{b \in \Opnn_{\nmb}^{\bgst} \;\big|\;\; 
    b \,\textrm{ satisfies (PB2') and (PB3')}} 
\laurin 
\] 
As conditions (PB2') and (PB3') are obviously finitely many open conditions 
and $U_{\nmb}^{\bgst} \neq \varnothing$ 
(since $\bsp \in U_{\nmb}^{\bgst}$), 
the set $U_{\nmb}^{\bgst}$ is open and dense in $\Base$. 
\ePf 


\bLem 
\label{pro-divisor terminates}
Let $\pn{\sD_{\rC(\gnccl),\nmb}}_{\nmb \geq 1}$ 
be the output of Algorithm \ref{algorithm} applied to {\rm (\ref{SY})} 
for $\Pnt = \gnccl$. 
This sequence becomes stationary. 
\eLem 

\bPf 
\textbf{Step 1:} \'Etale part of $\pn{\sD_{\rC(\gnccl),\nmb}}_{\nmb \geq 1}$. \\ 
By Claim \ref{BaseChange-for-Alg} 
we have 
\,$\pn{\sD_{\rC(\gnccl),\nmb}}_{\nmb \geq 1} 
= \pn{\sD_{\rC(\gnc),\nmb}}_{\nmb \geq 1}$\, 
and 
\[ \sD_{\rC(\gnc),\nmb} \cut \rC(b) \,=\, \sD_{\rC(b),\nmb} 
\hspace{17mm} \forall \,b \in U_{\nmb} 
\laurin 
\] 
If \,$\sD^{\,\et}_{\rC(\gnc),\nmb}$\, corresponds to the tuple 
\,$\pn{{}^{\nmb}\nbr_{\pntt_{\gnc}}}_{\pntt_{\gnc} \in \vrt{\mdl_{\rC(\gnc)}}} 
\in \dsum_{\pntt_{\gnc} \in \vrt{\mdl_{\rC(\gnc)}}} \Zint$, 
and \,$\sD^{\,\et}_{\rC(b),\nmb}$\, corresponds to the tuple 
\,$\pn{{}^{\nmb}\nbr_{\pntt_{b}}}_{\pntt_{b} \in \vrt{\mdl_{\rC(b)}}} 
\in \dsum_{\pntt_{b} \in \vrt{\mdl_{\rC(b)}}} \Zint$, 
then the proof of Claim \ref{BaseChange-for-Alg} showed 
\[ 
\pn{{}^{\nmb}\nbr_{\pntt_{\gnc}}}_{\pntt_{\gnc} \in \vrt{\mdl_{\rC(\gnc)}}} 
\,=\, \pn{{}^{\nmb}\nbr_{\pntt_{b}\pn{\pntt_{\gnc}}}}_
{\pntt_{\gnc} \in \vrt{\mdl_{\rC(\gnc)}}} 
\] 
where $\pntt_{b}\pn{\pntt_{\gnc}}$ is any point 
$\pntt_{b} \in \vrt{\mdl_{\rC(b)}}$ 
such that $\pntt_{b} \in \ol{\st{\pntt_{\gnc}}}$, 
and ${}^{\nmb}\nbr_{\pntt_{b}} = {}^{\nmb}\nbr_{\pntt_{b}'}$ 
for all $\pntt_{b}, \pntt_{b}' \in \ol{\st{\pntt_{\gnc}}}$. 
Moreover by Claim \ref{Dmb stationary} we have 
\[ \pn{{}^{\nmb}\nbr_{\pntt_{b}}}_{\pntt_{b} \in \vrt{\mdl_{\rC(b)}}} 
\,=\, \pn{\nbr_{\pntt_{b}}}_{\pntt_{b} \in \vrt{\mdl_{\rC(b)}}} 
\hspace{17mm} \forall \,\nmb \geq \mty(b) 
\] 
where $\pn{\nbr_{\pntt_{b}}}_{\pntt_{b} \in \vrt{\mdl_{\rC(b)}}} 
\in \dsum_{\pntt_{b} \in \vrt{\mdl_{\rC(b)}}} \Zint$\, 
is the tuple corresponding to $\sD_{\rC(b)}$. 
Thus the sequence $\bigpn{\sD^{\,\et}_{\rC(\gnc),\nmb}}_{\nmb}$ 
is stationary for \,$\nmb \geq \mty(\bsp)$, 
where $\bsp \in \Base$ is a general point as in Lemma \ref{local basis over B}. 

\medskip 

\textbf{Step 2:} Infinitesimal part of $\pn{\sD_{\rC(\gnccl),\nmb}}_{\nmb \geq 1}$. \\ 
For every $\nmb \geq 1$ 
the relative divisor $\sD^{\,\inf}_{\rC(\gnccl),\nmb}$ 
corresponds to a tuple 
\[ 
\bigpn{{}^{\nmb}_{\;1} w^{n,\nu}_{\,\pntt}}_{n,\nu,\pntt} 
= \bigpn{\pn{{}^{\nmb}_{\;1} w^{n,\nu}_{i,\pntt}}_{i \geq 0}}_{n,\nu,\pntt} 
\in 
\dsum_{n \in N\pn{\mdl_{\rCgncclgst}}} 
\dsum_{\substack{0 < \nu < n \\ (\nu,p) = 1}} 
\dsum_{\pntt \in \Sm_{\gnccl}(n)} \Witt\bigpn{R_{k(\gnccl)}} 
\laurin 
\] 
which is the $(\nmr = 1)$-part of a solution 
$\bigpn{{}^{\nmb}_{\;\nmr} w^{n,\nu}_{\,\pntt}}_{\nmr,n,\nu,\pntt}$ 
of the following system of 
matrix equations 
\begin{align} 
\label{matrix} \mytag{MAT} 
& \Expah\biggpn{ \Bigpn{\bigbt{c_{n,\nu,\sig}^{\,\pntt}}}_{\sig,\pntt} \cdot 
\Bigpn{{}^{\nmb}_{\;\nmr} w^{n,\nu}_{\,\pntt}}_{\pntt,\nmr} \cdot 
\Bigpn{\bigbt{\pn{c_{\mu}}^{\nmr}}}_{\nmr,\mu} } 
\,=\, \Bigpn{u_{n,\nu,\sig,\mu}\pn{\gnccl}}_{\sig,\mu} \\ 
& \hspace{0mm} n \in N\pn{\mdl_{\rCgncclgst}}, 
\hspace{5mm} 0 < \nu < n, \,(\nu,p) = 1 \notag 
\end{align} 

\noindent 
where 
\[ 
u_{n,\nu,\sig,\mu} 
= \rmp{\sD_{\That_{n,\nu,\sig,\mu}}}_{\Base}
   \bigpn{\vrt{\mfn^*\That_{n,\nu,\sig,\mu}} \tms_{\Base} \Opn_{n,\nu,\sig,\mu}}^{-1} 
\in \Upf_{R}\bigpn{\sO_{\Base}\pn{\Opn_{n,\nu,\sig,\mu}}} 
\laurin 
\] 
Since the \,$u_{n,\nu,\sig,\mu} \in R \tens \sO_{\Base}(\Opn_{n,\nu,\sig,\mu})$\, 
are rational functions on $\Base$, 
we have \,$u_{n,\nu,\sig,\mu}\pn{\gnccl} \in k(\gnc)$. 
Then procedure (\ref{Alg}) from Algorithm \ref{algorithm} 
shows that also 
\[ {}^{\nmb}_{\;\nmr} w^{n,\nu}_{i,\pntt} \in R \tens k(\gnc) 
\hspace{20mm} \forall \,\nmr, \nmb, n, \nu, i, \pntt 
\] 
are defined over $k(\gnc)$. 

\bClm 
\label{bound lh(b)}
If the set 
$\bigst{\lgth\bigpn{\sD_{\rC(b)}^{\,\inf}} \;\big|\; b \in \Base}$ 
admits a bound, then the sequence 
$\bigpn{\sD_{\rC(\gnc),\nmb}^{\,\inf}}_{\nmb \geq 1}$ 
becomes stationary. 
\eClm 

\bBf 
For every closed point $b \in \Base$ there exists a number \,$\mty(b)$\, 
such that the sequence $\bigpn{\sD_{\rC(b),\nmb}^{\,\inf}}_{\nmb \geq 1}$ 
becomes stationary after \,$\mty(b)$\, elements. 
By Lemma \ref{bounds} we see that a bound of 
$\bigst{\lgth\bigpn{\sD_{\rC(b)}^{\,\inf}} \;\big|\; b \in \Base}$, 
if it exists, 
yields a bound for the set $\st{\mty(b) \;|\; b \in \Base}$. 
In this case the sequence $\bigpn{\sD_{\rC(\gnc),\nmb}^{\,\inf}}_{\nmb \geq 1}$ 
becomes stationary, since it is compatible with specialization 
($\see$Claim \ref{BaseChange-for-Alg}). 
\eBf 

\bClm 
\label{bound lh(m)}
If the set 
$\bigst{\lgth\bigpn{\sD_{\rC(\gnc),\nmb}^{\,\inf}} \;\big|\; \nmb \geq 1}$ 
admits a bound, then the sequence 
$\bigpn{\sD_{\rC(\gnc),\nmb}^{\,\inf}}_{\nmb \geq 1}$ 
becomes stationary. 
\eClm 

\bBf 
Suppose there is a bound $\lgth(\gnc)$ of 
$\bigst{\lgth\bigpn{\sD_{\rC(\gnc),\nmb}^{\,\inf}} \;\big|\; \nmb \geq 1}$. 
Then by Claim \ref{BaseChange-for-Alg} 
we have \,$\lgth\bigpn{\sD_{\rC(b),\nmb}^{\,\inf}} \leq \lgth(\gnc)$\, 
for all $\nmb \geq 1$ and all $b \in U_{\nmb} \subset \Base$. 
By Corollary \ref{finite subsystem} 
we may omit for every $b \in \Base$ certain critical $\mu \geq 1$ 
in the pro-basis of $\ACHmp{\rC(\gnccl)}{\mdl_{\rCgncclgst}}$, 
e.g.\ those $\mu$ such that $\Tha_{n,\nu,\sig,\mu}$ does not meet $\rC(b)$ 
properly. 
Then we obtain a modified sequence 
$\bigpn{{\sD}_{\rC(b),\nmb'}^{\,\inf}}_{\nmb' \geq 1}$ 
with limit $\sD_{\rC(b)}^{\,\inf}$ 
and such that \,$\lgth\bigpn{{\sD}_{\rC(b),\nmb'}^{\,\inf}} \leq \lgth(\gnc)$\, 
for all $\nmb' \geq 1$. 
This shows that $\bigst{\lgth\bigpn{\sD_{\rC(b)}^{\,\inf}} \;\big|\; b \in \Base}$ 
is bounded by $\lgth(\gnc)$ 
and hence $\bigpn{\sD_{\rC(\gnc),\nmb}^{\,\inf}}_{\nmb \geq 1}$ 
becomes stationary by Claim \ref{bound lh(b)}. 
\eBf  
 
%

\medskip 
Now suppose the sequence $\bigpn{\sD^{\,\inf}_{\rC(\gnc),\nmb}}_{\nmb}$ 
does not become stationary. 
Then, due to Claim \ref{bound lh(m)}, the set 
$\bigst{\lgth\bigpn{\sD_{\rC(\gnc),\nmb}^{\,\inf}} \;\big|\; \nmb \geq 1}$
is unbounded, and hence 
\bMyEqn 
\label{Wq} \mytag{Wq} 
\forall \,\iota \;\; \exists \,i_{\iota} \geq \iota 
\;\;\exists \,\nmb_{\iota} \geq \mty\pn{\iota} 
\;\;\exists \,\nmr_{\iota}, n_{\iota}, \nu_{\iota}, \pntt_{\iota}: 
\hspace{6mm} {}^{\nmb_{\iota}}_{\;\nmr_{\iota}} 
w^{n_{\iota},\nu_{\iota}}_{i_{\iota},\pntt_{\iota}} \neq 0 
\eMyEqn 
where $\mty(l)$ is the function from Lemma \ref{bounds}. 

Due to Lemma \ref{Base property} we may assume that 
$\Base$ is defined over $\bF_{q}$ 
and is a dense open subscheme of a projective space. 
For every integer $e \geq 1$ the set $\Base(e)$ 
of $\bF_{q^e}$-rational points of $\Base$ is finite, 
therefore we find a common bound 
\,$\lgth(e)$\, for $\bigst{\lgth\bigpn{\sD_{\rC(b)}^{\,\inf}} \;\big|\; b \in \Base(e)}$. 
Then we have for every $e, \nmr, n, \nu, \pntt$: 
\bMyEqn 
\label{WZ} \mytag{WZ} 
{}^{\nmb}_{\;\nmr} w^{n,\nu}_{i,\pntt} (b) = 
\left\{ 
\begin{array}{ll} 
0 \hspace{5mm} \textrm{or} & \forall \,\nmb \geq \mty\bigpn{\lgth(e)}, \\ 
\textrm{not defined}  \hspace{15mm} & 
\forall \,i > \lgth(e), \;\forall \,b \in \Base(e) \laurin 
\end{array}
\right.  
\eMyEqn 
For any element \,$f \in R \tens k(\gnc)$\, define 
\[ \ety\pn{f} 
\,:=\, \sup \lrst{e \;\left|\; f(b) = 
\Big\{ 
\begin{array}{l} 
0 \hspace{5mm} \textrm{or} \\ 
\textrm{not defined}  
\end{array}
\hspace{3mm} \forall \,b \in \Base(e) 
\right. 
} 
\] 
and 
\[ 
\dty\pn{f} 
\,:=\, \min \lrst{ \deg \Es \;\left|\; 
\begin{array}{l} 
\textrm{$\Es$ an effective divisor on $\Base$}  \\ 
f \in R \tens \Gam\bigpn{\Base, \sO_{\Base}\pn{\Es}} 
\end{array} 
\right. 
} 
\laurin 
\] 
Note that \,$\ety(f) < \infty$, \,$\dty(f) < \infty$\, if $f \neq 0$, 
and there is a strictly increasing function \,$\hty: \Nat \ra \Nat$\, 
(automatically unbounded) 
such that 
\bMyEqn \label{ED} \mytag{ED} 
\forall \,f \in R \tens k(\gnc): \hspace{10mm} 
\ety\pn{f} \geq n \hspace{5mm} \Lra \hspace{5mm} \dty\pn{f} \geq \hty\pn{n} 
\laurin 
\eMyEqn  
Then (\ref{WZ}) is expressed by: $\;\forall \,e, \nmr, n, \nu, \pntt:$ 
\bMyEqn \label{WE} \mytag{WE} 
\ety\bigpn{{}^{\nmb}_{\;\nmr} w^{n,\nu}_{i,\pntt}} \geq e 
\hspace{15mm} \forall \,i > \lgth(e), \;\forall \,\nmb \geq \mty\bigpn{\lgth(e)}  
\laurin 
\eMyEqn 
Then by (\ref{WE}), in the notation of (\ref{Wq}), the set 
\[ \bigst{\ety\bigpn{{}^{\nmb_{\iota}}_{\;\nmr_{\iota}} 
w^{n_{\iota},\nu_{\iota}}_{i_{\iota},\pntt_{\iota}}} \;\big|\; \iota \geq 1} \subset \Nat 
\hspace{7mm} \textrm{is \,unbounded} 
\laurink 
\] 
and due to (\ref{ED}), 
the same is true for the set 
\,$\bigst{\dty\bigpn{{}^{\nmb_{\iota}}_{\;\nmr_{\iota}} 
w^{n_{\iota},\nu_{\iota}}_{i_{\iota},\pntt_{\iota}}} \;\big|\; \iota \geq 1} \subset \Nat$. 
On the other hand, 
the procedure (\ref{Alg}) from Algorithm \ref{algorithm} 
for solving (\ref{matrix}) 
shows that any common bound for $\dty\pn{u_{n,\nu,\sig,\mu}}$ 
yields a common bound for $\dty\pn{{}^{\nmb}_{\;\nmr} w^{n,\nu}_{i,\pntt}}$, 
a contradiction. 
Therefore, in order to conclude the proof of Lemma \ref{pro-divisor terminates}, 
it suffices to show 

\bClm 
\label{bound d(u)}
The set 
$\bigst{\dty\pn{u_{n,\nu,\sig,\mu}} \;\big|\; n,\nu,\sig,\mu}$ 
is bounded. 
\eClm 

\bBf 
The pro-basis $\st{\tha_{n,\nu,\sig,\mu}}_{n,\nu,\sig,\mu} \in \sK_X^*$ 
is defined over $k$ and hence over a finite extension $\bF_{q^{\oma+1}}$ 
of $\bF_q$ for some $\oma \geq \mu$. 
Then \,$\Tha_{n,\nu,\sig,\mu} = \dv\pn{\tha_{n,\nu,\sig,\mu}}_X$\, 
is defined over $\bF_{q^{\oma+1}}$ and 
\[ \Trace \Tha_{n,\nu,\sig,\mu} 
\,=\, \sum_{j=0}^{\oma} \Frob_q^{\,j} \Tha_{n,\nu,\sig,\mu} 
\,=\, \dv \Biggpn{ \prod_{j=0}^{\oma} \pn{\Frob_q^{\,j}}^* {\tha_{n,\nu,\sig,\mu}} }_X 
\] 
where $\Frob_q$ is the $q$-power Frobenius, 
is defined over $\bF_q$. 

Let \,$X \inj \Prj^m$\, be an embedding into a projective space, 
and \,$\mdll = \Zero(\eps)$\, an effective divisor on $\Prj^m$, 
defined over $\bF_q$, of degree $\leq \deg \mdl$, 
meeting $X$ properly and such that $\mdl \leq \mdll \cut X$. 
The projective embedding of $X$ yields an embedding 
\,$\rC \inj \Prj^m \tms \Base$\, into a projective space over $\Base$, 
and hence a projective embedding \,$\rC(\gnc) \inj \Prj_{\fld(\gnc)}^m$. 
By Lemma \ref{degZf}
we find a function \,$\del_{\rC(\gnc)} \in \sK_{\rC(\gnc)}^*$ with 
\[ 
\del_{\rC(\gnc)} 
\,\equiv\, \mfn^* \prod_{j=0}^{\oma} \pn{\Frob_q^{\,j}}^* {\tha_{n,\nu,\sig,\mu}} 
\mod \mdl_{\rCgncgst} 
\] 
such that the divisor of zeroes and divisor of poles of $\del_{\rC(\gnc)}$ 
are of degree 
\[ \deg\Zero\pn{\del_{\rC(\gnc)}} \,=\, \deg\Pole\pn{\del_{\rC(\gnc)}} 
\,\leq\, \deg\bigpn{\mdll \cut \rC(\gnc)} 
\,\leq\, \deg\mdl \cdot \deg\rC(\gnc) 
\laurin 
\] 
Since $\mdll$ is a divisor on $\Prj^m$, 
i.e.\ the defining function $\eps$ of $\mdll$ lives on $\Prj^m \supset X$, 
the function $\del_{\rC(\gnc)} \in \sK_{\rC(\gnc)}$ 
also comes from $X$ 
(due to construction as in the proof of Lemma \ref{degZf}), 
i.e.\ there is a function $\del \in \sK_X^*$ 
with $\del_{\rC(\gnc)} = \mfn^*\del|_{\rC(\gnc)}$. 
Then the divisor of zeroes and divisor of poles of $\del$ 
are of degree 
\[ \deg\Zero(\del) \,=\, \deg\Pole(\del) 
\,\leq\, \deg\mdl \cdot \deg X 
\laurin 
\] 
The divisor of $\del$ 
\[ \Del \,:=\, \dv(\del) 
\,=\, \Zero(\del) - \Pole(\del) 
\] 
is also defined over $\bF_q$. 
Then 
there is a dense open $U \subset \Base$ such that 
\begin{align*} 
{\iot{\Trace \Tha_{n,\nu,\sig,\mu}}}_* 
  \bigpn{\Trace \Tha_{n,\nu,\sig,\mu} \cdot \rC(b)} 
&= {\iot{\Del}}_* \bigpn{\Del \cdot \rC(b)} 
\hspace{20mm} \forall \,b \in U 
\\ 
& \in \ACHm{X}{\mdl} 
\end{align*} 
and hence by compatibility of skeleton divisors ($\see$Remark \ref{Char_Em}) 
\begin{align*} 
\bigpair{\sD_{\wt{\Trace \Tha_{n,\nu,\sig,\mu}}}, 
\wt{\Trace \Tha_{n,\nu,\sig,\mu}} \cut \rC(b)}_{\vrt{\wt{\Trace \Tha_{n,\nu,\sig,\mu}}},\mdl} 
&= \bigpair{\sD_{\wt{\Del}}, \wt{\Del} \cut \rC(b)}_{\vrt{\wt{\Del}},\mdl} \\ 
\parallel \hspace{45mm} & \hspace{15mm} \parallel \\ 
u_{\,\Trace \Tha_{n,\nu,\sig,\mu}} \pn{b} \hspace{30mm} 
& \hspace{13.8mm} u_{\Del} \pn{b} 
\hspace{20mm} \forall \,b \in U 
\end{align*} 
which yields \;$u_{\,\Trace \Tha_{n,\nu,\sig,\mu}} = u_{\Del}$. 
Write 
\[ \Tha^{\,j} \,:=\, \Frob_q^{\,j} \Tha_{n,\nu,\sig,\mu} 
\;\laurink 
\] 
then by Lemma \ref{Trace commutes} it holds 
\begin{eqnarray*} 
u_{\,\Trace \Tha_{n,\nu,\sig,\mu}} 
& = & \rmp{\sD_{\wt{\Trace\Tha_{n,\nu,\sig,\mu}}}}_{\Base}
 \Bigpn{\mfn^*\wt{\Trace\Tha_{n,\nu,\sig,\mu}} \tms_{\Base} \Opn_{n,\nu,\sig,\mu}}^{-1} \\ 
& = & \prod_{j = 0}^{\oma} \rmp{\sD_{\That^{\,j}}}_{\Base}
      \Bigpn{\mfn^*\That^{\,j} \tms_{\Base} \Opn_{n,\nu,\sig,\mu}}^{-1} \\ 
& = & \prod_{j = 0}^{\oma} u_{\,\Tha^{\,j}} 
\end{eqnarray*}
and hence 
\bMyEqn \label{EU} \mytag{EU} 
\dty\bigpn{u_{n,\nu,\sig,\mu}} 
\,\leq\, \dty\Biggpn{\prod_{j = 0}^{\oma} u_{\,\Tha^{\,j}}} 
\,=\, \dty\bigpn{u_{\Del}} 
\laurin 
\eMyEqn 
Now there exist only finitely many effective divisors $\Es$ in $X$ 
defined over $\bF_q$ and of degree $\leq \deg \mdl \cdot \deg X$. 
Therefore the set 
\[ \lrst{\, \dty\bigpn{u_{\Del}} \; \left| \; 
\begin{array}{l} 
\Del = \Es_1 - \Es_2, \hspace{3mm} \Es_i \,/\, \bF_q \\ 
\deg \Es_i \leq \deg \mdl \cdot \deg X 
\end{array} 
\right. }
\] 
is bounded, 
which by (\ref{EU}) implies that 
$\bigst{\dty\bigpn{u_{n,\nu,\sig,\mu}} \;\big|\; n,\nu,\sig,\mu}$ 
is bounded as well. 
\eBf 

\medskip 
The proof of Claim \ref{bound d(u)} finishes the proof of 
Lemma \ref{pro-divisor terminates}. 
\ePf 

\bLem 
\label{Trace commutes}
Let $\pn{\sD_C}_C \in \varprojlim \Emo{C}{\mdl_{\Cgst}}(R)$, $R$ a finite $k$-ring, 
be a skeleton relative divisor on $X$ defined over $\bF_q$. 
For any divisor $\Tha_{\mu} \subset X$ defined over $\bF_{q^{\oma+1}}$, 
we form its trace as 
\[ \Trace \Tha_{\mu} \,=\, \sum_{j=0}^{\oma} \Frob_q^{\,j} \Tha_{\mu} 
\laurink 
\] 
where $\Frob_q$ is the geometric $q$-power Frobenius. 
Then for $\Tha_{\mu} = \dv(\tha_{\mu})$ 
and some open dense $\Opn \subset \Base$ 
it holds 
\[ \rmp{\sD_{\wt{\Trace\Tha_{\mu}}}}_{\Base}
      \Bigpn{\mfn^*\wt{\Trace\Tha_{\mu}} \tms_{\Base} \Opn} \\ 
 \,=\, \prod_{j = 0}^{\oma} \rmp{\sD_{\wt{\Frob_q^{\,j} \Tha_{\mu}}}}_{\Base}
      \Bigpn{\mfn^*\wt{\Frob_q^{\,j} \Tha_{\mu}} \tms_{\Base} \Opn} 
\laurin 
\] 
\eLem 

\bPf 
We show the statement by restriction to the fibres over $\Base$. 
By compatibility of skeleton divisors ($\see$Remark \ref{Char_Em}) 
we 
have 
\[ \bigpair{\sD_{\rC(b)}, \Trace\Tha_{\mu} \cut \rC(b)}_{\rC(b),\mdl} 
= \bigpair{\sD_{\wt{\Trace\Tha_{\mu}}}, \wt{\Trace\Tha_{\mu}} \cut \rC(b)}
_{\vrt{\wt{\Trace\Tha_{\mu}}},\mdl} 
\laurin 
\] 
We 
obtain 
\begin{eqnarray*} 
\rmp{\sD_{\wt{\Trace\Tha_{\mu}}}}_{\Base}
      \bigpn{\mfn^*\wt{\Trace\Tha_{\mu}}}(b) 
& = &  \rmp{\sD_{\wt{\Trace\Tha_{\mu}}}}
            \Bigpn{\wt{\Trace\Tha_{\mu}} \cut \rC(b)} \\ 
& = &  \rmp{\sD_{\rC(b)}}
            \Bigpn{\Trace\Tha_{\mu} \cut \rC(b)} \\ 
& = &  \rmp{\sD_{\rC(b)}}
            \Biggpn{\sum_{j=0}^{\oma} \Frob_q^{\,j} \Tha_{\mu} \cut \rC(b)} \\ 
& = &  \prod_{j=0}^{\oma} 
            \rmp{\sD_{\rC(b)}} 
            \Bigpn{\Frob_q^{\,j} \Tha_{\mu} \cut \rC(b)} \\ 
& = &  \prod_{j=0}^{\oma} 
            \rmp{\sD_{\wt{\Frob_q^{\,j} \Tha_{\mu}}}} 
            \biggpn{\wt{\Frob_q^{\,j} \Tha_{\mu}} \cut \rC(b)} \\ 
& = &  \prod_{j=0}^{\oma} 
            \rmp{\sD_{\wt{\Frob_q^{\,j} \Tha_{\mu}}}}_{\Base} 
            \biggpn{\mfn^*\wt{\Frob_q^{\,j} \Tha_{\mu}}}(b) \\ 
\end{eqnarray*}
\hfill 
\ePf 


\bLem
\label{Secn-independent}
Let \,$\pn{\sD_C}_C \in \varprojlim \Emo{C}{\mdl_{\Cgst}}(R)$\, be a 
skeleton divisor on $X$, where $R$ is a finite $k$-ring. 
Suppose \,$\sD_{\rC} \in \Divf_{\rC}^{\vrt{\mdl \cut \rC}}(R)$ 
is a relative Cartier divisor on $\rC \ra \Base$ 
satisfying \,$\sD_{\rC} \cut \rC(b) = \sD_{\rC(b)}$\, for all $b \in \Base$. 
Then there is a unique relative Cartier divisor 
$\sD_{X} \in \Divf_{X}^{\vrt{\mdl}}(R)$ such that $\sD_{\rC} = \mfn^*\sD_X$. 
\eLem 

\bPf 
The unique relative divisor $\sD_{\rC}$ on $\rC \ra \Base$ 
with the required property 
was constructed in the proof of Lemma \ref{locRelDiv} 
from the skeleton divisor $\pn{\sD_C}_C$ on $X$ 
via rigidity ($\see$Lemma \ref{Basis-determines_p}) 
by means of a pro-basis coming from $X$. 
This implies that the restriction of the pseudo-divisor 
$(\sO_{\rC}(\sD_{\rC}), \vrt{\sD_{\rC}}, s_{\sD_{\rC}} = \textrm{canonical section})$
associated with $\sD_{\rC}$ ($\seecite$\cite[Def.~2.2.1]{F})
to any fibre of $\mfn: \rC \ra X$ is constant. 

If $\secn: X \supset U \ra \rC$ is a rational section of $\mfn: \rC \ra X$, 
the restriction of $\mfn$ to the closure $Y := \ol{\secn U}$ 
of the image of $\secn$ in $\rC$ 
gives a birational morphism $\mfn|_Y: Y \ra X$. 
If $Y$ is not suitable, replace it by its suitabilization. 
Then the induced map 
$\mfn|_Y^*: \Divf_X \lra \Divf_Y$ is an isomorphism, 
since $\Divf_X$ is a birational invariant among suitable proper schemes 
($\see$Proof of Theorem \ref{SkelThm}, Step 3). 
Thus 
$\sD_{\rC}|_Y$ corresponds to a relative divisor $\sD_X$ on $X$. 
Since $\sD_{\rC}$ is constant on the fibres of $\mfn$, 
the divisor $\sD_X$ is independent of the choice of section $\secn$. 
Therefore $\sD_{\rC} = \mfn^* \sD_X$ 
comes from a unique relative Cartier divisor on $X$. 
\ePf

\newpage 

\section{Roitman Theorem with Modulus} 
\label{RoitmanModulus}

The goal of this section is to generalize the famous Roitman Theorem 
from \cite{Roit} to a situation with modulus. 
Let $X$ be a projective variety over a perfect field $\fld$. 
Let $\mdl$ be a modulus for $X$, 
as in Section \ref{ChowMod}. 

\bPnt 
\label{exactSeq}
Let $\Lm{X}{\mdl}$ be the affine part of the Albanese $\Albm{X}{\mdl}$ 
of $X$ of modulus $\mdl$. 
The Abel-Jacobi map \;$\abjacm{X}{\mdl}$\; of $X$ of modulus $\mdl$ 
induces a commutative diagram with exact rows 
\[ \xymatrix{ 
0 \ar[r] & \ACHm{X}{\mdl} \ar[r] \ar[d]^{\abjacmaff{X}{\mdl}} & 
\CHmo{X}{\mdl} \ar[r] \ar[d]^{\abjacm{X}{\mdl}} & 
\CHmo{X}{0 \mdl} \ar[r] \ar[d]^{\abjacm{X}{0 \mdl}} & 0 \phantom{\laurin} \\ 
0 \ar[r] & \Lm{X}{\mdl}\pn{k} \ar[r] & \Albm{X}{\mdl}\pn{k} \ar[r] & 
\Alba{X}\pn{k} \ar[r] & 0 \laurin 
}
\] 
\ePnt 

\bThm 
\label{affineRoitmanThm}
Assume the base field $\fld$ is finite or an algebraic closure of a finite field 
or algebraically closed of characteristic 0. 
Let $X$ be a suitable projective variety over $\fld$. 
If $k$ is not algebraically closed, 
we assume $X$ to be smooth outside of $\vrt{\mdl_X}$. 
Then the affine Abel-Jacobi map 
is  an isomorphism: 
\bDpl  \abjacmaff{X}{\mdl}: \ACHm{X}{\mdl} \iso \Lm{X}{\mdl}\pn{\fld} 
\laurin 
\eDpl 
\eThm 

\bPf 
Let $C$ be a curve in $X$. 
Since the Abel-Jacobi maps for $X$ and $\Ct$ 
are compatible with the universal maps for $X$ and $\Ct$ 
by Corollary \ref{Alb=QuotCH}, 
the diagram of functoriality maps 
\[ \xymatrix{ 
\CHmo{\Ct}{\mdl_{\Ct}} \ar[d] \ar[r]^-{\sim} & \Albm{\Ct}{\mdl_{\Ct}}\pn{k} \ar[d] \\ 
\CHmo{X}{\mdl} \ar[r] & \Albm{X}{\mdl}\pn{k} 
}
\] 
commutes, and hence so does the diagram 
\[ \xymatrix{ 
\ACHm{\Ct}{\mdl_{\Ct}} \ar[d] \ar[r]^-{\sim} & \Lm{\Ct}{\mdl_{\Ct}}\pn{k} \ar[d] \\ 
\ACHm{X}{\mdl} \ar[r] & \Lm{X}{\mdl}\pn{k} \laurin 
}
\] 
Now $\Lm{\Ct}{\mdl_{\Ct}}$ 
resp.\ $\Lm{X}{\mdl}$ 
is the Cartier dual of $\Fmo{\Ct}{\mdl_{\Ct}}$ resp.\ $\Fmor{X}{\mdl}$. 
Therefore the affine Abel-Jacobi map \,$\abjacmaff{\Ct}{\mdl_{\Ct}}$ of $\Ct$ 
is given by 
\[ \abjacmaff{\Ct}{\mdl_{\Ct}}: \ACHm{\Ct}{\mdl_{\Ct}} \iso \Lm{\Ct}{\mdl_{\Ct}}\pn{k} \;, 
\hspace{10mm} 
\kap \lmt \pair{\llul,\kap}_{\Ct,\mdl_{\Ct}} 
\] 
where \,$\pair{\llul,\lull}_{\Ct,\mdl_{\Ct}}: 
\Fmo{\Ct}{\mdl_{\Ct}} \tms \ACHm{\Ct}{\mdl_{\Ct}} \ra \Gm$\, 
is the pairing from Proposition \ref{duality_curve}. 
As $\ACHm{X}{\mdl}$ is the inductive limit of the $\AQHm{C}{\mdl_{\Cgst}}$ 
(see Proposition \ref{CHm_indLim}) 
and the pairing 
$\pair{\llul,\lull}_{X,\mdl}: \Fmor{X}{\mdl} \tms \ACHm{X}{\mdl} \lra \Gm$ 
from Proposition \ref{Pairing_X} 
was constructed via restriction to curves, 
this shows that the affine Abel-Jacobi map \,$\abjacmaff{X}{\mdl}$ of $X$ 
is also given by 
\[ \abjacmaff{X}{\mdl}: \ACHm{X}{\mdl} \lra \Lm{X}{\mdl}\pn{k} \;, 
\hspace{10mm} 
\kap \lmt \pair{\llul,\kap}_{X,\mdl}
\] 
and factors as 
\[ \ACHm{X}{\mdl} \iso \varinjlim_C \AQHm{C}{\mdl} 
\iso \varinjlim_C \Lm{C}{\mdl}\pn{k} \iso \Lm{X}{\mdl}\pn{k}
\] 
where the last map follows from Theorem \ref{SkelThm} by Cartier duality, 
and we used descent of the base field (Subsection \ref{sub: Descent}). 
\ePf 

\bDef
\label{strongSuitable}
$X$ is called \emph{strongly suitable for $\mdl$} 
if $X$ is suitable and it satisfies the condition of Remark \ref{Bloch}, 
namely 
$ \CHm{X}{0 \mdl} = \CHO{X} 
$.  
\eDef


\bLem 
\label{ajm-exSeq}
Assume $\fld$ is a finite field or an algebraic closure of a finite field 
or algebraically closed of characteristic 0. 
Let $X$ be a projective variety over $\fld$, 
strongly suitable for an effective divisor $\mdl$. 
If $k$ is not algebraically closed, 
we assume $X$ to be smooth outside of $\vrt{\mdl_X}$. 
Then we have 
\,$ \Kajm{X}{\mdl} = \Kaj{X} 
$,
where 
\begin{align*}
	\Kaj{X} & := \ker\bigpn{\abjac{X}: \CHOo{X} \lra \Alba{X}(k)} \\ 
	\Kajm{X}{\mdl} & := \ker\bigpn{\abjacm{X}{\mdl}: \CHmo{X}{\mdl} \lra \Albm{X}{\mdl}(k)} 
	\laurin 
\end{align*}
\eLem 

\bPf 
Due to \Point \ref{exactSeq} and the Snake Lemma, 
we obtain the following commutative kernel-cokernel diagram 
with exact rows and exact snake sequence 
\[ \hspace{0mm} 
\begin{tikzpicture}[text height=1.5ex, text depth=0.25ex] 
\node (a0) at (0,0) {$0$}; 
\node (a1) [right=of a0] {$\ACHm{X}{\mdl}$}; 
\node (a2) [right=of a1] {$\CHmo{X}{\mdl}$}; 
\node (a3) [right=of a2] {$\CHOo{X}$}; 
\node (a4) [right=of a3] {$0$}; 
\node (b0) [below=of a0] {$0$}; 
\node (b1) [below=of a1] {$\Lm{X}{\mdl}\pn{k}$}; 
\node (b2) [below=of a2] {$\Albm{X}{\mdl}\pn{k}$}; 
\node (b3) [below=of a3] {$\Alba{X}\pn{k}$}; 
\node (b4) [below=of a4] {$0$}; 
\node (k1) [above=of a1] {$0$}; 
\node (k2) [above=of a2] {$\Kajm{X}{\mdl}$}; 
\node (k3) [above=of a3] {$\Kaj{X}$}; 
\node (c1) [below=of b1] {$0$}; 
\node (c2) [below=of b2] {$0$}; 
\node (c3) [below=of b3] {$0$}; 
\draw[->] 
(k1) edge (k2) 
(k2) edge (k3) 
(k3) edge[out=0, in=180] (c1) 
(c1) edge (c2) 
(c2) edge (c3) 
(a0) edge (a1) 
(a1) edge (a2) 
(a2) edge (a3) 
(a3) edge (a4) 
(b0) edge (b1) 
(b1) edge (b2) 
(b2) edge (b3) 
(b3) edge (b4) 
(k1) edge (a1) 
(k2) edge (a2) 
(k3) edge (a3) 
(a1) edge (b1) 
(a2) edge (b2) 
(a3) edge (b3) 
(b1) edge (c1) 
(b2) edge (c2) 
(b3) edge (c3); 
\end{tikzpicture} 
\] 
where the Abel-Jacobi maps \,$\abjac{X}$\, and \,$\abjacm{X}{\mdl}$\, 
are surjective by Corollary \ref{Alb=QuotCH}, 
while the map on affine parts \,$\abjacmaff{X}{\mdl}$\, 
is an isomorphism due to Theorem \ref{affineRoitmanThm}. 
\ePf 

\bLem 
\label{ChowTorsion}
Assume $k$ is an algebraic closure of a finite field. 
Then the groups $\CHmo{X}{\mdl}$ and $\Albm{X}{\mdl}(k)$ 
are both torsion. 
\eLem 

\bPf 
For an algebraic group $G$, 
every $k$-rational point $\pntt \in G(k)$ 
has a finite residue field $\fld\pn{\pntt} = \finfld$, 
and $G\pn{\finfld}$ is finite, hence torsion. 
Thus $G(k)$ is torsion. 

Let $c \in \CHmo{X}{\mdl}$ and $z \in \ZOm{X}{\mdl}$ a representative. 
As $k$ is algebraically closed, 
there is an irreducible curve $C$ in $X$ containing $z$. 
Then $c$ has a preimage in $\CHmo{C}{\mdl}$. 
By \Point \ref{CHm(C)_algGrp} it holds 
\,$\CHmo{C}{\mdl} = \Jacm{C}{\mdl}(k)$, 
and according to the argument above this group is torsion. 
So $c$ is the image of a torsion element, 
showing that $\CHmo{X}{\mdl}$ is torsion. 
\ePf 

\bThm[Roitman Theorem with Modulus over an Algebraic Closure of a Finite Field] 
\label{RoitmanThm over finite field}
Assume $k$ is an algebraic closure of a finite field. 
Let $X$ be a projective variety over $\fld$, 
strongly suitable for an effective divisor $\mdl$. 
Then the Abel-Jacobi map of $X$ of modulus $\mdl$ 
is an isomorphism: 
\[  \abjacm{X}{\mdl}: \CHmo{X}{\mdl} \iso \Albm{X}{\mdl}\pn{\fld} \laurin
\] 
\eThm 

\bPf 
Since $X$ is suitable, it is in particular normal 
by property (S2) of Definition \ref{suitable_Def}. 
We have $\CHOo{X} = \CHO{X}^{\tor}$ and $\Alba{X}(k) = \Alba{X}(k)^{\tor}$ 
by Lemma \ref{ChowTorsion}. 
Then $\Kaj{X} = \ker\bigpn{\CHO{X}^{\tor} \lra \Alba{X}(k)^{\tor}} = 0$ by \cite[Thm.~6.1]{Gei}. 
According to the exact snake sequence of Lemma \ref{ajm-exSeq}, 
the group $\Kajm{X}{\mdl}$ then also vanishes, 
giving the statement. 
\ePf 

\bThm[Roitman Theorem with Modulus] 
\label{RoitmanThm}
Assume $k$ is an algebraic closure of a finite field 
or algebraically closed of characteristic 0. 
Let $X$ be a projective variety over $\fld$, 
strongly suitable for an effective divisor $\mdl$. 
Then the Abel-Jacobi map $\abjacm{X}{\mdl}$ 
is an isomorphism on torsion parts: 
\[  \CHm{X}{\mdl}^{\tor} \iso \Albm{X}{\mdl}\pn{\fld}^{\tor} \laurin
\] 
\eThm 

\bPf 
Taking torsion the diagram from \Point \ref{exactSeq} yields 
with Remark \ref{Bloch} 
\[ \hspace{-10mm} 
\xymatrix{ 
0 \ar[r] & \ACHm{X}{\mdl}^{\tor} \ar[r] \ar[d] & 
\CHm{X}{\mdl}^{\tor} \ar[r] \ar[d] & 
\CHO{X}^{\tor} \ar[r] \ar[d] & 
\ACHm{X}{\mdl} \tens \Qrat/\Zint \ar[d] \\ 
0 \ar[r] & \Lm{X}{\mdl}\pn{k}^{\tor} \ar[r] & 
\Albm{X}{\mdl}\pn{k}^{\tor} \ar[r] & 
\Alba{X}\pn{k}^{\tor} \ar[r] & 
\Lm{X}{\mdl}\pn{k} \tens \Qrat/\Zint \laurin 
}
\] 
If $X$ is suitable, it is in particular normal 
by property (S2) of Definition \ref{suitable_Def}. 
The Abel-Jacobi map $\abjac{X}$ is an isomorphism on torsion parts 
by \cite[Thm.~6.1]{Gei}, 
so the fourth vertical map $\CHO{X}^{\tor} \iso \Alba{X}(k)^{\tor}$ 
is an isomorphism. 
The second and fifth vertical maps 
are isomorphisms 
due to the isomorphism $\ACHm{X}{\mdl} \iso \Lm{X}{\mdl}\pn{k}$ 
from Lemma \ref{ajm-exSeq}. 
Then the third vertical map $\CHm{X}{\mdl}^{\tor} \lra \Albm{X}{\mdl}(k)^{\tor}$ 
is an isomorphism by the Five Lemma. 
\ePf 

\bCor[Albanese Kernel over a Finite Field] 
\label{AlbKernel-finFld}
Let $X$ be a smooth projective geometrically irreducible variety 
over a finite field $\fld$. 
Let $\NS{X}$ be the N\'eron-Severi group of $X$, 
regarded as a $\fld$-group scheme, 
and $\mu_n = \ker\pn{\llul^n: \Gm \ra \Gm}$ the $\fld$-group of 
$n^{\textrm{th}}$ roots of unity. 
Define 
\bMyEqn 
\label{SX} \mytag{SX} 
\Kaj{X} = \Hom\Bigpn{\varinjlim_{n} \Homabk\pn{\mu_n, \NS{X}}, \Qrat/\Zint} 
\laurin 
\eMyEqn 
The Abel-Jacobi map \,$\abjacm{X}{\mdl}$\, of $X$ of modulus $\mdl$ 
fits into an exact sequence 
\bMyEqn 
\label{SCA} \mytag{SCA} 
0 \lra \Kaj{X} \lra \CHmo{X}{\mdl} \lra \Albm{X}{\mdl}\pn{\fld} \lra 0 
\laurin 
\eMyEqn 
In particular, if the N\'eron-Severi group of $X$ is torsion-free, then 
the Abel-Jacobi map of $X$ of modulus $\mdl$ 
is an isomorphism of finite groups: 
\[  \abjacm{X}{\mdl}: \CHmo{X}{\mdl} \iso \Albm{X}{\mdl}\pn{\fld} 
\laurin
\] 
\eCor 

\bPf 
The Abel-Jacobi map with modulus $\abjacm{X}{\mdl}$ 
is surjective by Corollary \ref{Alb=QuotCH}, 
and has the kernel $\Kaj{X} = \ker(\abjac{X})$ 
by Lemma \ref{ajm-exSeq}. 
This yields the exact sequence (\ref{SCA}). 
The description (\ref{SX}) of the Albanese kernel $\Kaj{X}$ 
is found in \cite[Prop.~9]{KS83}. 
For the second statement it suffices to remark that 
$\pn{\NS{X}}^{\tor} = 0$\, implies \,$\Kaj{X} = 0$. 
\ePf

\newpage 

\section{Geometric Class Field Theory} 
\label{ReciprocityMap}

Let $\fld = \bF_q$ be a finite field, 
$\clfld$ an algebraic closure, 
$X$ a geometrically irreducible projective variety over $\fld$. 
Let $\mdl$ be a modulus for $X$, 
as in Section \ref{ChowMod}. 
In order to keep the notation simple 
we assume that the map $\va: X \ra \Va$, 
which is used to define the modulus $\mdl$ for $X$,  
is an embedding (although this would not be necessary). 
We assume $X \setminus \mdl$ to be smooth. 
Moreover we assume X to be suitable. 
This can always be achieved without changing $X \setminus \mdl$, 
according to Proposition \ref{suitable_exists}.\footnote{ 
If we do not want to require $X \setminus \mdl$ to be smooth, 
let $S$ be a reduced effective divisor on $\Va$ containing $\Sing(X)$, 
then replace $\mdl$ by $\mdl + 0 S$, 
which makes sense inside functorial expressions 
such as $\CHm{X}{\mdl}$ due to Definition \ref{DefCHoMod}. 
}

\bPnt 
\label{local Galois groups}
For every irreducible curve $C \subset X$ intersecting $\mdl$ properly 
write 
\[ \mdl_{\Ct} \,= \sum_{\pnt \in \vrt{\mdl \cut \Ct}} n_{C,\pnt} \bt{\pnt} 
\] 
and define 
\[ G_{C,\pnt} \,:=\, \Gal\bigpn{\wh{K}_{C,\pnt}^{\ab} \big| \wh{K}_{C,\pnt}} 
\,=\hspace{-2mm} \varprojlim_{\hspace{5mm}L | \wh{K}_{C,\pnt}} \hspace{-2mm} \Gal\bigpn{L \big| \wh{K}_{C,\pnt}}
\] 
where $\wh{K}_{C,\pnt}^{\ab}$ is the maximal abelian extension 
of the completion $\wh{K}_{C,\pnt}$ of the function field $K_{C}$ 
of $C$ at $\pnt$, 
and the projective limit ranges over all finite abelian extensions 
$L$ over $\wh{K}_{C,\pnt}$. 
Moreover, for $\nu \in \Nat$ define the \emph{ramification groups} 
\[ G_{C,\pnt}^{\,\nu} \,:=\hspace{-2mm} \varprojlim_{\hspace{5mm}L | \wh{K}_{C,\pnt}} \hspace{-2mm} \Gal\bigpn{L \big| \wh{K}_{C,\pnt}}^{\nu}
\] 
as the projective limit of the ramification groups of the finite quotients 
(defined as e.g.\ in \cite[IV, \S1]{S66}). 
There is a canonical homomorphism 
\[ G_{C,\pnt} \lra \fundGab{X \setminus \mdl} 
\] 
to the abelianized fundamental group $\fundGab{X \setminus \mdl}$ 
classifying abelian finite \'etale covers of $X \setminus \mdl$. 
\ePnt 

\bDef 
\label{bounded ramification}
A rational map \;$\morphism: Y \dra X$\; 
of irreducible projective varieties over $\fld$ 
is an \emph{abelian covering} 
if it induces an abelian Galois extension $K_Y|K_X$ of function fields. 
$\morphism$ is said to be of \emph{ramification bounded by $\mdl$}, 
if it is finite and \'etale over $X \setminus \mdl$, 
and for every curve $C$ in $X$ intersecting $\mdl$ properly, 
for every $\pnt \in \vrt{\mdl_{\Ct}}$ and for every integer $\nu \geq n_{C,\pnt}$ 
the higher ramification groups $G_{C,\pnt}^{\,\nu}$ have zero-image in 
$\Gal\pn{Y | X} := \Gal\pn{K_{Y}|K_{X}} = \Aut\pn{Y_U | U}$. 
\eDef 

\bDef 
\label{fundGroupMod}
The \emph{abelian fundamental group of $X$ of modulus $\mdl$} 
is the cokernel 
\[ \fundGabm{X}{\mdl} 
\,=\, \coker\Biggpn{\dsum_{C \iscp \mdl} 
\dsum_{\pnt \in \vrt{\mdl \cut \Ct}} G_{C,\pnt}^{\,n_{C,\pnt}} 
\lra \fundGab{X \setminus \mdl} } 
\] 
where $C$ ranges over all irreducible curves $C$ in $X$ 
intersecting $\mdl$ properly 
(cf.\ \cite[Section 3]{HH}). 
Then by definition $\fundGabm{X}{\mdl}$ 
classifies abelian coverings of $X$ of ramification bounded by $\mdl$. 
\eDef 

\bRmk 
\label{fundGroupUnr}
Let $\fundGabm{X}{0 \mdl}$ be the group defined by 
\[ \fundGabm{X}{0 \mdl} 
\,=\, \coker\Biggpn{\dsum_{C \iscp \mdl} 
\dsum_{\pnt \in \vrt{\mdl \cut \Ct}} I_{C,\pnt} \lra \fundGab{X \setminus \mdl} } 
\] 
where $C$ ranges over all irreducible curves $C$ in $X$ 
intersecting $\mdl$ properly, 
and $I_{C,\pnt} = G_{C,\pnt}^{\,0}$ is the inertia group of $G_{C,\pnt}$ 
for all $C,\pnt$. 

Suppose $X$ is smooth. 
By definition, the abelianized fundamental group $\fundGab{X}$ of $X$ 
classifies unramified abelian coverings of $X$. 
Therefore $\fundGab{X}$ coincides with $\fundGabm{X}{0 \mdl}$. 
\eRmk 

\bPnt 
\label{afundGgeom}
According to the descriptions from 
Definition \ref{fundGroupMod} and Remark \ref{fundGroupUnr} 
we have a short exact sequence 
\[ 0 \lra \afundGabm{X}{\mdl} \lra \fundGabm{X}{\mdl} 
\lra \fundGabm{X}{0 \mdl} \lra 0 
\] 
where the kernel \,$\afundGabm{X}{\mdl}$ 
of the quotient map \,$\fundGabm{X}{\mdl} \lra \fundGabm{X}{0 \mdl}$ 
is given by 
\[ \afundGabm{X}{\mdl} 
\,=\, \im\Biggpn{\dsum_{C \iscp \mdl} 
\dsum_{\pnt \in \vrt{\mdl \cut \Ct}} 
\frac{I_{C,\pnt}}{G_{C,\pnt}^{\,n_{C,\pnt}}} } 
\laurin 
\] 
\ePnt 

\bPrp 
\label{CHm-to-fundGm}
There is a canonical homomorphism with dense image 
\[ \recm{X}{\mdl}: \CHm{X}{\mdl} \lra \fundGabm{X}{\mdl} 
\laurin 
\] 
In other words: an abelian covering \,$Y \dra X$\, 
is of ramification bounded by $\mdl$ 
if and only if \,$\Gal\pn{Y | X}$ is a quotient of \,$\CHm{X}{\mdl}$. 
\ePrp 

\bPf 
(Cf.~\cite[Section~3]{Hi14}.) 
Consider the Wiesend class group ($\see$\cite[Section 7]{KeSc09})
\[ \WCG{X \setminus \mdl} 
\,=\, \coker\Biggpn{\dsum_{C \iscp \mdl} K_C^* \lra \ZOm{X}{\mdl} 
\oplus \dsum_{C \iscp \mdl} 
\dsum_{\pnt \in \vrt{\mdl \cut \Ct}} \wh{K}_{C,\pnt}^*} 
\laurin 
\] 
One easily receives a canonical isomorphism 
\begin{align*}
\CHm{X}{\mdl} 
&= \coker\Biggpn{\dsum_{C \iscp \mdl} 
\bigcap_{\pnt \in \vrt{\mdl \cut \Ct}} U_{\Ct,\pnt}^{\,n_{\pnt}} 
\lra \ZOm{X}{\mdl}  } \\ 
&\cong \coker\Biggpn{\dsum_{C \iscp \mdl}
\dsum_{\pnt \in \vrt{\mdl \cut \Ct}} \wh{U}_{C,\pnt}^{\,n_{\pnt}} 
\lra \WCG{X \setminus \mdl} } 
\end{align*}
where \,$U_{\Ct,\pnt}^{\,n} = 1 + \fm_{\Ct,\pnt}^{n}$\, 
denotes the $n^{\textrm{th}}$ higher unit group in $K_{\Ct}^*$ for $n > 0$, 
$U_{\Ct,\pnt}^{\,0} = \sO_{\Ct,\pnt}^{*}$, 
respectively \,$\wh{U}_{C,\pnt}^{\,n} = 1 + \mc_{C,\pnt}^{n}$\, 
denotes the $n^{\textrm{th}}$ higher unit group in $ \wh{K}_{C,\pnt}^*$ 
and $\wh{U}_{C,\pnt}^{\,0} = \Oc_{C,\pnt}^{*}$. 
Here the Approximation Lemma ($\see$\cite[I, \S3]{S66}) 
is used for each $\Ct$ in this expression.  
According to Wiesend's class field theory ($\see$\cite[Prop.\ 7.5]{KeSc09}), 
there is a homomorphism (called \emph{reciprocity map}) 
\[ \rec{X \setminus \mdl}: \WCG{X \setminus \mdl} \lra \fundGab{X \setminus \mdl} 
\] 
induced by the Artin reciprocity map,  
and $\rec{X \setminus \mdl}$ has dense image, see \cite[Prop.~4.10]{Ke}. 
By local class field theory ($\seecite$\cite[XV, \S2, Thm.~2 and Rmk.]{S66}), 
the local Artin map yields isomorphisms 
\bMyEqn \label{UG} \mytag{UG} 
\wh{U}_{C,\pnt}^{\,n} \iso G_{C,\pnt}^{\,n} 
\laurin 
\eMyEqn 
Thus the reciprocity map $\rec{X \setminus \mdl}$ 
yields the desired homomorphism $\recm{X}{\mdl}$ with dense image 
on the quotients. 
\ePf 

\bThm 
\label{unramified}
The reciprocity map from Proposition \ref{CHm-to-fundGm} 
induces an isomorphism of finite groups 
\[ \recmo{X}{0 \mdl}: \CHmo{X}{0 \mdl} \iso \fundGgeom{X}{0 \mdl} 
\laurin 
\] 
\eThm 

\bPf 
Due to the log definitions, we have the following equalities: 
\[ \CHm{X}{\mdl_{\red}} \,=\, \WCGt{X \setminus \mdl} 
\] 
and 
\[ \fundGabm{X}{\mdl_{\red}} \,=\, \fundGabt{X \setminus \mdl} 
\] 
where $\WCGt{X \setminus \mdl}$ is the tame class group of Wiesend 
and $\fundGabt{X \setminus \mdl}$ the abelianized tame fundamental group. 
Moreover, we have an isomorphism of finite groups 
\[ \WCGto{X \setminus \mdl} \,\iso\, \fundGabto{X \setminus \mdl} 
\] 
($\seecite$\cite[Thm.~8.3]{KeSc09}). 
Using the isomorphism (\ref{UG}), we obtain 
\begin{eqnarray*} 
\frac{\CHmo{X}{\mdl_{\red}}}{\dsum_{C \iscp \mdl} 
\dsum_{\pnt \in \vrt{\mdl \cut \Ct}} U_{C,\pnt}^{\,0}} 
 & \iso 
 & \frac{\fundGgeom{X}{\mdl_{\red}}}{\dsum_{C \iscp \mdl} 
\dsum_{\pnt \in \vrt{\mdl \cut \Ct}} G_{C,\pnt}^{\,0}} \\ 
 \parallel \hspace{18mm} & & \hspace{18mm} \parallel \\ 
\CHmo{X}{0 \mdl} \hspace{7mm} & & \hspace{8mm}\fundGgeom{X}{0 \mdl} 
\laurin 
\end{eqnarray*}
The objects in the lower row of the diagram are independent of 
log and non-log version, 
therefore the assertion holds for both versions. 
\ePf 

\vspace{\vs} 

The following theorem relies 
on Lang's class field theory for function fields over finite fields, 
as described in \cite[VI]{S59}. 

\bThm 
\label{reciprocityLaw}
Let \;$\morphism: Y \dra X$\; be an abelian covering 
of suitable projective geometrically connected varieties 
over the finite field $\fld$, 
finite and \'etale outside a closed proper subset $S$ of $X$. 
Then there exists a modulus $\mdl$ for $X$ 
giving a bound for the ramification of \;$Y \dra X$, 
and we have a canonical isomorphism 
\[ 
\Gal\pn{Y|X} \;\cong\; 
\frac{\Albm{X}{\mdl}\pn{\fld}}{\morphism_* \Albm{Y}{\mdl_{\Ygst}}\pn{\fld}} 
\] 
where $\Gal\pn{Y|X}$ denotes the Galois group of the extension of function fields 
$K_Y | K_X$. 
\eThm 

\bPf 
Every abelian covering 
that arises from a ``geometric situation'' 
(i.e.\ does not arise from an extension of the base field \footnote{ 
$Y = X \tens_{\fld} \extfld$ would not be geometrically connected 
for $\extfld \supsetneqq \fld$, 
as $\extfld \tens_{\fld} \clfld$ is not simple.}) 
is the pull-back of a separable isogeny $H \ra G$ 
via a rational map $\phe: X \dra G$ 
(see \cite[VI, No.~8, Cor.\ of Prop.~7]{S59}). 
Using the universal mapping property of the Albanese variety with modulus, 
there is a homomorphism $g: \Albm{X}{\mdl} \lra G$ such that 
$\phe$ factors as $X \dra \Albm{X}{\mdl} \lra G$, 
where $\mdl$ is such that $\mdl_X \geq \modu\pn{\phe}$. 
Replacing $H \ra G$ by the pull-back $g^*H =: H_{Y|X} \lra \Albm{X}{\mdl}$, 
we may thus assume that $\morphism: Y \dra X$ is the pull-back of 
a separable isogeny $i_{Y|X}: H_{Y|X} \lra \Albm{X}{\mdl}$. 
We have $\Gal\pn{Y|X} = \ker\pn{i_{Y|X}}$. 
The 
isogeny $i_{Y|X}$ is a quotient of the 
``$q$-power Frobenius minus identity'' 
$\wp := \Frob_q - \id$ 
(see \cite[VI, No.~6, Prop.~6]{S59}), 
i.e.\ there is a homomorphism $h_{Y|X}: \Albm{X}{\mdl} \lra H_{Y|X}$ 
with $\wp = i_{Y|X} \circ h_{Y|X}$. 
Let $X_D \dra X$ be the pull-back of 
$\wp: \Albm{X}{\mdl} \lra \Albm{X}{\mdl}$. 
Then $\Gal\pn{X_{\mdl}|X} = \ker\pn{\wp} = \Albm{X}{\mdl}\pn{\fld}$. 

The universal property of $\Albm{Y}{\mdl_{\Ygst}}$ implies that 
the map $Y \dra H_{Y|X}$ factors through 
$\albm{Y}{\mdl_{\Ygst}}: Y \dra \Albm{Y}{\mdl_{\Ygst}}$. 
Let $H_{X_{\mdl}|Y}$ be the fibre product of 
$\Albm{X}{\mdl}$ and $\Albm{Y}{\mdl_{\Ygst}}$ over $H_{Y|X}$. 
Due to the universal property of the fibre-product, 
the commutative square 
\[ \xymatrix{ X_{\mdl} \ar[d] \ar[r]  &  \Albm{X}{\mdl} \ar[d] \\ 
Y \ar[r]  &  H_{Y|X} 
} 
\] 
factors as 
\[ \xymatrix{ X_{\mdl} \ar[d] \ar[r]  &  
H_{X_{\mdl}|Y} \ar[d]_{i_{X_{\mdl}|Y}} \ar[r] \ar@{}[dr]|{\square}  &  
\Albm{X}{\mdl} \ar[d] \\ 
Y \ar[r]  &  \Albm{Y}{\mdl_{\Ygst}} \ar[r]  &  H_{Y|X} \laurin 
} 
\] 
Here again the isogeny \,$i_{X_{\mdl}|Y}$ is a quotient of 
\,$\wp: \Albm{Y}{\mdl_{\Ygst}} \lra \Albm{Y}{\mdl_{\Ygst}}$. \\ 
Let \,$Y_{\mdl} \dra Y$\, be its pull-back. 

The functoriality of $\Albm{X}{\mdl}$ yields a commutative square 
\[ \xymatrix{ Y \ar[rr]^-{\albm{Y}{\mdl_{\Ygst}}} \ar[d]_{\morphism}  &&  
\Albm{Y}{\mdl_{\Ygst}} \ar[d]^{\morphism_*} \phantom{\laurin} \\ 
X \ar[rr]^-{\albm{X}{\mdl}}  &&  \Albm{X}{\mdl} \laurin 
} 
\] 
We obtain the following commutative diagram 
\[ \xymatrix{ 
Y_{\mdl} \ar[d]_{\morphism_{\mdl}} \ar[r]  &  
\Albm{Y}{\mdl_{\Ygst}} \ar[d]_{h_{X_{\mdl}|Y}} 
\ar@{-}@/^3pc/[dd]^(0.3){\Albm{Y}{\mdl_{\Ygst}}(\fld)}  &  \\ 
X_{\mdl} \ar[d] \ar[r]  & 
H_{X_{\mdl} | Y} \ar[d]_{i_{X_{\mdl}|Y}} \ar[r]^(0.39){j_{Y|X}} 
\ar@{-}@/^1.5pc/[d]^(0.35){\Gal(X_{\mdl}|Y)}  &  
\Albm{X}{\mdl} \phantom{\laurin} 
\ar[d]_{h_{Y|X}} \ar@{-}@/^3pc/[dd]^(0.3){\Albm{X}{\mdl}(\fld)} \\ 
Y \ar[d]_{\morphism} \ar[r]  &  
\Albm{Y}{\mdl_{\Ygst}} \ar[d]_{\morphism_*} \ar[r]  &  
H_{Y|X} \ar[d]_{i_{Y|X}} \ar@{-}@/^1.5pc/[d]^(0.4){\Gal(Y|X)} \\ 
X \ar[r]  &  \Albm{X}{\mdl} \ar@{=}[r]  &  \Albm{X}{\mdl} 
}
\] 
(where a caption on an arch denotes the Galois group of the covering). 



\bClm 
\label{Alb-factors}
$j_{Y|X} \circ h_{X_{\mdl}|Y} = \morphism_*$ 
(the functoriality map: $\Albm{Y}{\mdl_{\Ygst}} \lra \Albm{X}{\mdl}$) 
\eClm 

\noindent 
This yields the second isomorphism of the statement: 
\begin{eqnarray*}
\Gal\pn{Y|X} \; = \; \frac{\Gal\pn{X_{\mdl} | X}}{\Gal\pn{X_{\mdl} | Y}} 
 & = & \frac{\Gal\pn{X_{\mdl} | X}}
                   {\Gal\pn{Y_{\mdl} | Y} \big/ \Gal\pn{Y_{\mdl} | X_{\mdl}}} \\ 
 \parallel \hspace{12mm} & & \hspace{18mm} \parallel \\ 
 h_{Y|X} \Albm{X}{\mdl}\pn{\fld} 
 & = & \frac{\Albm{X}{\mdl}\pn{\fld}} 
                  {j_{Y|X} h_{X_{\mdl}|Y} \Albm{Y}{\mdl_{\Ygst}}\pn{\fld}} \\ 
 &  &  \hspace{18mm} \parallel \\ 
 &  &  \hspace{7mm} 
 \frac{\Albm{X}{\mdl}\pn{\fld}}{\morphism_* \Albm{Y}{\mdl_{\Ygst}}\pn{\fld}} 
\;\laurin 
\end{eqnarray*}
\hfill 
\ePf 

\bBf[Proof of Claim \ref{Alb-factors}] 
As the diagram commutes, we have 
\begin{align*} 
i_{Y|X} \; h_{Y|X} \; j_{Y|X} \; h_{X_{\mdl}|Y}  & \,=\,  
\morphism_* \; i_{X_{\mdl}|Y} \; h_{X_{\mdl}|Y}  \\ 
\parallel \hspace{12mm} &  \hspace{10.5mm} \parallel \\ 
\wp \;\; j_{Y|X} \; h_{X_{\mdl}|Y}  & \,=\, 
\; \morphism_* \; \wp 
\end{align*} 
Now $\wp := \Frob_q - \id$ 
commutes with regular maps defined over $\fld = \bF_q$. 
Thus $\wp$ in particular commutes with $j_{Y|X} \; h_{X_{\mdl}|Y}$, 
i.e.\ we have 
\bDpl  j_{Y|X} \; h_{X_{\mdl}|Y} \; \wp \; = \; \morphism_* \; \wp 
\eDpl 
and $\wp$ is surjective 
($\seecite$\cite[VI, No.~4, Prop.~3]{S59}). 
Hence \,$j_{Y|X} \; h_{X_{\mdl}|Y} = \morphism_*$, 
which proves the claim. 
\eBf 

\vspace{\vs} 

The following special case of Theorem \ref{reciprocityLaw} 
was already  treated in its proof: 

\bCor[Existence Theorem] 
\label{Alb=Gal}
There is a canonical isomorphism 
\[ 
\Albm{X}{\mdl}\pn{\fld} \;\iso\; \Gal\pn{X_{\mdl}|X} 
\] 
where $X_{\mdl} \dra X$ is the pull-back of the 
``$q$-power Frobenius minus identity'' morphism 
\;$\wp = \Frob_q - \id: \Albm{X}{\mdl} \lra \Albm{X}{\mdl}$. 
\eCor 

\bCor 
\label{Gal_maxAb}
Suppose $X$ is smooth and irreducible. 
Taking the limit over all moduli $\mdl$ for $X$ we obtain 
\[ 
\varprojlim_{\mdl} \Albm{X}{\mdl}\pn{\fld} \;\iso\; 
\Gal\bigpn{K^{\ab}_{X} \big| K_{X}\clfld} 
\] 
where $\Gal\bigpn{K^{\ab}_{\Xo} \big| K_{X}\clfld}$ is the geometric Galois group 
of the maximal abelian extension $K^{\ab}_{X}$ of the function field $K_{X}$ 
of $X$. 
\eCor 

\bLem 
\label{fundGm-sur-Albm}
There is a canonical epimorphism 
\[ \fundalbm{X}{\mdl}: \fundGgeom{X}{\mdl} \lsur \Albm{X}{\mdl}\pn{\fld} 
\] 
where $\fundGgeom{X}{\mdl}$ is the 
abelian geometric fundamental group of $X$ of modulus $\mdl$ 
which classifies abelian \'etale covers of $X \setminus \mdl$ 
of ramification bounded by $\mdl$ 
that do not arise from extending the base field. 
\eLem 

\bPf 
The abelian covering $X_{\mdl} \dra X$ from Corollary \ref{Alb=Gal} 
is obviously an \'etale covering of $X \setminus \mdl$ 
arising from a geometric situation. 
Since $\Gal\pn{X_{\mdl} | X} \cong \Albm{X}{\mdl}\pn{\fld}$ 
is a quotient of $\CHmo{X}{\mdl}$, 
this covering is of ramification bounded by $\mdl$. 
Thus $\Gal\pn{X_{\mdl} | X} \cong \Albm{X}{\mdl}\pn{\fld}$ 
is a quotient of $\fundGgeom{X}{\mdl}$. 
\ePf 

\bThm[Reciprocity Law] 
\label{reciprocityDiagram}
The reciprocity map from Proposition \ref{CHm-to-fundGm} 
yields an isomorphism of finite groups: 
\[ \recmo{X}{\mdl}: \CHmo{X}{\mdl} \iso \fundGgeom{X}{\mdl} 
\laurin 
\] 
\eThm 

\bPf 
Consider the diagram 
\[ \xymatrix{ 
0 \ar[r] & \ACHm{X}{\mdl} \ar[r] \ar[d]^{\arecmo{X}{\mdl}} 
\ar@/_4pc/[dd]_(0.66){\abjacmaff{X}{\mdl}} & 
\CHmo{X}{\mdl} \ar[r] \ar[d]^{\recmo{X}{\mdl}} & 
\CHmo{X}{0 \mdl} \ar[r] \ar[d]^{\reco{X}} & 0 \phantom{\laurin} \\ 
0 \ar[r] & \afundGgeom{X}{\mdl} \ar[r] \ar[d]^{\afundalbm{X}{\mdl}} & 
\fundGgeom{X}{\mdl} \ar[r] \ar[d]^{\fundalbm{X}{\mdl}} & 
\fundGgeom{X}{0 \mdl} \ar[r] \ar[d]^{\fundalb{X}} & 0 \phantom{\laurin} \\ 
0 \ar[r] & \Lm{X}{\mdl}\pn{k} \ar[r] & \Albm{X}{\mdl}\pn{k} \ar[r] & 
\Alba{X}\pn{k} \ar[r] & 0 \laurin 
}
\] 

The composition $\fundalbm{X}{\mdl} \circ \recmo{X}{\mdl}$ coincides with 
$\abjacm{X}{\mdl}$ according to \cite[VI, No.\ 24, Thm.\ 2]{S59}. 
Then $\afundalbm{X}{\mdl} \circ \arecmo{X}{\mdl} = \abjacmaff{X}{\mdl}$ 
is injective by Theorem \ref{affineRoitmanThm}, 
hence $\arecmo{X}{\mdl}$ is injective. 
As $\Lm{X}{\mdl}\pn{k}$ is finite, 
$\ACHm{X}{\mdl}$ is finite as well. 
In the right column $\reco{X}$ is an isomorphism of finite groups 
due to Wiesend's class field theory, see Theorem \ref{unramified}. 
In particular $\CHmo{X}{0 \mdl}$ is finite. 
Then $\CHmo{X}{\mdl}$, as an extension of finite groups, is finite. 
The map $\recmo{X}{\mdl}$ is dense by Proposition \ref{CHm-to-fundGm}, 
but $\CHmo{X}{\mdl}$ is finite, hence $\recmo{X}{\mdl}$ is surjective. 
As in addition $\reco{X}$ is injective, 
$\arecmo{X}{\mdl}$ is surjective by the Snake Lemma. 
So $\arecmo{X}{\mdl}$ is injective and surjective, 
thus an isomorphism. 
Then $\recmo{X}{\mdl}$ is an isomorphism by the Five Lemma. 
While $\CHmo{X}{\mdl}$ is a finite group, 
the same is true for the isomorphic group $\fundGgeom{X}{\mdl}$. 
\ePf 

\bCor 
\label{reciprocityLimit}
Let $S$ be the support of an effective divisor on $\Va$ 
such that $X \setminus S$ is smooth. 
Taking the limit over all moduli $\mdl$ for $X$ with support in $S$ 
we obtain 
\[ \varprojlim_{\substack{\mdl \\ \vrt{\mdl} \subset S}} \CHmo{X}{\mdl} 
\;\iso \varprojlim_{\substack{\mdl \\ \vrt{\mdl} \subset S}} \fundGgeom{X}{\mdl} 
\;\iso\; \fundGgeo{X \setminus S} 
\] 
where $\fundGgeo{X \setminus S}$ is the 
abelian geometric fundamental group of $X \setminus S$ 
which classifies abelian \'etale covers of $X \setminus S$ 
that do not arise from extending the base field. 
\eCor 

\bCor 
\label{maxAbGalois}
Suppose $X$ is smooth and irreducible. 
Taking the limit over all moduli $\mdl$ for $X$ we obtain 
\[ \varprojlim_{\mdl} \CHmo{X}{\mdl} \;\iso\; 
\varprojlim_{\mdl} \fundGgeom{X}{\mdl} \;\iso\; 
\Gal\bigpn{K^{\ab}_{X} \big| K_{X}\clfld} 
\] 
where $\Gal\bigpn{K^{\ab}_{\Xo} \big| K_{X}\clfld}$ is the geometric Galois group 
of the maximal abelian extension $K^{\ab}_{X}$ of the function field $K_{X}$ 
of $X$. 
\eCor 

\bCor 
\label{reciprocityQuotient}
Any abelian covering \;$\morphism: Y \dra X$\; 
of ramification bounded by $\mdl$ induces canonical isomorphisms 
\[ 
\frac{\CHmo{X}{\mdl}}{\morphism_* \CHmo{Y}{\mdl}} 
\;\iso\; \frac{\fundGgeom{X}{\mdl}}{\morphism_*\,\fundGgeom{Y}{\mdl}} 
\;\iso\; \Gal\pn{Y|X} 
\] 
where $\Gal\pn{Y|X}$ denotes the Galois group of the extension of function fields 
$K_Y | K_X$. 
\eCor 

\bPf 
Due to functoriality of $\CHm{X}{\mdl}$ and $\fundGgeom{X}{\mdl}$. 
\ePf 

\bCor 
\label{fundGm-Albm-exSeq}
The canonical map 
\;$\fundalbm{X}{\mdl}: \fundGgeom{X}{\mdl} \lsur \Albm{X}{\mdl}\pn{\fld}$\; 
from Lemma \ref{fundGm-sur-Albm} 
fits into an exact sequence 
\[ 0 \lra \Kaj{X} \lra \fundGgeom{X}{\mdl} \lra \Albm{X}{\mdl}\pn{\fld} \lra 0 
\laurin 
\] 
(For the definition of $\Kaj{X}$ see Corollary \ref{AlbKernel-finFld}.) 
\eCor 

\bPf 
Follows from the exact sequence of Corollary \ref{AlbKernel-finFld} 
and the Reciprocity Law \ref{reciprocityDiagram}. 
\ePf

\newpage 

\appendix
\section{Leitfaden}
\label{Leitfaden} 

The main results are presented here by giving a summary of each section. 
For precise definitions and statements the reader is referred to the main article. 

\medskip 
\textbf{Section~\ref{Sec:SuitableVar}:} Suitable Varieties. \\  
In this work 
we will make use of the generalized Albanese varieties with modulus from \cite{Ru13}, 
which are defined for smooth projective varieties over a perfect field. 
For any projective variety $X$ over a perfect field $\fld$ 
we seek a birational projective morphism $\mfn: Y \ra X$ such that 
$Y$ admits such a generalized Albanese. 
As we do not assume the existence of a resolution of singularities, 
we ask for weaker assumptions on $Y$ that still enable the construction from \cite{Ru13}. 
We find that the following conditions are sufficient: 
(S1) every rational map from $Y$ to an abelian variety $A$ 
extends to a morphism from $Y$ to $A$,  and 
(S2) $Y$ is normal. 
Those varieties will be called \emph{suitable} 
($\see$Definition \ref{suitable_Def}), 
and a birational projective morphism $\mfn: Y \ra X$ such that $Y$ is suitable 
will be called a \emph{suitabilization} of $X$. 
Every projective variety $X$ over a perfect field 
admits a suitabilization ($\see$Prop.\ \ref{suitable_exists}). 

\bigskip 
\textbf{Section \ref{ChowMod}:} Relative Chow Groups with Modulus. \\ 
Regarding Chow groups with modulus, 
there exist a \emph{log} version and a \emph{non-log} version. 
Each version has its advantages and disadvantages, 
but similar statements hold for both of them. 
We will formulate the theory in a unifying way 
that matches both cases at the same time. 

Let $X$ be a projective variety over a field. 
A \emph{modulus} $\mdl$ for $X$ is an effective Cartier divisor on $X$ 
in the non-log theory. 
In the log-theory, a modulus for $X$ is a Cartier divisor 
on a locally factorial projective variety $V$ that $X$ admits a morphism to 
(e.g.\ $V = X$ if $X$ itself is locally factorial, 
or $X \ra V$ could be an embedding into a projective space $V = \Prj^r$). 
This is necessary because we need $D$ to be both, 
a Cartier and a Weil divisor, see notion (\ref{L_C}) below. 

The \emph{relative Chow group of $0$-cycles on $X$ of modulus $\mdl$} 
is 
\[ \CHm{X}{\mdl} = \coker\Bigpn{\ZLm{X}{\mdl} \xra{\dv} \ZOm{X}{\mdl}} 
\] 
($\see$Definition \ref{DefCHoMod}, cf.\ also Notation \ref{DefC1}),  
where \,$\ZOm{X}{\mdl} = \Z_0(X \setminus \vrt{\mdl \cut X})$ 
is the group of $0$-cycles on $X$ that do not meet 
$\mdl \cut X$ (= the pull-back of $D$ to $X$), 
and \;$\ZLm{X}{\mdl} \subset 
\underset{C \subset X} \dsum K_{C}^*$\; 
is the subgroup generated by 
\[ \Rlnpm{X}{\mdl} 
= \lrst{ \pn{C,f} \left| 
	\begin{array}{l}
		C\textrm{ a curve in $X$ meeting } \vrt{\mdl \cut X} \textrm{ properly}  \\ 
		f \in K_{C}^* \hspace{2mm} \textrm{ s.t.\ } \hspace{2mm} 
		f \equiv 1 \mod \mdl_{\Ct} 
	\end{array}
	\right.
} 
\laurin 
\] 
Here $K_C$ is the function field of $C$, 
$\Ct$ the normalization of $C$, 
and $\mdl_{\Ct}$ is defined 
in the \emph{non-log} theory as 
\bMyEqn 
\label{NL_C} \mytag{NL} 
\mdl_{\Ct} := {\mdl \cut \Ct} 
\laurink 
\eMyEqn 
whereas in the \emph{log} theory 
\bMyEqn 
\label{L_C} \mytag{L} 
\mdl_{\Ct} := \pn{\mdl \cut \Ct}_{\red} + {\pn{\mdl-\mdl_{\red}} \cut \Ct} 
\eMyEqn 
where $()_{\red}$ denotes the reduced part of a divisor. 
The advantage of the non-log notion is its simpler definition 
and less technical effort in the case that $X$ is not locally factorial: 
$\mdl$ is a Cartier divisor on $X$ without any restrictions. 
On the other hand, the log notion is motivated by class field theory: 
the definition is made in such a way 
that the following implication is satisfied: 
\bMyEqn
\label{DCR} \mytag{DCR}
\mdl \;\textrm{ reduced} 
\hspace{8mm} \Lra \hspace{8mm} 
\mdl_{\Ct} \;\textrm{ reduced} 
\hspace{8mm} \forall \,C \subset X, \,C \iscp \mdl 
\eMyEqn 
(here the symbol $\iscp$ stands for \emph{proper intersection}). 
In class field theory, ``$\mdl$ reduced'' 
corresponds to ``tame ramification'', 
whereas ``$\mdl$ with higher multiplicities'' 
corresponds to ``wild ramification''. 
The implication says that this distinction is compatible with restriction to curves. 

There is a natural map \,$\CHm{X}{\mdl} \lra \CHm{X}{0 \mdl}$\, 
from the relative Chow group of modulus $\mdl$ 
to the Chow group of ``zero-modulus'' $0 \mdl$, 
which is formed in the same way as $\CHm{X}{\mdl}$, 
but with $\Rlnpm{X}{\mdl}$ replaced by 
\[ \Rlnpm{X}{0 \mdl} = \left\{ \bigpn{C,f}\left|\begin{array}{l}
	C\textrm{ a curve in $X$ meeting } \vrt{\mdl \cut X} \textrm{ properly} \\ 
	f \in \sK_C^* \hspace{2mm} \textrm{ s.t.\ } \hspace{2mm} 
	f \in \sO_{\Ct,\pnt}^* \hspace{6mm} \forall\, \pnt \in \mdl \cut \Ct 
\end{array}\right.\right\} 
\laurin 
\] 
If $X$ is smooth, then $\CHm{X}{0 \mdl}$ coincides with 
the usual Chow group $\CHO{X}$ 
by Bloch's ``easy moving lemma'' ($\see$Proposition \ref{Bloch}). 

For singular curves $C$ in $X$ we introduce 
a relative Chow group of $0$-cycles 
$\QHm{C}{\mdl}$ as a quotient of $\CHm{\Ct}{\mdl}$ 
modulo the kernel of the push-forward of cycles 
$\nu_*: \ZOm{\Ct}{\mdl} \lra \ZOm{C}{\mdl}$, 
where $\nu: \Ct \ra C$ is the normalization 
(Definition \ref{Def_CHmSing}). 
The \emph{affine part of \,$\CHm{X}{\mdl}$} is 
\[ \ACHm{X}{\mdl}= \ker\Bigpn{\CHm{X}{\mdl} \lra \CHm{X}{0 \mdl}} 
\] 
and similarly for the affine part $\AQHm{C}{\mdl}$ of $\QHm{C}{\mdl}$ 
($\see$Definition \ref{Def_ACHm}). 

If the base field is algebraically closed, 
or if $X$ is smooth outside $\vrt{\mdl_X}$ and the base field is finite, 
we obtain descriptions 
\begin{align*} 
	\CHm{X}{\mdl} &= \varinjlim_{C \iscp \mdl} \QHm{C}{\mdl_{\Cgst}} \\ 
	\ACHm{X}{\mdl} &= \varinjlim_{C \iscp \mdl} \AQHm{C}{\mdl_{\Cgst}} 
\end{align*} 
where $C$ ranges over all curves $C$ in $X$ 
that meet $\vrt{\mdl_X}$ properly 
($\see$Proposition \ref{CHm_indLim}). 

\newpage 
Now assume for simplicity that the base field $k$ is 
algebraically closed. 
Fix a projective embedding (so the degree of subschemes is defined). 
For the computation of $\ACHm{X}{\mdl}$ 
we can restrict to complete intersection (=: ci) curves of degree 
$\leq \degbnd:= \pn{\deg\mdl}^{\dim X - 1} \deg X$: 
\[ \ACHm{X}{\mdl}  =  
\sum_{\substack{C \subset X \,\textrm{ci curve} \\ C \iscp \mdl \\ \deg C \leq \degbnd}} 
\im\bigpn{\ACHm{C}{\mdl_{\Cgst}}} 
\] 
($\see$Proposition \ref{ACH bounded intersection}), 
and even to infinitesimal neighbourhoods of $\mdl$ on those curves 
($\see$Theorem \ref{ACH generated by truncations}).

\bigskip 
\textbf{Section \ref{Modulus}:} Modulus of a Rational Map to a Torsor. \\ 
Let $X$ be a normal proper variety over a perfect field $k$. 
Let $K$ be a discrete valuation field 
with normalized valuation $\val$ 
and residue field $k$. 
\bDef
\label{filtration_Witt}
We define the following filtration on the additive group of $K$: 
\[ \fil_{n}\Ga(K) = 
\bigst{ f \in K  \;\big|\;  \val\pn{f} \geq -n } 
\laurin 
\] 
Suppose $\chr(K) = p > 0$. 
Brylinski's filtration on the group of Witt vectors $\Witt_r$ of length $r$ 
with coefficients in $K$ 
from \cite[No.~1, Prop.~1]{Br} is 
\[ \fil_{n}\Witt_r(K) = 
\left\{ \lrpn{f_{r-1},\ldots,f_{0}} 
\;\left|\;  
\begin{array}{l} 
	f_i \in K, \quad \val\pn{f_i} \geq -n/p^i \\ 
	\forall \;0 \leq i \leq r-1 
\end{array}
\right.\right\} 
\laurin 
\] 
Matsuda's filtration from \cite[3.1]{Mda}, 
but shifted by 1, is 
\[ \filb_{n}\Witt_r(K) 
= \fil_{n-1}\Witt_r(K) + \Ver^{r-r'} \fil_{n}\Witt_{r'}(K) 
\] 
where \,$r' := \min\lrst{\ord_p(n), r-1}$ 
and $V$ is the Verschiebung. 

Let $\filf_{n}\Witt_r(K)$ be the saturation of $\fil_{n}\Witt_r(K)$ 
by means of the Frobenius $\Frob$ 
(see \cite[2.2]{KR10}), 
respectively $\filfb_{n}\Witt_r(K)$ the saturation of $\filb_{n}\Witt_r(K)$ 
(see \cite[4.7]{KR10}). 
Then $\filf_{n}\Witt_r(K)$ will be regarded as the ``\emph{log} filtration of $\Witt_r(K)$'', 
and $\filfb_{n}\Witt_r(K)$ as the ``\emph{non-log} filtration of $\Witt_r(K)$''. 
We will use geometric global versions $\filf_{\mdl}\Witt_r(\sK_X)$ resp.\ 
$\filfb_{\mdl}\Witt_r(\sK_X)$ of these filtrations, see \cite[Def.~3.2 and 3.8]{Ru13}, 
where $\mdl$ is an effective divisor on $X$. 
\eDef 

The \emph{modulus of a rational map $\phe$ from $X$ to 
	a torsor $P$ under a commutative smooth connected algebraic group} 
from \cite{KR10} is an effective divisor on $X$ 
\[ \modu\pn{\phe} = \sum_{\codim(\pnt)=1} \modu_{\pnt}(\phe) \;\ol{\st{\pnt}} 
\] 
(here the sum ranges over points of codimension 1, 
and $\ol{\st{\pnt}}$ denotes the closure of $y$ in $X$) 
with support on the locus where $\phe$ is not defined. 
As a rough idea, 
for $\pnt \in X$ of codimension $1$ with $\phe \notin P\pn{\sO_{X,\pnt}}$ 
the multiplicities are defined 
in the \emph{log} theory as 
\[ \modu_{\pnt}(\phe) 
= 1 + 
\begin{array}{c} 
	\textrm{``pole order at $\pnt$ of the unipotent part of $\phe$} \\ 
	\textrm{w.r.t.\ the \emph{log} filtration of $\Witt_r\pn{\sK_{X,\pnt}}$''} 
\end{array} 
\] 
whereas in the \emph{non-log} theory as 
\[ \modu_{\pnt}(\phe) 
= 
\begin{array}{c} 
	\textrm{``pole order at $\pnt$ of the unipotent part of $\phe$} \\ 
	\textrm{w.r.t.\ the \emph{non-log} filtration of $\Witt_r\pn{\sK_{X,\pnt}}$''} 
\end{array} 
\] 
($\see$Definition \ref{DefMod}). 
In particular, if the target of $\phe$ has a non-trivial unipotent part, 
then \,$\modu\pn{\phe}$ has always higher multiplicity in the log theory, 
while it might be reduced in the non-log theory. 
Thus the log definition is compatible with the modulus for curves 
in the sense of \cite{S59}, 
whereas the non-log definition might not always see the difference between 
tame and wild ramification.

\bigskip 
\textbf{Section \ref{subsec:AlbMod}:} Albanese Variety with Modulus. \\ 
Here we recall the construction from \cite{Ru13} 
as a higher dimensional analog 
to the generalized Jacobian with modulus of Rosenlicht-Serre: 
Let $X$ be a suitable projective variety over a perfect field, 
and $\mdl$ a modulus for $X$. 
Consider a rational map $\phe: X \dra G$ from $X$ 
to a commutative smooth connected algebraic group $G$.  
We say \emph{$\phe$ factors through $\CHm{X}{\mdl}$}, 
if its associated map on 0-cycles of degree 0 
factors through $\CHmo{X}{\mdl}$ 
(= degree $0$ part of $\CHm{X}{\mdl}$). 
A rational map $\phe: X \dra G$ of modulus $\leq \mdl$ 
always factors through $\CHm{X}{\mdl}$ 
($\see$Corollary \ref{factorCHm}). 
The category of all rational maps $\phe: X \dra G$ 
with \,$\modu\pn{\phe} \leq \mdl$\, 
admits a universal object 
\[ \albm{X}{\mdl}: X \dra \Albm{X}{\mdl} 
\] 
($\see$Theorem \ref{AlbMod-construction_2}), 
the algebraic group $\Albm{X}{\mdl}$ is called 
the \emph{Albanese variety of $X$ of modulus $\mdl$}. 
The universal map \,$\albm{X}{\mdl}$ factors through a canonical epimorphism 
\[ \abjacm{X}{\mdl}: \CHmo{X}{\mdl} \lsur \Albm{X}{\mdl}(\fld) 
\] 
($\see$Corollary \ref{Alb=QuotCH}), 
called \emph{Abel-Jacobi map on $X$ of modulus $\mdl$}. 

Supplements for suitable varieties and non-log version 
(that extend \cite{Ru13}) are added.

\bigskip 
\textbf{Section~\ref{AbelJacobiMap}:} Duality of Divisors and 0-Cycles. \\ 
Let $\Lm{X}{\mdl}$ be the affine part of the smooth connected 
algebraic group $\Albm{X}{\mdl}$ 
(in the canonical decomposition due to \cite{C}) 
and $\Fmor{X}{\mdl}$ its Cartier dual. 
Here $\Fmor{X}{\mdl}$ is represented by a formal subgroup of 
the group-sheaf $\Divf_X$ of relative Cartier divisors on $X$ 
($\see$Definition \ref{DefFm(X,D)}). 
We construct a pairing between $\Fmor{X}{\mdl}$ and 
the affine part of the relative Chow group with modulus $\ACHm{X}{\mdl}$ 
\bMyEqn 
\label{FA} \mytag{FA} 
\pair{\llul,\lull}_{X,\mdl}: \Fmor{X}{\mdl} \tms \ACHm{X}{\mdl} \lra \Gm 
\eMyEqn 
($\see$Definition \ref{Def_perfPair} and Proposition \ref{Pairing_X}). 
For a curve $X = C$, this pairing is induced by the local symbol 
$\pn{\llul,\lull}_{\pnt}: G\pn{K_C} \tms K_C^* \lra G\pn{k}$ 
from \cite[III, \S1]{S59}:
\[ \Bigpair{\sD,\bigbt{\dv(f)}}_{\Crv,\mdl} = \,\prod_{\pnt \in \vrt{\mdl}} \,\pn{\sD,f}_{\pnt} 
\] 
($\see$Proposition \ref{duality_curve}), 
for higher dimensional $X$ the pairing $\pair{\llul,\lull}_{X,\mdl}$ 
is obtained by restricting to curves 
($\see$Remark \ref{restrict to curves}). 

On a smooth proper curve $C$ we have the following explicit description 
\[ \Lm{C}{\mdl}(\fld) 
\;=\; \ACHm{C}{\mdl} 
\;=\; \frac{\prod_{\pnt \in \vrt{\mdl}} \fld(\pnt)^*}{\extfld^*} \tms 
\prod_{\pnt \in \vrt{\mdl}} \frac{1+\fm_\pnt}{1+\fm_\pnt^{n_\pnt}} 
\] 
($\see$\Point \ref{Curve-business}), 
where $\extfld := \H^0\pn{C,\sO_C}$ 
and $\mdl = \sum_{\pnt \in \vrt{\mdl}} n_{\pnt} \bt{\pnt}$. 

Assume the base field is an algebraic closure of a finite field 
or algebraically closed of characteristic 0. 
Extending the notion of $\Fmor{C}{\mdl}$ in such a way 
that it is defined for singular curves $C$ in $X$ 
($\see$\Point \ref{dual_ACHm}), 
we can consider compatible systems of relative Cartier divisors on curves 
$\pn{\sD_C}_C \in \varprojlim_C \Emo{C}{\mdl_{\Cgst}}$ 
($\see$Definition \ref{Def_skeletonDiv}), 
which we will call \emph{skeleton relative divisors} 
in analogy to the terminology \emph{2-skeleton sheaves} from \cite{EK}.
By means of a descent of the base field, we can carry over 
constructions and results to the case of a finite base field 
(Subsection \ref{sub: Descent}).

\bigskip 
\textbf{Section \ref{sec: Skeleton}:} Skeleton Theorem. \\ 
With the notion of a 
\emph{(pro-)basis of $\ACHm{C}{\mdl}$} 
($\see$Definition \ref{basis-ACHm} resp.\ \ref{pro-basis-ACHm}) 
we obtain: 
A relative Cartier divisor $\sD \in \Fmor{C}{\mdl}$ 
is uniquely determined by the values 
$\bigst{\bigpair{\sD, \alp_{\idx}}_{C,\mdl}}_{\idx}$ 
of the pairing (\ref{FA}) on a (pro-)basis 
$\st{\alp_{\idx}}_{\idx}$ of $\ACHm{C}{\mdl}$ 
($\see$Rigidity Lemma \ref{Basis-determines_0} resp.\ \ref{Basis-determines_p}). 

Assume that the base field $\fld$ is finite, 
an algebraic closure of a finite field 
or algebraically closed of characteristic $0$. 
We show that the natural homomorphism 
from $\Fmor{X}{\mdl}$ to $\varprojlim_C \Emor{C}{\mdl_{\Cgst}}$, 
given by pull-back of a relative Cartier divisor on $X$ 
to the various curves $C$ in $X$, 
is an isomorphism 
\[ \skelm{X}{\mdl}: \Fmor{X}{\mdl} \iso \varprojlim_C \Emor{C}{\mdl_{\Cgst}} 
\] 
($\see$Theorem \ref{SkelThm}, ``Skeleton Theorem''). 
Taking the limit over all moduli $\mdl$ for $X$, 
this constitutes a 1-1 correspondence between relative Cartier divisors on $X$ 
and compatible systems of relative Cartier divisors 
on the various curves of $X$ ($\see$Definition \ref{Def_skeletonDiv}). 
This is the most relevant (and most difficult) part of this work. 

\newpage

\bigskip 
\textbf{Section \ref{RoitmanModulus}:} Roitman Theorem with Modulus. \\ 
Let $X$ be a suitable projective variety, $\mdl$ a modulus for $X$. 
Assume that the base field $\fld$ is finite, 
an algebraic closure of a finite field 
or algebraically closed of characteristic $0$. 
If $\fld$ is not algebraically closed assume moreover 
that $X$ is smooth outside $\vrt{\mdl_X}$. 
Using the pairing (\ref{FA}), 
the Skeleton Theorem implies that the Abel-Jacobi map with modulus 
induces an isomorphism between affine parts 
\[ \abjacmaff{X}{\mdl}: \ACHm{X}{\mdl} \iso \Lm{X}{\mdl}(\fld) 
\] 
($\see$Theorem \ref{affineRoitmanThm}, ``Affine Roitman Theorem''). 
With this tool we prove: 
Let $X$ be strongly suitable 
($\see$Definition \ref{strongSuitable}, e.g.\ smooth)
projective variety, 
assume that the base field $\fld$ is an algebraic closure of a finite field 
or algebraically closed of characteristic $0$: 
Then the Abel-Jacobi map with modulus 
is an isomorphism on torsion parts: 
\[ \CHm{X}{\mdl}^{\tor} \iso \Albm{X}{\mdl}(\fld)^{\tor} 
\] 
($\see$Theorem \ref{RoitmanThm}, ``Roitman Theorem with Modulus''). 

According to the unramified class field theory of Kato and Saito \cite{KS83} 
and a statement about Albanese kernels with modulus (Lemma \ref{ajm-exSeq}) 
we obtain: 
Let $X$ be a smooth projective geometrically irreducible variety 
over a finite field $\fld$. 
Then we have an exact sequence 
\[ 0 \lra \Kaj{X} \lra \CHmo{X}{\mdl} \xra{\abjacm{X}{\mdl}} \Albm{X}{\mdl}\pn{\fld} 
\lra 0 
\] 
($\see$Corollary \ref{AlbKernel-finFld}), 
where $\Kaj{X}$ is 
the Pontrjagin dual of a subgroup 
of the torsion part of the N\'eron-Severi group of $X$ 
(defined in (\ref{SX})). 
If the N\'eron-Severi group of $X$ is torsion free, 
or if $\fld$ is an algebraic closure of a finite field, 
then the Abel-Jacobi map with modulus is an isomorphism: 
\[ \abjacm{X}{\mdl}: \CHmo{X}{\mdl} \iso \Albm{X}{\mdl}(\fld) 
\]  
($\see$Corollary \ref{AlbKernel-finFld} 
resp.\ Theorem \ref{RoitmanThm over finite field}).

\bigskip 
\textbf{Section~\ref{ReciprocityMap}:} Geometric Class Field Theory. \\ 
Let $\fld = \bF_q$ be a finite field 
and $X$ be a suitable projective 
variety over $\fld$. 
Let $\mdl$ be a modulus for $X$ 
and assume that $X$ is smooth outside $\vrt{\mdl_X}$. 

Let $\morphism: Y \dra X$ be an abelian finite \'etale covering 
over $X \setminus \mdl_X$ of suitable projective varieties over $\fld$. 
We say \,$Y \dra X$\, is of \emph{ramification bounded by $\mdl$}, 
if for every curve $C$ in $X$ meeting $\mdl_X$ properly, 
the divisor $\mdl_{\Ct} = \sum_{\pnt \in \vrt{\mdl \cut \Ct}} n_{C,\pnt} \bt{\pnt}$ 
encodes an upper bound 
for the upper indices of the higher ramification groups of 
$\Gal\bigpn{\wh{K}_{C,\pnt}^{\ab} \big| \wh{K}_{C,\pnt}}$ 
whose images do not vanish in $\Gal\pn{K_{Y}|K_{X}}$ 
($\see$Definition \ref{bounded ramification}). 
The \emph{abelian fundamental group $\fundGabm{X}{\mdl}$ 
of $X$ of modulus $\mdl$} 
classifies abelian coverings of $X$ of ramification bounded by $\mdl$ 
($\see$Definition \ref{fundGroupMod}), 
and the geometric part $\fundGgeom{X}{\mdl}$ classifies 
those that arise from a ``geometric situation'', 
i.e.\ not from extending the base field. 

By Lang's class field theory of function fields over finite fields, 
for any abelian covering $Y \dra X$ 
of smooth projective geometrically connected varieties over $\fld$ 
there is an effective divisor $\mdl$ for $X$ 
yielding a bound for its ramification 
and a canonical isomorphism 
\[ 
\frac{\Albm{X}{\mdl}\pn{\fld}}{\morphism_* \Albm{Y}{\mdl}\pn{\fld}} 
\;\iso\; \Gal\pn{K_Y|K_X} 
\] 
($\see$Theorem \ref{reciprocityLaw}). 
In particular, if $X_{\mdl} \dra X$ denotes the covering of $X$ 
that arises as a pull-back of 
\,$\wp = \Frob_q - \id: \Albm{X}{\mdl} \lra \Albm{X}{\mdl}$, 
the ``$q$-power Frobenius minus identity'', 
we have a canonical isomorphism 
\[ \Albm{X}{\mdl}\pn{\fld} \iso \Gal\pn{K_{X_{\mdl}}|K_X} 
\] 
($\see$Corollary \ref{Alb=Gal}, ``Existence Theorem''). 
Taking the limit over all effective divisors $\mdl$ on $X$ 
we obtain a canonical isomorphism 
\[ \varprojlim_{\mdl} \Albm{X}{\mdl}\pn{\fld} 
\iso \Gal\pn{K^{\ab}_{X} | K_{X} \clfld} 
\] 
($\see$Corollary \ref{Gal_maxAb}), 
where $\Gal\pn{K^{\ab}_{\Xo} | K_{X} \clfld}$ is the geometric Galois group 
of the maximal abelian extension $K^{\ab}_{X}$ of the function field $K_{X}$ 
of $X$. 

Using the natural relations between the groups $\CHmo{X}{\mdl}$, 
$\fundGgeom{X}{\mdl}$ and $\Albm{X}{\mdl}$ 
($\see$Propositions \ref{CHm-to-fundGm} and \ref{fundGm-sur-Albm}), 
the affine Roitman Theorem \ref{affineRoitmanThm} 
and the unramified class field theory of Wiesend \cite{KeSc09} 
yield the following canonical isomorphism of finite groups: 
\[ \recmo{X}{\mdl}: \CHmo{X}{\mdl} \iso \fundGgeom{X}{\mdl} 
\] 
($\see$Theorem \ref{reciprocityDiagram}, ``Reciprocity Law''). 
Let $S$ be the support of an effective divisor on $X$ 
such that $X \setminus S$ is smooth. 
Taking the limit over all effective divisors $\mdl$ on $X$ 
with support in $S$ 
we obtain 
\[ \varprojlim_{\substack{\mdl \\ \vrt{\mdl} \subset S}} \CHmo{X}{\mdl} 
\;\iso \varprojlim_{\substack{\mdl \\ \vrt{\mdl} \subset S}} \fundGgeom{X}{\mdl} 
\;\cong\; \fundGgeo{X \setminus S} 
\] 
(Corollary \ref{reciprocityLimit}). 

\begin{flushright} 
e-mail: \texttt{henrik.russell@gmail.com} 
\end{flushright}

\end{document}